\documentclass[12pt,reqno,twoside,fleqn,twoside,a4paper]{amsart}%
\usepackage{amscd,amsmath,amssymb,amsthm,palatino}
\usepackage[mathcal]{euler}

\numberwithin{equation}{section}
\makeatletter
\renewcommand{\section}{\@startsection {section}{1}{\z@}%
                                   {-3.5ex \@plus -1ex \@minus -.2ex}%
                                   {.5\linespacing}%
                                   {\normalfont\scshape\centering}}
\makeatother

\newtheorem{thm}{Theorem}[section]

\newtheorem{lem}[thm]{Lemma}
\newtheorem{cor}[thm]{Corollary}
\newtheorem{prop}[thm]{Proposition}
  
\theoremstyle{definition}
\newtheorem*{definition}{Definition}

\def\beq#1\eeq{\begin{equation}#1\end{equation}}
\makeatletter
\@ifpackageloaded{euler}{
 \newcommand{\fchoice}[2]{#1}
 \newcommand{\onto}{\to\mkern-14mu\to}
 \def\rightarrowfill@#1{\m@th\setboxz@h{$#1\relbar$}\ht\z@\z@
   $#1\copy\z@\mkern-6mu\cleaders
   \hbox{$#1\mkern-2mu\box\z@\mkern-2mu$}\hfill
   \mkern-6mu\mathord\rightarrow$}
 \def\leftarrowfill@#1{\m@th\setboxz@h{$#1\relbar$}\ht\z@\z@
   $#1\mathord\leftarrow\mkern-6mu\cleaders
   \hbox{$#1\mkern-2mu\copy\z@\mkern-2mu$}\hfill
   \mkern-6mu\box\z@$}
 \def\B@R#1#2{\raisebox{-.07ex}{$#1#2$}\mkern-6mu}
 \renewcommand{\hbar}{{\mspace{1mu}\mathpalette\B@R{\mathchar'26}h}}
}{
 \newcommand{\fchoice}[2]{#2}
 \DeclareMathSymbol{\onto}{\mathrel}{AMSa}{"10}
 \renewcommand{\hbar}{{\mathchar'26\mkern-9muh}}
}
\makeatother
\newcommand{\A}{\mathcal{A}}
\renewcommand{\AA}{\mathbb{A}}
\newcommand{\As}{\A_s}
\newcommand{\C}{\mathcal{C}}
\newcommand{\M}{\mathcal{M}}
\newcommand{\co}{\mathbb{C}}
\newcommand{\I}{\mathcal{I}}
\newcommand{\N}{\mathbb{N}}
\newcommand{\cs}{\mbox{\upshape C}\ensuremath{{}^*}} 
\newcommand{\csr}{{\mathrm C}^*_{\mathrm r}}
\newcommand{\csm}{\cs_{\mathrm{max}}}
\newcommand{\R}{\mathbb{R}}
\newcommand{\Z}{\mathbb{Z}} 
\newcommand{\into}{\hookrightarrow}
\newcommand{\isom}{\mathrel{\widetilde\longrightarrow}}
\newcommand{\Po}{\mathcal{P}}
\newcommand{\Gh}{\Gamma_{\!\mathrm{hol}}}
\newcommand{\Gc}{\Gamma_{\!\mathrm{c}}}
\newcommand{\Gz}{\Gamma_{\! 0}}
\newcommand{\Gi}{\Gamma^\infty}
\newcommand{\Gci}{\Gamma_{\mathrm c}^\infty}
\newcommand{\Gzi}{\Gamma_0^\infty}
\newcommand{\Hi}{\mathcal H}
\newcommand{\jb}{{\bar{\jmath}}} 
\newcommand{\ib}{{\bar{\imath}}}  
\DeclareMathOperator{\End}{End}

\DeclareMathOperator{\Image}{Im}
 \renewcommand{\Im}{\Image}
\DeclareMathOperator{\id}{id}

\DeclareMathOperator{\td}{td}

\DeclareMathOperator{\ind}{ind}
\DeclareMathOperator{\Spec}{Spec}
\newcommand{\Kahler}{K\"ahler} 
\newcommand{\norm}[1]{\lVert#1\rVert}
\newcommand{\Norm}[1]{\left\|#1\right\|}
\newcommand{\abs}[1]{\lvert#1\rvert}
\newcommand{\Abs}[1]{\left|#1\right|}
\newcommand{\intM}{\int_\M\mathchoice{\mskip-2\thinmuskip}%
                   {\mskip-\thinmuskip}{}{}}
\newcommand{\TM}{T\fchoice{}{\!}\M}
\newcommand{\Mt}{{\widetilde\M}}
\newcommand{\G}{\Gamma}
\newcommand{\g}{\gamma}
\renewcommand{\ss}{{\sigma_s}}
\newcommand{\h}{\mathfrak H}
\newcommand{\hs}{\h_s}
\newcommand{\Hs}{\Hi_s}
\newcommand{\Li}{\mathcal L}
\newcommand{\Lb}{\mathfrak L}
\newcommand{\Ls}{\Lb^s}
\newcommand{\LL}{\Lb}
\newcommand{\K}{\mathcal K}
\newcommand{\Khat}{\hat K}
\newcommand{\B}{\mathcal{B}}
\newcommand{\Bs}{{\B_s}}
\newcommand{\Bsm}{{\B_s^{\max}}}
\newcommand{\Hsm}{\Hs^{\max}}
\newcommand{\Lsm}{\Ls_{\max}}
\newcommand{\Asm}{\As^{\max}}

\newcommand{\Miscenko}{Mi\v{s}\v{c}enko}
\DeclareMathOperator{\dom}{dom}
\newcommand{\braket}[3]{\left<#1\vphantom{#2}\right|#2\left|#3\vphantom{#2}\right>}
\newcommand{\db}{\nabla_{\!\bar\partial}}
\newcommand{\dbb}{\nabla_{\!\partial}}
\newcommand{\Cb}{\C_{\mathrm b}}
\newcommand{\Gb}{\Gamma_{\!\mathrm b}}
\newcommand{\Or}{\mathcal O}
\newcommand{\dv}{d^{\mathrm v}}
\renewcommand{\dh}{d^{\mathrm h}}
\DeclareMathOperator{\tot}{tot}
\newcommand{\one}{{}^{\scriptscriptstyle\mathrm I}}
\newcommand{\two}{{}^{\scriptscriptstyle\mathrm{II}}}
\newcommand{\Qb}{\mathsf{Q}}
\newcommand{\T}{\mathbb{T}}

\newcommand{\algtensor}{\otimes_{\mathrm{Alg}}}
\newcommand{\bimod}{\ensuremath{\mathsf{bimod}}}
\newcommand{\Thol}{T_{\mathrm{hol}}}

\newcommand{\HH}{H}

\newcommand{\Ebar}{\underline{E}}
\newcommand{\kbar}{\underline{k}}
\newcommand{\inner}{\mathbin{\raise1.5pt\hbox{$\lrcorner$}}}
\DeclareMathOperator{\tr}{tr}
\newcommand{\ntr}{\widetilde{\smash\tr\vphantom{\raise.5pt\hbox{r}}}}
\DeclareMathOperator{\vol}{vol}
\providecommand{\href}[2]{#2}

\title{Quantization of Multiply Connected Manifolds}
\date{}
\author{Eli Hawkins}

\begin{document}
\begin{flushright}
\vspace*{-0.5in}
\begin{tabular}{l}
\textsf{\small SISSA 33/2003/FM}\\
\textsf{\small math.QA/0304246}\\
\end{tabular}
\vspace{0.25in}
\end{flushright}
\vfil
\maketitle
\begin{center}
\vspace{-4ex}
\emph{\small Scuola Internazionale Superiore di Studi Avanzati }\\
\emph{\small Via Beirut 4, I-34014 Trieste, Italy}\\
{\small mrmuon@mac.com}\\
\end{center}
\begin{abstract}
The standard (Berezin-Toeplitz) geometric quantization of a compact \Kahler\ manifold is restricted by integrality conditions. These restrictions can be circumvented by passing to the universal covering space, provided that the lift of the symplectic form is exact.  I relate this construction to the Baum-Connes assembly map and prove that it gives a strict quantization of the manifold. I also propose a further generalization, classify the required structure, and provide a means of computing the resulting algebras. These constructions involve twisted group \cs-algebras of the fundamental group which are determined by a group cocycle constructed from the cohomology class of the symplectic form.
\end{abstract}
\vfil\vfil
\thispagestyle{empty}
\newpage

\tableofcontents
\section{Introduction}
As the name suggests, ``geometric quantization'' was originally 
conceived as a means of constructing models of quantum physics. The 
idea was to construct quantum observables from the classical phase 
space. Mathematically, the phase space is a symplectic manifold and 
the quantum observables are operators on a Hilbert space. 

It has not thus far been altogether successful for this original purpose. 
The truly challenging problems of quantization remain unsolved. 
Geometric quantization has not elucidated quantum field theory and 
it has not been successfully applied to General Relativity to yield a 
viable theory of quantum gravity. Instead, geometric quantization has 
had a more fruitful mathematical life in representation theory 
because many of the Hilbert spaces that it constructs naturally carry 
unitary group representations. Much of the literature on geometric 
quantization is thus entirely unconcerned with algebras of ``observable'' 
operators.

More recently, a new use for geometric quantization has arisen in the 
context of Noncommutative Geometry. In Noncommutative 
Geometry, geometric concepts are reexpressed in an algebraic form. In 
particular, the topology of a space is expressed through an algebra 
of continuous functions on that space. This is done in such a way that 
geometry can be generalized by allowing this algebra ``of functions'' 
to be noncommutative. 

Geometric quantization comes into this as a potential means of 
constructing noncommutative spaces (but see \cite{haw6}). On a given manifold, the geometric
quantization construction depends upon a parameter $\hbar$ 
(Planck's constant). As $\hbar\to0$, classical physics was supposed to 
be recovered. Quantum observable operators are constructed from 
functions on the classical phase space; as $\hbar\to 0$, the 
composition of operators converges to the multiplication of 
functions. Because of this, algebras of quantum observable operators 
approximate (in an appropriate sense) the algebra of functions on the 
classical phase space. From the perspective of Noncommutative 
Geometry, \cs-algebras are noncommutative topological spaces; thus
what geometric quantization does is to produce noncommutative topologies 
which approximate the topology of the given symplectic manifold.
\medskip

When it comes to quantization, you can't always get what you want, or 
may expect. When quantization was first approached 
mathematically, it was assumed that the commutator of  
operators, $[a,b]_-:=ab-ba$, should be related to the Poisson bracket of
functions by, 
\beq
[Q(f),Q(g)]_- = -i\hbar Q(\{f,g\})
\mbox,\label{commutator}
\eeq
where $Q$ is the map that takes functions to corresponding operators. This 
assumption \eqref{commutator} was ill-founded. The Groenwald--van\,Hove ``no go'' theorem (see \cite{che2,g-s1}) showed that this  equation cannot be satisfied exactly for all smooth functions. Now, this equation is only assumed to hold to leading order 
as $\hbar\to 0$.

Similarly, it was initially assumed that quantization ought to work for 
all values of $\hbar$ in some interval containing $0$, see 
\cite{rie1}. This too was overoptimistic. Geometric quantization 
constructions involve an $\hbar$-dependent family of line bundles 
whose curvatures diverge as $\hbar\to0$. However, the curvature of a 
line bundle must belong to an integral cohomology class. For a 
compact manifold, this is very restrictive; such line 
bundles \emph{cannot} vary continuously. 
\medskip

Seemingly in contradiction to this is the ``noncommutative
torus''\/,  probably the most frequently cited example of a noncommutative
space.  
\begin{definition}
For $\theta\in\R$, the noncommutative torus algebra $\T_\theta$ is the \cs-algebra generated by two unitary operators $U$ and $V$ with the single relation, 
\beq
VU=e^{-2\pi i\theta}UV
\mbox.\label{nct}\eeq
\end{definition}
Note that this only depends on $\theta$ modulo $1$. The parameter $\theta$ plays the role of $\hbar$. As $\theta\to0$, this becomes the algebra of functions on an ordinary torus. This construction is not limited to discrete values of $\theta$, it works for \emph{any} value of $\theta\in\R$.
Geometric quantization looks on in jealousy! Can the construction of 
the noncommutative torus be generalized?

The torus is a quotient of the plane, $\R^2$, by an action of the 
group $\Z^2$. Likewise, the noncommutative torus is, in a sense, a quotient of a 
noncommutative plane by $\Z^2$. Formally, the algebra of functions on 
the plane can be generated by two commuting unbounded operators, the 
coordinates $x$ and $y$. We can ``quantize'' the plane by replacing 
commutativity with the relation,
\beq
[x,y]_-= -i\theta/2\pi
\mbox.\label{NCplane}\eeq
Now, working in terms of algebras, the quotient by $\Z^2$ is constructed by restricting to $\Z^2$-invariant ``functions''\/. In other words, restrict attention to period $1$ functions of $x$ and $y$. Any such periodic function is a function of $e^{2\pi i x}$ and $e^{2\pi i y}$. Thus, the algebra of ``functions'' on the quotient of this noncommutative plane by $\Z^2$ is generated by the two unitary operators $U=e^{2\pi i x}$ and $V=e^{2\pi iy}$. These satisfy eq.~\eqref{nct}. 

The obstruction to quantization for all values of $\hbar$ is in degree $2$ cohomology. However, for the plane $H^2(\R^2)=0$. This is why there is no obstruction to quantizing the plane for arbitrary $\hbar$. The construction of the noncommutative torus was generalized by Klimek and Lesniewski \cite{k-l2,k-l4} to higher-genus Riemann surfaces. These have the hyperbolic plane as their universal covering space. They first quantized the hyperbolic plane and then ``quotiented'' by the appropriate fundamental groups to produce quantum Riemann surfaces.

Natsume and Nest took an entirely different approach to ``quantum surfaces''\/. They started from the observation that the algebras, $\T_\theta$, are twisted group \cs-algebras of the group $\Z^2$. They proposed that quantum Riemann surfaces should be at least Morita equivalent to twisted group \cs-algebras of the fundamental group of the surface in question. Based on this principle, they constructed algebras for quantum Riemann surfaces. 

My main results in this paper will be to elucidate the relationship between the quotient construction and the twisted group \cs-algebras of the fundamental group, and to show that the construction gives strict deformation quantizations for all compact \Kahler\ manifolds to which it applies.

\subsection{Notation}
Before proceeding further, let's pause and review some of the relatively standard notations that I will be using.

If $\M$ is a smooth manifold, then $\Gamma^k(\M,E)$ denotes the space of $k$-times continuously differentiable sections of a vector bundle $E$. The tangent bundle is $\TM$. The $p$-form bundle is $\wedge^p\M := \wedge^p T^*\M$, the $p$'th exterior power of the cotangent bundle. $\Omega^p(\M) := \Gamma^\infty(\M,\wedge^p\M)$ is the space of smooth $p$-forms, and $\Omega^p(\M,E) := \Gamma^p(\M,E\otimes\wedge^p\M)$ is the space of smooth , $E$-valued $p$-forms.

If $\M$ happens to be a complex manifold, then the holomorphic tangent bundle, $\Thol\M$, is the tangent bundle in the sense of complex manifolds. The holomorphic $p$-form bundle is $\wedge^{p,0}\M := \wedge^p\Thol^*\M$. The antiholomorphic $p$-form bundle is the complex conjugate bundle, $\wedge^{0,p}\M := \overline{\wedge^{p,0}\M}$. The $(p,q)$-form bundle is $\wedge^{p,q}\M := \wedge^{p,0}\M\otimes\wedge^{0,q}\M$, and $\Omega^{p,q}(\M) := \Gamma^\infty(\M,\wedge^{p,q}\M)$ is the space of smooth $(p,q)$-forms. Finally, $\Omega^{p,q}(\M,E) := \Gamma^\infty(\M, E\otimes\wedge^{p,q}\M)$ is the space of smooth $E$-valued $(p,q)$-forms.

I use the convention that the curvature of a connection, $\nabla$, is $i\nabla^2$. When referring to a bundle, I usually mean a bundle with a specific connection.

All cohomology is taken with complex coefficients, unless otherwise specified. Square brackets will be used to denote a cohomology class (or $K$-theory or $K$-homology class). For instance, $[\omega] \in H^2(\M)$ is the cohomology class of the symplectic form $\omega\in\Omega^2(\M)$. Square brackets will also be used to denote the basis elements of a group algebra. In Section \ref{injectivity}, I follow algebraic geometry notation and use square brackets for the line bundle associated to a divisor. Yet another use for square brackets is the commutator-anticommutator notation:
\[
[a,b]_\pm := ab \pm ba
\mbox.\]

$\T$ is the (multiplicative) group of complex numbers of modulus $1$. The identity element of a group in multiplicative notation is usually denoted $e$.

If $A$ is a subspace of an algebra and $B$ is a subspace of a module, then $AB$ denotes the linear span of all products of elements of $A$ with elements of $B$.

A tensor product should always be understood to be taken in the most restrictive relevant category --- vector spaces, algebras, Hilbert spaces, Hilbert \cs-modules, etc. A tensor product of two \cs-algebras is always the reduced \cs-algebraic tensor product here. A homomorphism between \cs-algebras is always assumed to be a $*$-homomorphism. % All \cs-algebras here are assumed to be separable.

If $H$ is a Hilbert space, then $\Li(H)$ and $\K(H)$ denote the algebras of bounded and compact operators, respectively.

I use Dirac's bra-ket notation for inner products. The inner product in a Hilbert space or Hilbert \cs-module is denoted as $\langle\psi\vert\varphi\rangle$. This is antilinear in the \emph{first} argument and linear in the second. For inner products involving operators, $\braket\psi{a}\varphi$ is the same as $\langle\psi\vert a\varphi\rangle$. A vector is also sometimes written as $\lvert\varphi\rangle$. The notation is associative; $\lvert\psi\rangle\langle\varphi\rvert$ denotes the operator such that,
\[
\left(\lvert\psi\rangle\langle\varphi\rvert\right)\lvert\chi\rangle 
= \lvert\psi\rangle\langle\varphi\vert\chi\rangle
\mbox.\]

If $X\subset \M$ is a subset of a smooth manifold, but not necessarily a manifold itself, then I regard a function on $X$ to be smooth (infinitely differentiable) if it is the restriction of a smooth function over $\M$. For example, if $X$ is an open subset, the space $\C_0^\infty(X) := \C_0(X) \cap \C^\infty(X) $ is the subspace of smooth functions in $\C_0^\infty(\M)$ which vanish identically on $\M\smallsetminus X$.

\subsection{Quantization Defined}
\label{quantization.definition}
The intuitive notion that an algebra is continuously deformed takes a mathematically precise statement with the notion of a continuous field of \cs-algebras, \cite{dix1}. This is essentially a bundle of \cs-algebras, but with the condition of local triviality considerably relaxed.

This definition is the one given by Kirchberg and Wasserman in \cite{k-w1}. It is equivalent to that given by Dixmier \cite{dix1}, provided that the base space, $X$, is locally compact and Hausdorff.
\begin{definition} 
\label{field.def}
A \emph{continuous field of \cs-algebras,} $\A$, over a locally compact Hausdorff space, $X$, is given by a set of \cs-algebras $\{\A_x\}_{x\in X}$, a \cs-algebra $\Gz(X,\A)$, and a set of surjective homomorphisms, $\{\Po_x:\Gz(X,\A) \onto \A_x \}_{x\in X}$, such that:
\begin{enumerate}
\item
$\forall a\in \Gz(X,\A)$, $\norm a = \sup_{x\in X} \Norm{\Po_x(a)}$
\item
$\Gz(X,\A)$ is a $\C_0(X)$-module such that for $f\in\C_0(X)$, $a\in\Gz(X,\A)$, and $x\in X$, $\Po_x(fa) = f(x) \Po_x(a)$.
\item
For $a\in \Gz(X,\A)$, $x \mapsto \Norm{\Po_x(a)}$ defines a function in $\C_0(X)$.
\end{enumerate}
A \emph{continuous field of Banach spaces} over $X$ is defined identically, but with ``Banach space'' replacing ``\cs-algebra'' and ``contractive map'' replacing ``$*$-hom\-o\-morph\-ism''\/.
\end{definition}

The direct product of the evaluation maps, $\Po_x$, identifies $\Gz(X,\A)$ with a subspace of the direct product $\prod_{x\in X}\A_x$. A section (element of the direct product) is \emph{continuous} if, over any compact subset of $X$, it coincides with some element of $\Gz(X,\A)$. As the notation suggests, $\Gz(X,\A)$ is the closure of the space of compactly supported continuous sections. For $X$ compact the subscript $0$ is superfluous. From this, we can also define the space, $\Gb(X,\A)$, of bounded continuous sections. This is also a \cs-algebra (or respectively, Banach space). 

An upper (lower) semicontinuous field of \cs-algebras or Banach spaces is defined identically, but for the last condition. The functions defined by norms are only required to be upper (lower) semicontinuous, rather than continuous.
\medskip

Quantization of a manifold, $\M$, should involve a continuous field of \cs-algebras over the set $\I\subseteq\R$ of allowed values of $\hbar$, such that the fiber over $\hbar=0$ is $\C_0(\M)$.

The quantization map, $Q$, should take functions on $\M$ to sections 
of the continuous field, in such a way that the value of the section 
$Q(f)$ at $\hbar=0$ is $f$ itself. This implicitly encompasses the 
condition that the quantization map should asymptotically preserve 
multiplication in the limit as $\hbar\to0$.

The commutator $[Q(f),Q(g)]_-$ for two smooth functions $f,g\in\C^\infty_0(\M)$ should be given by eq.~\eqref{commutator} to leading order as $\hbar\to\infty$. For a symplectic manifold the Poisson bracket is determined by the symplectic structure. In general it is a Lie bracket on $\C^\infty_0(\M)$ which is a derivation on both arguments.

This leads to the following definition of a strict quantization,  essentially equivalent to that given Landsman in \cite{lan2}.
\begin{definition} 
\label{quantization.def}
A \emph{strict quantization} of a smooth manifold, $\M$, with a Poisson bracket 
consists of a continuous field of \cs-algebras, $\A$, over some closed subset
$\I\subseteq\R$ which is dense at $0$, and a continuous linear map $Q:\C_0(\M)\to\Gb(\I,\A)$. These must satisfy the following conditions, 
\begin{enumerate}
     \item For each $\hbar\in\I$, the image of $\Po_\hbar \circ Q : \C_0(\M) \to \A_\hbar$ generates the algebra $\A_\hbar$ as a \cs-algebra.  
  \item For the evaluation $\Po_0:\Gb(\I,\A)\to\A_0$ at $\hbar=0\in\I$, the composition $\Po_0\circ Q$ is an isomorphism, $\Po_0\circ Q : \C(\M) \isom \A_0$.
    \item For any $f,g\in\C^{\infty}_0(\M)$,
    \[
    [Q(f),Q(g)]_- -i\hbar Q\{f,g\}\in \hbar \ker\Po_0
    \mbox.\]
\end{enumerate}
\end{definition}
This is equivalent to Landsman's definition except for the first condition. He makes the more restrictive assumption that $\Po_\hbar\circ Q$ is surjective, although he notes that this condition may be excessive.

A great deal of structure is compressed into the quantization map, $Q$. There is really a family of quantization maps, $\Po_\hbar\circ Q : \C_0(\M)\to\A_\hbar$ for $\hbar\in\I$. The algebras $\A_\hbar$ are determined by these maps. Even the continuous field structure is uniquely determined by the condition that $Q$ gives continuous sections (see \cite{lan2}). 

\subsection{Standard Construction}
\label{classical}
Let $\M$ be a manifold with symplectic form $\omega$. The symplectic structure is a closed 2-form, $\omega\in \Omega^2(\M)$, which is invertible if viewed as a map from vectors to 1-forms. A symplectic manifold always has a Poisson structure. For any function $f\in\C^1(\M)$, an associated vector-field, $\xi_f \in \G(\M,\TM)$, is uniquely defined by $\xi_f \inner \omega = df$. The Poisson bracket 
of two functions $f,g\in\C^1(\M)$ is defined by,
\beq\label{Poisson}
\{f,g\} := \xi_f(g) = -\xi_g(f) = -\omega(\xi_f,\xi_g)
\mbox.\eeq

All standard geometric quantization constructions involve a line bundle $L$ over $\M$ with a connection and inner product, and curvature $\omega$. This is often called a ``quantum'' line bundle, but ``quantization line bundle'' seems a more appropriate term. The existence of such a line bundle implies the integrality condition that $[\omega]\in2\pi\, H^2(\M,\Z)$. The tensor power line bundle $L^{\otimes N}$ has curvature $N\omega$. The inner product on fibers of $L^{\otimes N}$ can be integrated against the symplectic volume form, $\omega^n/n!$, ($2n=\dim\M$) to give an inner product on sections. Completion of the space of sections of $L^{\otimes N}$ with respect to the associated norm gives the Hilbert space $L^2(\M,L^{\otimes N})$. The principal Hilbert spaces employed in geometric quantization are subspaces of each $L^2(\M,L^{\otimes N})$ defined by some ``polarization'' condition.

I shall only be concerned with \Kahler\ polarization, so suppose that $\M$ is a \Kahler\ manifold. That is, $\M$ is simultaneously and compatibly a Riemannian, complex, and symplectic manifold. The polarization condition is holomorphicity. The Hil\-bert space that we are interested in is thus the space of square-integrable holomorphic sections, $L^2_{\mathrm{hol}}(\M,L^{\otimes N})$. Let $\Pi_N$ be the orthogonal projection from $L^2(\M,L^{\otimes N})$ onto the subspace $L^2_{\mathrm{hol}}(\M,L^{\otimes })$. The simplest quantization map in this setting is the so-called Toeplitz quantization map, $T_N: L^\infty(\M) \to \Li[L^2_{\mathrm{hol}}(\M,L^{\otimes N})]$ defined concisely by,
\[
T_N(f) := \Pi_N f
\mbox.\]
That is, for any essentially bounded function $f$, the action of the operator $T_N(f)$ on $L^2_{\mathrm{hol}}(\M,L^{\otimes N})$ is defined by first multiplying a holomorphic section $\psi \in L^2_{\mathrm{hol}}(\M,L^{\otimes N})$ by $f$ (to give a non-holomorphic section $f\psi \in L^2(\M,L^{\otimes N})$) and then projecting back orthogonally to the subspace of holomorphic sections. Although $T_N$ is defined on $L^\infty(\M)$, we are really only interested in applying it to continuous functions $f\in \C(\M)$. Despite the name, the Toeplitz quantization construction for a \Kahler\ manifold is due to Berezin \cite{ber1}.

If $\M$ is compact and connected, then the Hilbert spaces of holomorphic sections are each finite-di\-men\-sion\-al, and each $T_N$ is a surjection onto the full matrix algebra on this space (see \cite{b-m-s}).  Along with $\C(\M)$, these matrix algebras can be assembled into a continuous field of \cs-algebras over $\hat\N := \{1,2,\ldots,\infty\}$ (the one-point compactification of the natural numbers). The continuous field structure is defined by the condition that for any $f\in\C(\M)$, the section given by $N\mapsto T_N(f)$ and $\infty\mapsto f$ is continuous. 

This defines a quantization of $\M$; see, \cite{b-m-s}. We must identify $N$ with $\hbar^{-1}$, so that $\hat\N$ is identified homeomorphically with $\{0,N^{-1} \mid N\in\N\} \subset \R$.

This can be generalized very slightly by replacing $L^{\otimes N}$ with $L_N:=L^{\otimes N}\otimes L_0$ for some holomorphic bundle, $L_0$. In the standard construction, $L_0$ is trivial, but I prefer not to make that assumption and will base the constructions in this paper on the more general form. We can think of $N\mapsto L_N$ as defining a straight line through the lattice of possible holomorphic line bundles. Restricting to $L_0$ trivial would be like only considering lines which pass through the origin.

 It is most natural to view $L_0$ as having a $\co$-valued inner product on its \emph{sections} (a pre-Hilbert space structure). Such an inner product can also be viewed as a vector bundle homomorphism ${\bar L}_0\otimes L_0 \to \wedge^{2n}\M$. In the standard case, $L_0=\co\times\M$ and the inner product is the obvious one multiplied by the canonical volume form, but there are other choices which may be important. For instance, if we choose $L_0=\wedge^{n,0}\M$, then because $\wedge^{2n}\M=\wedge^{0,n}\M\otimes\wedge^{n,0}\M$, the inner product is intrinsically defined independently of a metric.

The parameter space $\hat\N$ is certainly not homeomorphic to an interval; it 
is almost discrete. This is inevitable because the Hilbert spaces and algebras are finite-dimensional. The dimension is an integer and cannot vary continuously.

\subsection{Prelude}
So, integrality conditions for the curvature of a line bundle obstruct us from a continuously parametrized quantization of a compact manifold, $\M$. On the other hand, over a noncompact manifold it may be possible to define ``$L^{\otimes s}$'' for arbitrary, noninteger $s$.  Such a bundle should have curvature $s\omega$. If these ``tensor power'' bundles exist for all $s$, then the cohomology class of $s\omega$ must be integral for all $s$. This means that $\omega$ must be cohomologically trivial --- an exact $2$-form. If $\omega$ is exact, then we can in fact construct a line bundle with the desired curvature, $s\omega$, for any $s\in\R$.

This cannot happen on a compact symplectic manifold.  The invertibility of a symplectic form implies that $\omega^n := \omega\wedge\dots\wedge\omega \in \Omega^{2n}(\M)$ is nowhere-vanishing.  Thus, for $\M$ compact $\intM\omega^n \neq0$, and so $[\omega]\neq 0\in H^2(\M)$.

If the universal covering space, $\Mt$, is noncompact, then the symplectic form of $\Mt$ might be exact. If so, it will be possible to construct quantization line bundles over $\Mt$ for arbitrary $s$. Let $\pi:\Mt\onto\M$ be the universal covering; the symplectic form of $\Mt$ is the pullback $\widetilde\omega := \pi^*\omega$. The necessary condition can be stated in several ways. First, $\widetilde\omega$ is exact. This means that $[\omega] \in H^2(\M)$ is contained in the kernel of the functorial map $\pi^*: H^2(\M) \to H^2(\Mt)$. Equivalently, $[\omega]$ pairs trivially with the image of $\pi_*:H_2(\Mt) \to H_2(\M)$, but by the Whitehead theorem $H_2(\Mt,\Z)=\pi_2(\Mt)=\pi_2(\M)$, so $[\omega]$ must pair trivially with $\pi_2(\M)$. This gives the most directly geometrical restatement: $\widetilde\omega$ is exact if and only if $\int_{S^2}\omega =0$ for any embedded $2$-sphere, $S^2\subset\M$.

If this condition is satisfied, then we can construct quantization line bundles and Toeplitz maps for $\Mt$. The space of continuous functions on $\M$ can be thought of as a subalgebra $\C(\M)\subset\C_{\mathrm b}(\Mt)$ of the bounded continuous functions on $\Mt$; it is precisely the subalgebra of functions which are invariant under the action of the fundamental group, $\G:=\pi_1(\M)$, on $\Mt$. The key idea here, following the noncommutative torus and Klimek--Lesniewski, is that the Toeplitz quantization maps for $\Mt$ can be applied to $\C(\M)$. I shall prove in this paper that this ``covering construction'' yields a quantization of $\M$.

For a compact, connected symplectic manifold, the \cs-algebra generated by the image of a standard Toeplitz map is easily characterized. As I have indicated, it is just a full matrix algebra so it is specified by its dimension. For the generalized version which I am considering here, characterization becomes more complicated. The image of an $\Mt$-Toeplitz map consists of bounded operators on a Hilbert space, and the algebra we are interested in must be $\G$-invariant. However, the algebra of all $\G$-invariant operators on a Hilbert space is a von\,Neumann algebra. Surely the algebra in question is something smaller than this! 

Indeed it is. A clearer picture of this algebra can be obtained by using a more refined structure than a Hilbert space. I shall show that the algebra generated by the image of the $\Mt$-Toeplitz map is isomorphic to the algebra of compact operators on a Hilbert \cs-module of a twisted group \cs-algebra of the fundamental group of $\M$.

My approach to analyzing the covering construction and proving that it gives a quantization will be based upon a different, but equivalent, construction. This is a special case of a generalization of the Toeplitz quantization construction. Rather than only considering quantization line bundles with curvature proportional to $\omega$, I define a \emph{quantization bundle} to be a bundle of Hilbert \cs-modules over $\M$ with curvature proportional to $\omega$. The construction of a Toeplitz map from a quantization bundle is formally identical to the construction of the Toeplitz map, $T_N$, from $L^{\otimes N}$. There exists a family of quantization bundles, $\Ls$, such that the Toeplitz quantization is equivalent to the above covering construction. 

The very existence of these quantization bundles shows that quantization bundles can circumvent integrality conditions. Also notable is the property of another related family of quantization bundles, $\Lsm$; any quantization bundle with curvature $s\omega$ can be constructed from $\Lsm$.

This also draws a connection with the quantum surfaces of Natsume and Nest. They propose quantum surface algebras which are Morita equivalent to reduced, twisted group \cs-algebras of the fundamental group, $\G$. The fibers of my quantization bundles are Hilbert \cs-modules of reduced, twisted group \cs-algebras of $\G$. A simple calculation using an index theorem suggests that my Toeplitz algebras for a Riemann surface are isomorphic to those of Natsume and Nest, at least for irrational parameter values.
\medskip

The quantization line bundles which can be constructed over $\Mt$ are not quite $\G$-equivariant. There is instead a projective action of $\G$ on sections. 
In Section \ref{twisting}, I review the definitions of projective representations and twisted group \cs-algebras. I construct the group cocycles which describe these projective actions of $\G$. I show that the reduced twisted group \cs-algebras constructed from such a family of cocycles form a continuous field of \cs-algebras.

In Section \ref{Hilbert.modules}, I review the definitions and basic properties of Hilbert \cs-mod\-ules and related constructions. I discuss the properties of bundles of Hilbert \cs-mod\-ules and Dirac operators on bundles of Hilbert \cs-modules. Finally, I construct the quantization bundles related to the covering space quantization construction.

In Section \ref{Toeplitz}, I present the Toeplitz construction for a quantization bundle over a compact manifold and analyze some of its properties. The primary tool in this analysis is a Dirac-type operator, the Dolbeault operator. When the curvature of a quantization bundle is sufficiently great, the space of holomorphic sections is the kernel of a Dolbeault operator. I show that this space is a Hilbert \cs-module and that the Toeplitz maps are asymptotically multiplicative.

In Section \ref{topological}. I briefly discuss the relationship of some of these constructions with topological constructions --- specifically, $KK$-theory, indices, Baum-Connes assembly maps, and $L^2$-index theorems.

Hilbert \cs-modules play a central role in the theory of Morita equivalence. In Section \ref{Morita}, I review two definitions of Morita equivalence and present a category in which Morita equivalence is isomorphism. I discuss some properties of my Toeplitz construction which are related to Morita equivalence. In particular, I give a method of computing the Toeplitz algebras.

In Section \ref{quantization}, I prove the equivalence between the covering  construction and the Toeplitz construction with the quantization bundle $\Ls$. I prove that these indeed give a quantization of $\M$.

The previous work to which this is most closely related was concerned with surfaces, so I discuss the $2$-dimensional examples in Section \ref{examples}. I prove that when my constructions are applied to a flat torus, the result is a standard noncommutative torus. I also give a partial result indicating that for a higher genus Riemann surface, my construction gives algebras isomorphic to those constructed by Natsume and Nest.

\section{The Twist}
\label{twisting}
Let $\M$ be smooth, connected manifold. I will eventually assume that $\M$ is a compact \Kahler\ manifold with symplectic form $\omega$, but in this section I only need $\omega\in \Omega^2(\M)$ to be a closed $2$-form.

Let $\pi:\Mt\onto\M$ be the universal covering and $\G:=\pi_1(\M)$ 
the fundamental group. Assume that the lift, $\widetilde\omega:=\pi^*\omega$, is an exact form.  (If $\M$ is compact, this implies that $\G$ is infinite.)
This means that there exists some $1$-form, $A\in\Omega^1(\Mt)$, such that $\widetilde\omega=dA$. Fix some choice of $A$, and define the bundle $L^s$ over $\Mt$ to be $\co\times\Mt$ with the connection 
\[
\nabla_s:=d-isA
\mbox.\] 
This has curvature $s\widetilde\omega$; that is $\nabla_s^2 = -is\widetilde\omega$. In fact, $L^s$ has the properties we would expect of an ``$s$'th'' tensor power of $L^1$, but $s$ is not restricted to be an integer.  

I will assume that $A$ is real. This choice simplifies the inner products. However, when $\M$ is \Kahler, it is possible to choose an equivalent $A\in\Omega^{1,0}(\Mt)$ which exhibits $L^s$ as a  holomorphic bundle.

Although $\G$ acts on $\Mt$, the bundle $L^s$ is not quite equivariant. If $L^s$ were $\G$-equivariant, it would simply be the lift of a bundle from $\M$, and that would imply an integrality condition on $s\omega$. The problem is that the obvious action of $\G$ on $L^s$ is incompatible with the connection. Nevertheless, this can be corrected to a \emph{projective} action of $\G$ on sections of $L^s$ which is compatible with the connection. 

\subsection{Projective Representations and Twisted Algebras}
Let's first recall some definitions and properties associated with projective representations of a group.

A projective group representation only respects the group product modulo multiplication by complex numbers. If $U$ is a projective  representation of $\G$\/, then $U(\g_1)U(\g_2)$ is only proportional to (rather than equal to) $U(\g_1\g_2)$. The map $\sigma:\G\times\G\to\co^\times := \co\smallsetminus\{0\}$ defined by $\sigma(\g_1,\g_2):= U(\g_1)U(\g_2)U(\g_1\g_2)^{-1}$ 
satisfies the group cocycle property,
\beq
\sigma(\g_2,\g_3)\sigma(\g_1\g_2,\g_3)^{-1}\sigma(\g_1,\g_2\g_3)\sigma(\g_1,\g_2)^{-1}
= 1
\mbox.\label{cocycle.identity}\eeq
Thus $\sigma$ defines a group cohomology class in $H^2(\G;\co^\times)$; I shall refer to $\sigma$ as the twist cocycle of the projective representation $U$. 

A projective-unitary representation on a Hilbert space $\HH$ does induce a true representation of $\G$ on the algebras $\Li(\HH)$ and $\K(\HH)$. An element $\g\in\G$ acts by the automorphism $a\mapsto U(\g)aU(\g)^*$. The extraneous phases all cancel out.

Two projective representations should be considered equivalent if they are projectively equivalent. That is, $U$ is equivalent to the projective representation obtained by multiplying $U(\g)$ by a factor  $\rho(\g)$. This has the effect of multiplying $\sigma$ by the coboundary of $\rho$. Thus, equivalent projective representations determine the same cohomology class in $H^2(\G,\co^\times)$. Equivalent projective-unitary representations induce the same actions on $\Li(\HH)$ and $\K(\HH)$.
\medskip

From a discrete group, $\G$, and a $2$-cocycle, $\sigma:\G\times\G\to\co^\times$, we can construct a twisted group algebra, $\co[\G,\sigma]$. This is a complex vector space with a basis indexed by $\G$, just as for the ordinary group algebra $\co[\G]$. However, the product in $\co[\G,\sigma]$ is defined by,
\beq
[\g_1][\g_2] := \sigma(\g_1,\g_2)[\g_1\g_2]
\mbox,\label{twist}\eeq
for any $\g_1,\g_2\in\G$. The cocycle identity, eq.~\eqref{cocycle.identity}, is equivalent to the associativity of this algebra. 

If the cocycle $\sigma$ is valued in the unitary group $\T\subset\co^\times$, then $\co[\G,\sigma]$ is a $*$-algebra if we define the basis elements to be unitary. That is, define the involution to be the unique antilinear map ${}^*:\co[\G,\sigma] \to \co[\G,\sigma]$, such that $[\g]^*=[\g]^{-1}$, $\forall \g\in\G$. Cohomologous cocycles give isomorphic algebras; the isomorphism is given by multiplying each $[\g]$ by a phase factor.

As in the untwisted case, a twisted group algebra can be completed to a (reduced or maximal) twisted group \cs-algebra. In the general framework of twisted crossed products described in \cite{p-r}, these are reduced and maximal twisted crossed products of the group $\G$ with the \cs-algebra $\co$.

The maximal norm of $a\in \co[\G,\sigma]$ is the supremum of norms of images of $a$ in $*$-representations of $\co[\G,\sigma]$. The maximal \cs-algebra $\cs_{\max}(\G,\sigma)$ is the completion of $\co[\G,\sigma]$ in this norm. Note that a $*$-representation of $\co[\G,\sigma]$ is the same thing as a projective-unitary representation of $\G$ with twist $\sigma$.

As a vector space, $\co[\G,\sigma]$ can be completed to the Hilbert space $l^2(\G)$. The product in $\co[\G,\sigma]$ extends to an action of $\co[\G,\sigma]$ on $l^2(\G)$. This is the regular representation of $\co[\G,\ss]$ (or equivalently, the $\ss$-twisted regular representation of $\G$\label{reg.rep}). The reduced \cs-algebra, $\csr(\G,\sigma)$ is the completion of $\co[\G,\sigma]$ in the regular representation. 
\medskip

Consider the linear map $\tau:\co[\G,\sigma]\to\co$ defined by $\tau(1)=1$ and $\tau[\g]=0$ for $\g\neq e\in\G$.
\begin{lem}\label{trace}
$\tau$ is well defined and extends to a faithful, tracial state on $\csr(\G,\sigma)$. 
\end{lem}
\begin{proof}
As an operator in $\csr(\G,\sigma)$, $[e]$ is defined by its action on the basis elements $[\g]\in l^2(\G)$. By eq.~\eqref{twist},
\[
[e][\g]=\sigma(e,\g)[\g]
\mbox,\]
but by the cocycle condition (associativity) $\sigma(e,\g)=\sigma(e,e)$. So $[e]$ is proportional to the identity $1\in\csr(\G,\sigma)$. Thus $\tau$ is well defined on $\co[\G,\sigma]$, because with $\tau[e] = \sigma(e,e)$ we have the value of $\tau$ on all basis elements.

Normalization is the condition $\tau(1)=1$ in the definition.

To see that $\tau$ is a trace, consider basis elements. For $\g,\g'\in\G$, $\tau([\g\g'])$ and $\tau([\g'\g])$ are each nonzero if and only if $\g'=\g^{-1}$. Now,
\[
[\g][\g^{-1}] = \sigma(\g,\g^{-1}) \in\co
\]
is central. So,
\[
[\g][\g^{-1}] = ([\g][\g^{-1}])[\g]^{-1}[\g] = [\g]^{-1}([\g][\g^{-1}])[\g] = [\g^{-1}][\g]
\mbox.\]
Thus, $\tau([\g][\g'])$ and $\tau([\g'][\g])$ are equal if nonvanishing, and hence always equal. Therefore, $\tau$ is a trace.

The inner product $\tau(a^*b)$ for $a,b\in\co[\G,\sigma]$ is just the standard inner product for $l^2(\G)$ since $\tau([\g]^*[\g])=1$ and $\tau([\g]^*[\g'])=0$ for $\g\neq\g'$ (this shows positivity of $\tau$). Therefore, the regular representation is just the GNS representation constructed from $\tau$. This shows that $\tau$ extends to $\csr(\G,\sigma)$ and is faithful.
\end{proof}

\subsection{The Cocycle}
\label{cocycle}
Returning to the structures described on p.~\pageref{twist}, we seek to define a projective right action of $\G$ on sections of $L^s$. The reason for preferring a \emph{right} action will become apparent later. The obvious right action is given by the pullback of sections; for any $\g\in\G$, $\psi\in\Gamma(\Mt,L^s)$, and $x\in\Mt$,
\[
(\g^*\psi)(x) := \psi(\g x)
\mbox.\]
However, this is incompatible with the connection $\nabla_s \equiv d-isA$,
\[
\g^*([d-isA]\psi) = (d-is\g^*A)\g^*\psi 
\neq (d-isA)\g^*\psi
\mbox.\]
The discrepancy here is $is(\g^*A-A)\g^*\psi$. Fortunately,
\begin{lem}
    For every $\g\in\G$, the $1$-form $\g^*A-A$ is exact.
\end{lem}
\begin{proof}
\[
d(\g^*A) = \g^*dA = \g^*(\pi^*\omega) = \pi^*\omega = dA
\]
So, $\g^*A-A$ must be closed.  By definition, the universal cover
is simply connected, thus $H^1(\Mt)=0$ and any closed $1$-form must be exact. 
\end{proof}

This means that we can choose functions $\phi_\g\in\C^\infty(\Mt)$ for
all $\g\in\G$, such that $d\phi_\g=A-\g^*A$. Such a choice will be fixed for the remainder.
\begin{definition}
    For any $\psi\in\Gamma(\Mt,L^s)$ and $\g\in\G$,
    \beq
    \psi\cdot\g := e^{is\phi_\g}\g^*\psi
    \mbox.\label{adef}\eeq 
\end{definition}
\begin{lem}\label{twisted}
    Equation \eqref{adef} defines a projective right action of $\G$ on 
    $\Gamma(\Mt,L^s)$ compatible with the connection, $\nabla_s \equiv d-isA$.
\end{lem}
\begin{proof}
On the one hand,
\[
(\psi\cdot\g_1)\cdot\g_2 = e^{is\phi_{\g_2}}\g_2^*(e^{is\phi_{\g_1}}\g_1^*\psi) 
= e^{is(\phi_{\g_2}+\g_2^*\phi_{\g_1})}(\g_1\g_2)^*\psi
\mbox.\]
On the other hand,
\[
\psi\cdot(\g_1\g_2) = e^{is\phi_{\g_1\g_2}}(\g_1\g_2)^*\psi
\mbox.\]
These differ by a factor of, 
\[
e^{is(\phi_{\g_2}+\g_2^*\phi_{\g_1}-\phi_{\g_1\g_2})}
\mbox,\] 
which is a constant on $\Mt$ because
\begin{align*}
d(\phi_{\g_2}+\g_2^*\phi_{\g_1}-\phi_{\g_1\g_2}) 
&= d\phi_{\g_2} + \g_2^*d\phi_{\g_1}-d\phi_{\g_1\g_2}\\
&= \g_2^*A-A+\g_2^*(\g_1^*A-A)-(\g_1\g_2)^*A+A\\
&= 0
\mbox.\end{align*}
Therefore this is a projective action.

For any $\psi\in\Gamma(\Mt,L^s)$ and $\g\in\G$,    
\begin{align*}
\nabla_s (\psi\cdot\g) &=  (d-isA)e^{is\phi_\g}\g^*\psi \\
   & = e^{is\phi_\g}(d-isd\phi_\g-isA)\g^*\psi
    \\&= e^{is\phi_\g}(d-is\g^*A+isA-isA)\g^*\psi
    \\&= e^{is\phi_\g}\g^*([d-isA]\psi) \\
&= (\nabla_s\psi)\cdot\g
\end{align*}
Therefore this projective action commutes with the connection.    
\end{proof}
\begin{definition}
\beq\label{c.def}
c(\g_1,\g_2) := \phi_{\g_2}+\g_2^*\phi_{\g_1}-\phi_{\g_1\g_2}
\eeq 
\beq\label{ss.def}
\ss(\g_1,\g_2):= e^{isc(\g_1,\g_2)}
\eeq
\end{definition}
\begin{cor}
    $c$ and $\ss$ are cocycles, and $\ss$ is the twist of the action \eqref{adef}.
\end{cor}

\begin{thm}\label{cohomology.isomorphism}
The cohomology class of $c$ is determined uniquely by the cohomology class of $\omega$. This correspondence gives an isomorphism from 
the kernel of $\pi^*:H^2(\M)\to H^2(\Mt)$ to $H^2(\G)$.
\end{thm}
\begin{proof} 
The spaces $\Omega^*(\Mt)$ are right $\G$-modules with the action by pull-back. Let  $E_0^{p,q} := C^p(\G,\Omega^q(\Mt))$ be the first quadrant double complex of group cohomology (bar) complexes with coefficients in the de\,Rham complex of $\Mt$. The vertical (degree $(0,1)$) differential, $\dv$, is the de\,Rham differential; the horizontal (degree $(1,0)$) differential, $\dh$, is the bar coboundary. Note that the vertical edge $E_0^{0,*}$ is simply the de\,Rham complex for $H^*(\Mt)$.

The left spectral sequence of this double complex starts by taking the  horizontal cohomology (in this case, the group cohomology). Because the action of $\G$ on $\Mt$ is proper, the complex $\Omega^*(\Mt)$ consists of free $\co[\G]$-modules. The group cohomology thus gives the $\G$-invariant subcomplex $\Omega^*(\M)\subset \Omega^*(\Mt)$ in horizontal degree $0$, and $0$ in higher degrees. In other words, the left spectral sequence collapses to one column $\one E_1^{0,*} = \Omega^*(\M)$. It therefore converges to  $\one E_\infty^{0,*} = \one E_2^{0,*} = H^*(\M)$. The obvious inclusion $\Omega^*(\M)\subset\Omega^*(\Mt) = E_0^{0,*} \subset E_0^{*,*}$ gives the identification of the de\,Rham cohomology $H^*(\M)$ with the cohomology of the total complex\footnote{Note that because the differentials commute, the total differential should be taken as  $d=\dv\pm\dh$ with the sign alternating by columns.} $\tot E_0^{*,*}$. The natural map $\pi^*: H^*(\M)\to H^*(\Mt)$ is given by the inclusion $\Omega^*(\M)\subset\Omega^*(\Mt)$. In terms of the complex $\tot E_0^{*,*}$, $\pi^*$ is given by the map $\tot E_0^{*,*} \to E_0^{0,*} = \Omega^*(\Mt)$ that extracts the first column part of a cochain.

The right spectral sequence must also converge to $H^*(\M)$. This sequence begins by taking the vertical cohomology, so
\[
\two E_1^{p,q} = C^p(\G,H^q(\Mt))
\mbox.\]
Because $\Mt$ is simply connected, $H^1(\Mt)=0$, so the $q=1$ row vanishes, $\two E_1^{*,1} =0$. Because $\Mt$ is connected, $H^0(\M)=\co$, and the bottom row is $\two E_1^{*,0} = C^*(\G)$, the group cohomology complex (not to be confused with the group \cs-algebra). In degree $0$, the cohomology is just the space of coboundaries; so there is an inclusion $C^*(\G) = \two E_1^{*,0} \subset E_0^{*,*}$. This gives a natural map in cohomology,
\[
k^*: H^*(\G) \to H(\tot E_0^{*,*}) = H^*(\M)
\mbox.\]
By this construction, it is clear that $\pi^*\circ k^* =0$ except in degree $0$, since the bottom row and right column of $E_0^{*,*}$ only intersect at $E_0^{0,0}$.

Now, we can reexamine the construction of $c$ in terms of the double complex $E_0^{*,*}$. First, $\omega\in\Omega^2(\M) \subset E_0^{0,2}$. The hypothesis that $[\pi^*\omega] = 0 \in H^2(\Mt)$ means that $\omega$ is trivial in the vertical cohomology; thus there exists $A\in E_0^{0,1}$ such that $\dv A = \omega$. The $\dh$-coboundary, $\dh A$, is the map $\g\mapsto A-\g^*A$. This is a $\dv$-cocycle because 
\[
\dv(\dh A) = \dh\dv A = \dh\omega = 0
\mbox.\]
Because $H^1(\Mt)=0$, the columns of $E_0^{*,*}$ are exact at degree $1$, therefore $\dh A$ is a $\dv$-coboundary and there exists $\phi\in E_0^{1,0}$ such that $\dv\phi = \dh A$. Finally, eq.~\eqref{c.def} is equivalent to $c=\dh\phi\in E_0^{2,0}$; this is explicitly a $\dh$ cocycle but it is also a $\dv$ cocycle because 
\[
\dv c = \dv\dh\phi =\dh\dv\phi = \dh\dh A =0
\mbox,\]
so $c\in \two E_1^{2,0} = C^2(\G)$. Schematically, the construction is:
\[
\begin{CD}
\omega \\
@AAA \\
A @>>> \dh A \\
@. @AAA \\
@. \phi @>>> c
\end{CD}
\]
By this construction, $c$ is cohomologous to $\omega$ in $\tot E_0^{*,*}$. This shows that $k^*[c] = [\omega]$. 

To check the uniqueness of the cohomology class $[c]\in H^2(\G)$ it is sufficient to compute the kernel of $k^*$ in degree $2$. So, suppose that $[\omega] =0 \in H^2(\M)$. This means that we can choose $A\in\Omega^1(\M) = \one E_1^{0,*}$, i.~e., such that $\dh A=0$. In this case, $\dv\phi = \dh A =0$, so $[\phi]\in \two E_1^{1,0} = C^1(\G)$ and $[c]=0\in H^2(\G)$. 

This shows that there is an exact sequence,
\beq
0 \to H^2(\G) \stackrel{k^*}{\longrightarrow} H^2(\M) \stackrel{\pi^*}{\longrightarrow} H^2(\Mt)
\mbox.\eeq
The construction of $[c]$ from $[\omega]$ is just the inverse of the isomorphism \[
k^*:H^2(\G) \isom \ker \pi^*
\mbox.\]
\end{proof}
The group cohomology can also be constructed as the cohomology of a classifying space, $B\G$. Since $\G\equiv\pi_1(\M)$, there is a canonical (up to homotopy) classifying map $k:\M\to B\G$. The above map $k^*$ is just the induced map on cohomology $k^*:H^*(\G)=H^*(B\G)\to H^*(\M)$. In particular, if $\Mt$ is contractible, then we can take $\M$ as a classifying space for $\G$ and the classifying map is the identity map. 

Also note that this proof is not specific to de\,Rham cohomology. The construction of $c$ from $\omega$ in terms of the double complex $E_0^{*,*}$ would work just as well with any coefficients.

\subsection{Continuous Field} 
\label{continuous.field}
The reduced, twisted group \cs-algebras are of such importance here that I shall denote them more concisely as $\Bs:=\csr(\G,\ss)$. Denote the tracial states (from Lem.~\ref{trace}) as $\tau_s:\Bs\to\co$.

The collection of twisted group \cs-algebras, $\{\Bs\}_{s\in\R}$, can be understood as a continuous field of \cs-algebras. In order to show this I shall use a more general result.

Recall that a function is upper semicontinuous if its value at any point $x$ is greater than or equal to its lim-sup at $x$. It is lower semicontinuous if its value at any point $x$ is less than or equal to its lim-inf at $x$. A supremum (infimum) of continuous functions is lower (upper) semicontinuous.
\begin{lem}\label{field}
Let $X$ be a locally compact space and $A$ a \cs-algebra such that $\C_0(X)$ is contained in the center of $A$. Let  $\rho:A\to\C_0(X)$ be a $\C_0(X)$-linear map.
\begin{enumerate}
\item 
For each $x\in X$, 
\[\A_x := A/\C_0(X\smallsetminus \{x\})A\]
is a well-defined \cs-algebra.  
\item
For each $x\in X$, $\rho$ uniquely determines a linear map $\rho_x: \A_x\to\co$ such that, $\rho_x\circ\Po_x = \Po_x\circ\rho$, where $\Po_x:A\onto \A_x$ is the quotient map.
\item
If $\C_0(X)A=A$ then the $\A_x$'s form an upper semicontinuous field of \cs-algebras with $A = \Gz(X,\A)$.
\item
If, furthermore, the maps $\rho_x$ are faithful states, then $\A$ is a continuous field of \cs-algebras.
\end{enumerate}
\end{lem}
\begin{proof}
Note that $\C_0(X\smallsetminus \{x\})\subset\C_0(X)$ is simply the (closed) ideal of functions which vanish at $x\in X$. Because $C_0(X)$ is contained in the center of $A$, $\C_0(X\smallsetminus \{x\})A$ is a closed, two-sided ideal and so the \cs-algebra $\A_x$ is well-defined. 

The induced map $\rho_x$ should fit in the exact/commutative diagram,
\[
\begin{CD}
0 @>>> \C_0(X\smallsetminus \{x\})A @>>> A @>{\Po_x}>> \A_x @>>> 0\\
@. @VVV @V{\rho}VV @V{\rho_x}VV \\
0 @>>> \C_0(X\smallsetminus \{x\}) @>>> \C_0(X) @>>> \co @>>> 0
\mbox.\end{CD}
\]
To see that $\rho_x$ is well defined, consider any $\alpha\in C_0(X\smallsetminus \{x\})$ and $a\in A$. By the $\C_0(X)$-linearity of $\rho$,
\[
\rho(\alpha a) = \alpha\,\rho(a) \in \C_0(X\smallsetminus \{x\})
\mbox.\]

For $a\in A$, let $N_a:X\to\R_+$ be the norm function defined by $N_a(x):=\Norm{\Po_x(a)}$. These norms can be computed as,
\[
\Norm{\Po_x(a)} = \inf \left\{\norm{a+b} \mid b\in \C_0(X\smallsetminus\{x\})A\right\}
\mbox.\]
For $b\in \C_{\mathrm b}(X,A)$, an $A$-valued function, let $\nu_b \in \C_{\mathrm b}(X)$ be the function defined by $\nu_b(x)=\norm{b(x)}$. If we assume that $\C_0(X)A=A$, then multiplication gives a surjective homomorphism $\C_0(X)\otimes A \onto A$. A surjective homomorphism extends to the multiplier algebras. The algebra $\C_{\mathrm b}(X,A)$ is contained in the multiplier algebra of $\C_0(X,A)= \C_0(X)\otimes A$, so we have a map from $\C_{\mathrm b}(X,A)$ to the multiplier algebra of $A$. Let $J \subset \C_{\mathrm b}(X,A)$ be the kernel. looking back at this construction, $J$ is the ideal generated by functions in $\C_{\mathrm b}(X \times X)$ which vanish along the diagonal. 

If we identify $a$ with the corresponding constant function in $\C_{\mathrm b}(X,A)$, then the function $N_a$ is given by,
\[
N_a = \inf \left\{\:\nu_{a + b} \mid b\in J\;\right\}
\mbox.\]
Since this is an infimum of continuous functions, it is an upper semicontinuous function. It also goes to $0$ at $\infty$ because $A=\C_0(X)A$. Furthermore, because  $A\into \C_{\mathrm b}(X,A)/J$ is injective, $\sup_{x\in X}N_a(x) = \norm a$. This verifies the axioms for an upper semicontinuous field of \cs-algebras. 

If the induced maps $\rho_x$ are faithful states, then these norms can also be computed as,
\[
N_a(x) \equiv \norm{\Po_x(a)} = \sup\left\{\Abs{\rho_x\left[\Po_x(a)b\right]}\big| b\geq0\in \A_x, \,
\rho_x(b)\leq 1\right\}
\mbox.\]
This can be written as,
\[
N_a = \sup\left\{\abs{\rho(ab)}\mid b\geq0\in A, 
\rho(b)\leq 1\right\}
\mbox.\]
This is a supremum of continuous functions, it is therefore a lower semicontinuous function. Since we already know that $N_a$ is upper semicontinuous, this means it is continuous. Therefore $\A$ is a continuous field.
\end{proof}
Using the KSGNS tensor product and the algebra of bounded-adjoint\-able operators (see Section \ref{Hilbert.modules}) the hypothesis on $\rho$ can be stated more succinctly as $\rho : A\to\C_0(X)$ must be a completely positive $\C_0(X)$-linear retraction such that the natural homomorphism $A\to \Li_{\C_0(X)}(A\otimes_\rho\C_0(X))$ is injective. However, the statement of the hypothesis as given above is more practical.

In the cases I am considering here, the map $\rho$ is actually a trace.
\medskip

Using the cocycle $c$, we can define an extension of $\G$,
\[
0\to \R\to \tilde\G \to \G \to 0
\mbox.\]
\begin{definition} 
Let $\tilde\G$ be the space $\G\times\R$ with the product,
\[
(\g_1,r_1)(\g_2,r_2) = \left(\g_1\g_2,r_1+r_2+c(\g_1,\g_2)\right)
\mbox.\]
\end{definition}
Note that this is more general than the situation at hand. In fact this works for any discrete group $\G$, real $2$-cocycle $c$, $\ss=e^{isc}$, and $\Bs=\cs(\G,\ss)$.
\begin{thm}\label{Bfield}
    The algebras $\Bs$ are the fibers of a continuous field of \cs-algebras 
    over $\R$ (which I denote as $\B$) with $\Gz(\R,\B)=\csr(\tilde\G)$.
\end{thm}
\begin{proof}
Because $\R\subset\tilde\G$ is an open subgroup (in fact, the identity component of $\tilde\G$) and central, $\cs(\R) \cong \C_0(\R)$ is contained in the center of $\csr(\tilde\G)$.

The algebra $\csr(\tilde\G)$ is the completion of $\C^\infty_{\mathrm c}(\tilde\G)$ in the regular representation.  The representation space, $L^2(\tilde\G)$, decomposes as a direct integral of copies of $l^2(\G)$ corresponding to the irreducible representations of $\R$, which are in turn the points of $\Spec[\cs(\R)]$.

For $s\in\R=\Spec[\cs(\R)]$, the fiber algebra is the image of $\csr(\tilde\G)$ in the corresponding representation on $l^2(\G)$. This is the regular representation of $\co[\G,\ss]$. Thus, the fiber algebra is $\Bs\equiv\csr(\G,\ss)$.

Define a map $\tau:\csr(\tilde\G)\to \cs(\R)\cong\C_0(\R)$ as follows. An element $b\in\csr(\tilde\G)$ is given by a distribution on $\tilde\G$. Restricting this to a distribution on the identity component, $\R\subset\tilde\G$, gives the element $\tau(b)\in\cs(\R) = \C_0(\R)$. This construction is $\R$-covariant, therefore $\tau$ is $\cs(\R)$-linear. The induced map $\Bs\to\co$ is simply $\tau_s$ as that is also defined by restriction to the identity component.

By Lem.~\ref{trace}, $\tau_s$ is a faithful state, so the requirements of Lem.~\ref{field} are satisfied, and $\B$ is a continuous field.
\end{proof}

This is not true of the maximal twisted \cs-algebras, $\csm(\G,\ss)$, because there is no general analogue of the family of faithful states $\tau_s$. In general, the maximal \cs-algebras merely form an upper semicontinuous field.

\section{Hilbert \cs-Modules}
\label{Hilbert.modules}
Let $B$ be a \cs-algebra. 
\begin{definition}
A \emph{Hilbert $B$-module}, $\HH$, (see \cite{lan} and references therein) is a
right $B$-module with a $B$-valued inner product, $\langle\,\cdot\,\vert\,\cdot\,\rangle:\HH\times\HH\to B$ such that:
\begin{enumerate}
\item
The inner product is antilinear in the first argument and linear in the second argument.
\item
For $\psi,\varphi\in\HH$, $\left(\langle\psi\vert\varphi\rangle\right)^* =
\langle\varphi\vert\psi\rangle$. 
\item
$\langle\psi\vert\psi\rangle\geq0$.
\item
For $b\in B$, $\langle\psi\vert\varphi b\rangle =
\langle\psi\vert\varphi\rangle b$.
\item
The semi-norm defined by
\beq\label{module.norm}
\norm\psi := \norm{\langle\psi\vert\psi\rangle}^{1/2}
\eeq
is a norm and $\HH$ is complete with respect to it.
\end{enumerate}
\end{definition}
The concept of Hilbert \cs-module generalizes that of Hilbert space; a Hilbert $\co$-module is precisely a Hilbert space. Another elementary but important example is that $B$ is itself a Hilbert $B$-module.

The algebra of bounded, linear operators on a Hilbert space generalizes to:
\begin{definition}
The algebra of \emph{bounded-adjoint\-able} operators on a Hilbert $B$-module $\HH$ is
\[
\Li_B(\HH) := \left\{a:\HH\to\HH \mid \exists a^*:\HH\to\HH,\; \forall\psi,\varphi\in\HH,\; \langle a^*(\psi)\vert\varphi\rangle = \langle\psi\vert a\varphi\rangle \right\}
\]
\end{definition}
This is a \cs-algebra. When $B=\co$ and $\HH$ is a Hilbert space it is $\Li_\co(\HH)=\Li(\HH)$, the algebra of bounded operators. The rather succinct ``adjointability'' condition implies that $a:\HH\to\HH$ is a $\co$-linear map and is $B$-linear in the sense that for all $\psi\in\HH$ and $b\in B$,
$a\lvert\psi b\rangle = (a\lvert\psi\rangle) b$.  I am using the term ``bounded-adjoint\-able'' rather than simply ``adjointable'' to avoid confusion with unbounded operators that have adjoints in the sense described below.

\begin{definition}
The algebra of \emph{compact} operators on a Hilbert $B$-module is the \cs-ideal $\K_B(\HH)\subseteq\Li_B(\HH)$ densely spanned by operators of the form,
\beq\label{rank1}
\lvert\psi\rangle\langle\varphi\rvert
\eeq
for $\psi,\varphi\in\HH$. 
\end{definition}
I am using the Dirac ``bra-ket'' notation here in which $\left(\lvert\psi\rangle\langle\varphi\rvert\right)\lvert\chi\rangle = \lvert\psi\rangle\left(\langle\varphi\vert\chi\rangle\right)$. The spaces of bounded-adjoint\-able and compact maps from one Hilbert $B$-module to another are defined identically.

The term ``compact'' is used because in the case when $B=\co$ and $\HH$ is a Hilbert space, $\K_\co(\HH)=\K(\HH)$ is the standard algebra of compact operators. I shall use the adjective ``compact'' in this generalized sense, despite the fact that such a compact operator is not compact in the sense of operators on a Banach space. In this same spirit, finite linear combinations of operators of the form \eqref{rank1} can be considered ``finite rank''\/.

In the case that $\HH=B$, $\K_B(B)=B$ and $\Li_B(B)$ is the multiplier algebra of $B$. Indeed, $\Li_B(\HH)$ is always the multiplier algebra of $\K_B(\HH)$. The algebras coincide, $\K_B(\HH)=\Li_B(\HH)$, if and only if $B$ is unital and $\HH$ is finitely generated.

The Kasparov stabilization theorem says that any countably generated Hilbert $B$-module, $\HH$, is isomorphic to a submodule of the Hilbert $B$-module $B^{\oplus\infty} = l^2(\N)\otimes B$. In fact it is a direct summand, $\HH \oplus B^{\oplus\infty} \cong B^{\oplus\infty}$.

\begin{definition}
A \emph{Hilbert $A$-$B$-bimodule}, $\HH$, is a Hilbert $B$-module with a 
$*$-homo\-morph\-ism $A\to\Li_B(\HH)$. A \emph{Hilbert $B$-bimodule} is a Hilbert $B$-$B$-bimodule.
\end{definition}
\medskip

An important, and closely related concept is that of a completely positive map. Let $A$ and $B$ be \cs-algebras. An element of a \cs-algebra is positive if it can be written as the square of a self-adjoint element. A linear map $\rho:A\to B$ is positive if it maps positive elements of $A$ to positive elements of $B$. Any linear map $\rho:A\to B$ extends trivially to a map of matrix algebras $\rho\otimes\id : M_m(A)\to M_m(B)$, given by applying $\rho$ entrywise.
\begin{definition}
A linear map $\rho:A\to B$ is \emph{completely positive} if for any $m\in\N$, $\rho\otimes\id :M_m(A)\to M_m(B)$ is positive.
\end{definition}
If either $A$ or $B$ is commutative, then positivity is equivalent to complete positivity. 
If $\HH$ is a Hilbert $A$-module, $\HH'$ a Hilbert $B$-module, and $\rho : A\to\Li_B(\HH)$ a completely positive and strict map then we have:
\begin{definition}
The \emph{KSGNS tensor product} (see \cite{lan}) $\HH\otimes_\rho\HH'$ is the completion of the algebraic tensor product $\HH\algtensor\HH'$ with the inner product defined by
\[
\langle\psi\otimes\psi'\vert\varphi\otimes\varphi'\rangle 
= \braket{\psi'}{\rho\left(\langle\psi\vert\varphi\rangle\right)}{\varphi'}
\mbox,\]
for $\psi,\varphi\in\HH$, $\psi',\varphi'\in\HH'$.
\end{definition}
Completion must be understood here to involve first quotienting by the subspace of elements with (semi-)norm \eqref{module.norm} equal to $0$, and then taking the closure. The condition of strictness is trivial for unital \cs-algebras.

This has several import special cases and properties. 
\begin{enumerate}
\item
If $\rho$ is actually a homomorphism, then this makes $\HH'$ a Hilbert $A$-$B$-bimodule. The KSGNS tensor product is then the (more elementary) tensor product over $A$, $\HH\otimes_\rho\HH'=\HH\otimes_A\HH'$.
\item
If $\HH'=B$ and $\rho:A\to B$ is a homomorphism then $\HH\otimes_\rho B = \rho_*(\HH)$ is the push-forward of $\HH$ by $\rho$ to a Hilbert $B$-module.
\item
The tensor product is a functor. There is a canonical homomorphism $\otimes_\rho H' : \Li_A(\HH)\to\Li_B(\HH\otimes_\rho\HH')$. Thus, if $\HH$ is a bimodule, then so is this tensor product.
\item
The tensor product is associative. Thus a general tensor product factorizes as $\HH\otimes_\rho\HH' = \HH\otimes_A (A\otimes_\rho B) \otimes_B \HH'$. Thus any tensor product can be constructed from two special cases: when $\rho$ is a homomorphism, and when the Hilbert \cs-modules are the algebras themselves.
\item
If $B=\co$ and $\rho:A\to\co$ is a state, then $A\otimes_\rho\co$ is the GNS Hilbert space and the canonical map $A\to\Li(A\otimes_\rho\co)$ is the GNS representation constructed from $\rho$.
\item
If $B=\co$ and $\rho$ is a faithful state, then $\HH\otimes_\rho\co$ is a completion of $\HH$ to a Hilbert space. There are inclusions,
\[
\K_A(\HH) \subseteq \Li_A(\HH) \subseteq \Li(\HH\otimes_\rho\co)
\mbox.\]
\end{enumerate}
\medskip

\begin{lem}\label{ModuleField}
Let $\B$ be a continuous field of unital \cs-algebras over a locally compact, Hausdorff space $X$.
Any countably generated Hilbert $\Gz(X,\B)$-module is isomorphic to the space $\Gz(X,\Hi)$ of $\C_0$-sections of a continuous field, $\Hi$, of Banach spaces over $X$; the fiber, $\Hi_x$, at $x\in X$ is a Hilbert $\B_x$-module.

There exists a continuous field of \cs-algebras, $\K_\B(\Hi)$, over $X$, with fiber $\K_{\B_x}(\Hi_x)$ over $x\in X$, and spaces of sections,
\[
\Gz[X,\K_\B(\Hi)] = \K_{\Gz(X,\B)}(\Gz[X,\Hi])
\]
and
\[
\Gb[X,\K_\B(\Hi)] \subseteq \Li_{\Gz(X,\B)}(\Gz[X,\Hi])
\mbox.\]
\end{lem}
\begin{proof}
Let $H$ be any Hilbert $\Gz(X,\B)$-module. By the Kasparov stabilization theorem, $H$ can be embedded as a direct summand of $\Gz(X,\B)^{\oplus\infty}$\/. If we identify $H$ with this submodule, then it is the image of a projection $e\in \Li_{\Gz(X,\B)}(\Gz[X,\B]^{\oplus\infty})$.

Let $\Po_x:\Gz(X,\B)\onto \B_x$ be the evaluation homomorphism at $x$. For any $x\in X$, the image of the projection $\Po_x(e)$ is the pushed forward Hilbert $\B_x$-module, $\Po_x^*(H) = H \otimes_{\Po_x}\B_x$. Call this $\Hi_x$. These are subspaces of the fibers of the continuous field $\B^{\oplus\infty}$. The space of $\C_0$-sections of $\B^{\oplus\infty}$ taking their values in these subspaces is $H=\Im e$. This is closed, therefore these subspaces form a continuous (sub)field of Banach spaces over $X$.

Let $\co^{\oplus\infty}=l^2(\N)$ be the canonical countably infinite-dimensional Hilbert space. The algebra of compact operators on $\Gz(X,\B^{\oplus\infty}) = \Gz(X,\B)\otimes \co^{\oplus\infty}$ is
\[
\K_{\Gz(X,\B)}(\Gz[X,\B^{\oplus\infty}]) 
= \Gz(X,\B)\otimes \K(\co^{\oplus\infty})
= \Gz[X,\B\otimes \K(\co^{\oplus\infty})]
\mbox.\]

The algebra of compact operators on $\Hi_x = \Im \Po_x(e)$ is 
\[
\K_{\B_x}(\Hi_x) = \Po_x(e) [\B_x\otimes\K(\co^{\oplus\infty})] \Po_x(e)
\mbox.\]
The subspace of $\C_0$-sections of $\B\otimes\K(\co^{\oplus\infty})$ taking their values in these subalgebras is 
\[
e\, \Gz[X,\B\otimes\K(\co^{\oplus\infty})] e = \K_{\Gz(X,\B)}(H)
\mbox.\]
This is a \cs-subalgebra and so these \cs-algebras form a continuous field, $\K_\B(\Hi)$, over $X$ as claimed.

Let $a \in \Gb[X,\K_B(\Hi)]$. We can multiply this section (pointwise) with a section in $\Gz[X,\K_B(\Hi)]$. This multiplication preserves continuity, and because $a$ is bounded, it preserves $\Gz[X,\K_B(\Hi)]$. Pointwise multiplication is associative, so $a$ is a multiplier of $\Gz[X,\K_B(\Hi)] = \K_{\Gz(X,\B)}(H)$. However, the multiplier algebra of $\K_{\Gz(X,\B)}(H)$ is $\Li_{\Gz(X,\B)}(H)$. So, $\Gb[X,\K_B(\Hi)] \subseteq \Li_{\Gz(X,\B)}(H)$.
\end{proof}
\medskip

As with Hilbert spaces, an unbounded operator on a Hilbert \cs-module is not really a linear map defined on the entire space. Let $\HH$ be a Hilbert $B$-module. An unbounded operator on $\HH$ is a $B$-linear map $\Xi:\dom \Xi \to \HH$, defined on some $B$-submodule $\dom \Xi\subset \HH$. The operator $\Xi$ is densely defined if $\dom \Xi\subset \HH$ is a dense subspace. The operator $\Xi$ is closed if $\dom \Xi$ is a Hilbert $B$-module with the ``graph'' inner product,
\[
\langle\psi\vert\varphi\rangle_{\dom \Xi} := \langle\psi\vert\varphi\rangle + \langle \Xi\psi\vert \Xi\varphi\rangle
\mbox.\]
Equivalently, $\Xi$ is closed if $\dom \Xi$ is complete under the corresponding ``graph norm''\/. 

The action of the adjoint $\Xi^*$ of $\Xi$ is defined by the condition that $\forall \varphi\in\dom \Xi$,
\[
\langle \Xi^*\psi\vert\varphi\rangle = \langle\psi\vert \Xi\varphi\rangle
\mbox.\]
The domain $\dom \Xi^*$ is the set of $\psi$ for which this condition can be satisfied.

Just as merely bounded operators are not satisfactory on a Hilbert \cs-module, unbounded operators require an additional condition in order to be adequately well behaved, analogous to the condition of adjointability. Again, see \cite{lan}.
\begin{definition}
A \emph{regular} operator $\Xi$ on $\HH$ is a densely defined, closed unbounded operator such that $\Xi^*$ is densely defined and $1+\Xi^*\Xi$ has dense image in $\HH$. 
\end{definition}
This implies that $1+\Xi^*\Xi$ is actually surjective.
Among other properties, regular operators have a well behaved functional calculus. If $\chi$ is a bounded, continuous function on the spectrum of a regular operator $\Xi$, then $\chi(\Xi)$ is bounded-adjoint\-able. In particular, a bounded regular operator is a bounded-adjoint\-able operator. The composition of a regular operator and a bounded-adjoint\-able operator is regular.

I will need to discuss whether the domains of two regular operators are equivalent, so it is necessary to clarify what this means. In the case of unbounded operators on a Hilbert space, the domains of two different operators may be equal as sets, but they will never be equal as Hilbert spaces (with the graph inner product). Equality as Hilbert spaces is too strong a condition. The domains may be equivalent as Hilbertian spaces --- that is, homeomorphic. We need an analogous concept for Hilbert \cs-modules.
\begin{definition}
\label{Hilbertian.def}
The category of Hilbertian $B$-modules is the collection of Hilbert $B$-modules with bounded-adjoint\-able maps as morphisms.
\end{definition}
An isomorphism of Hilbertian $B$-modules is thus an invertible bounded-ad\-joint\-able map. This is much less restrictive than isomorphism of Hilbert $B$-mod\-ules, which requires a unitary map.

A (unitary) isomorphism of Hilbert $B$-modules gives isomorphisms of the algebras of bounded-adjoint\-able (respectively, compact) operators as \cs-algebras. A (bounded-adjoint\-able) isomorphism of Hilbertian $B$-modules gives isomorphisms of the algebras of bounded-adjoint\-able (respectively, compact) operators as topological algebras, but not as \cs-algebras. 

As in the Hilbert space case, a \emph{core} of a closed unbounded operator $\Xi$ is a subspace $H_0\subseteq\dom \Xi$ such that $\Xi$ is the closure of the restriction $\Xi_{H_0}$. Equivalently, no other closed operator coincides with $\Xi$ over $H_0$. 
\begin{lem}\label{equivalence.equivalence}
Let $\Xi$ and $\Upsilon$ be regular operators on a Hilbert $B$-module $\HH$ with a common core. If the graph norms for $\Xi$ and $\Upsilon$ are equivalent over this core, then $\dom \Xi$ and $\dom \Upsilon$ are equal as Hilbertian $B$-modules.
\end{lem}
\begin{proof}
The domains $\dom\Xi$ and $\dom\Upsilon$ are the closures of the common core with the respective graph norms. Because the graph norms are equivalent, these domains are equal as subspaces of $\HH$.

The operator $1+\Xi^*\Xi$ is positive and thus has a (regular) square root. The graph inner product for $\dom \Xi$ can be written as,
\begin{align*}
\langle\psi\vert\psi\rangle_{\dom \Xi} 
&\equiv \langle\psi\vert\psi\rangle + \langle \Xi\psi\vert \Xi\psi\rangle \\
&= \braket\psi{(1+\Xi^*\Xi)}\psi \\
&= \left<(1+\Xi^*\Xi)^{1/2}\psi\bigm\vert(1+\Xi^*\Xi)^{1/2}\psi\right>
\mbox.\end{align*}
Thus $(1+\Xi^*\Xi)^{1/2} : \dom \Xi \isom \HH$ is an isomorphism of Hilbert $B$-modules.

Using this identification, and the analogous one for $\dom \Upsilon$, the identity map from $\dom \Xi$ to $\dom \Upsilon$ is equivalent to the operator,
\beq
\label{domain.isomorphism}
(1+\Upsilon^*\Upsilon)^{-1/2}(1+\Xi^*\Xi)^{1/2}
\eeq
on $\HH$. Because $\Upsilon$ is regular, $(1+\Upsilon^*\Upsilon)^{-1/2}$ is bounded-adjoint\-able; therefore \eqref{domain.isomorphism} is regular. The bounds on the ratio between the graph norms give bounds on the norms of \eqref{domain.isomorphism} and its inverse; therefore \eqref{domain.isomorphism} is bounded-adjoint\-able. By an identical argument, its inverse is bounded-adjoint\-able.
\end{proof}

\subsection{Bundles of Hilbert \cs-Modules}
\label{B.bundles}
In this section, $\M$ is a connected, smooth manifold. All derivatives here are defined using the appropriate norm topologies. For instance, if $B$ is a \cs-algebra, then the partial derivatives of a $B$-valued function are defined in the coordinate patches using norm limits in $B$. The space of smooth functions $\C^\infty(\M,B)$ consists of those functions for which arbitrary order partial derivatives exist in any smooth coordinate chart.

\begin{definition}
A (smooth) bundle of Hilbert $B$-modules $E$ over $\M$ is a (locally trivial) fiber bundle such that:
\begin{enumerate}
\item
The fibers of $E$ are Hilbert $B$-modules.
\item
The gluing functions are smooth functions valued in the algebra of bounded-adjoint\-able operators on the fiber. (This makes the space of smooth sections $\Gi(\M,E)$ meaningful.)
\item
For any two smooth sections $\psi,\varphi\in\Gi(\M,E)$, the inner product is a smooth $B$-valued function $\langle\psi\vert\varphi\rangle\in\C^\infty(\M,B)$.
\end{enumerate}
\end{definition}
If $\M$ is merely a topological space, then a bundle of Hilbert $B$-modules is defined in the same way, but with ``continuous'' replacing ``smooth''\/.
For $\M$ not necessarily compact, $\Gz(\M,E)$ is the set of continuous sections whose fiberwise norms give functions in $\C_0(\M)$.

An example of this is the \Miscenko-Fomenko ``line bundle'' \cite{m-f}, constructed from the universal covering of a manifold, $\M$. This is succinctly described as $\Mt\times_\G\csr(\G)$. It is a bundle of Hilbert \cs-modules of the reduced \cs-algebra, $\csr(\G)$, of the fundamental group. It is a quotient of $\Mt\times\csr(\G)$ by $\G$. This is locally trivial because the universal covering is locally trivial and the lift to $\Mt$ is $\Mt\times\csr(\G)$, which is globally trivial. It is referred to as a line bundle because the fibers are isomorphic to the algebra $\csr(\G)$ itself.

\begin{lem}\label{module.bundle}
$\Gz(\M,E)$ is a Hilbert $\C_0(\M,B)$-module. The algebras $\Li_B(E_x)$ for $x\in \M$ form a (locally trivial) bundle $\Li_B(E)$ over $\M$ and 
\[
\C_{\mathrm b}(\M) \subset  \Gb(\M,\Li_B[E]) \subseteq \Li_{\C_0(\M,B)}[\Gz(\M,E)]
\mbox.\]
\end{lem}
\begin{proof}
For $\psi,\varphi\in \Gz(\M,E)$, define the inner product to be the function,
\begin{align*}
\langle\psi\vert\varphi\rangle : \M &\to B \\
 x &\mapsto \langle\psi(x)\vert\varphi(x)\rangle
\mbox.
\end{align*}
This is well defined because the gluing functions preserve the inner products in the fibers. It is a continuous function, because it is continuous in every local trivialization. Since $\norm{\psi(x)},\norm{\varphi(x)}\in\C_0(\M)$, we have,
\[
\langle\psi\vert\varphi\rangle \in \C_0(\M,B) = \C_0(\M)\otimes B
\mbox.\]

Multiplication with $\C_0(\M,B)$ certainly preserves continuity and the fall-off of the norm; so, $\Gz(\M,E)$ is a $\C_0(\M,B)$-module. The inner product is linear with respect to this module structure, because this is true for the fibers and thus in local trivializations. This inner product inherits positivity from the inner products on the fibers. The associated (semi-)norm is
\[
\norm{\psi} = \sum_{x\in\M} \norm{\psi(x)}
\mbox,\]
which is indeed a norm. 

Since $\M$ is a manifold, there is a finite open cover $\{U_\alpha\}$ of $\M$ by $\dim\M+1$ contractible sets. By local triviality, the restriction of $E$ to any $U_\alpha$ is isomorphic to $U_\alpha\times \HH$ for some Hilbert $B$-module $\HH$. Thus $\Gz(U_\alpha,E) \cong \C_0(U_\alpha)\otimes\HH$ are identified as $B$-modules with the same inner product, but the right hand side is closed, so $\Gz(U_\alpha,E)$ is closed. A partition of unity shows that,
\[
\Gz(\M,E) = \sum_{\alpha = 1}^{\dim\M+1} \Gz(U_\alpha,E)
\mbox.\]
Therefore $\Gz(\M,E)$ is closed and is a Hilbert $\C_0(\M,B)$-module.

The local trivializations of $E$ give local trivializations for the collection of algebras, $\Li_B(E_x)$. The unitary gluing functions for $E$ give $*$-automorphism gluing functions, so this defines a bundle of \cs-algebras $\Li_B(E)$ over $\M$.

Let $a\in\Gb[\M,\Li_B(E)]$. There is an obvious pointwise multiplication with $\Gz(\M,E)$. This gives continuous sections over any local trivialization, therefore it gives globally continuous sections. Because $a$ is bounded, this also preserves the fall-off of the norm, so $a$ defines a linear map $a:\Gz(\M,E)\to\Gz(\M,E)$. This is bounded by the norm bound of $a$. The section $a^*\in\Gb[\M,\Li_B(E)]$, defined by $a^*(x)=[a(x)]^*$ is the adjoint of $a$, thus $a$ is bounded-adjoint\-able and $\Gb(\M,\Li_B[E]) \subseteq \Li_{\C_0(\M,B)}[\Gz(\M,E)]$.

The fibers of $\Li_B(E)$ are all unital algebras, so $\C_{\mathrm b}(\M) \subset \Gb(\M,\Li_B[E])$.
\end{proof}

Tensor products of bundles are defined, as usual, by taking tensor products of fibers in local trivializations and patching together.

A connection on a bundle of Hilbert $B$-modules, $E$, is a $B$-linear map, 
\[
\nabla: \Gi(\M,E)\to\Gi(\M,E\otimes T^*\M)
\mbox,\]
 which is given in local charts by the partial derivatives plus potentials which are smooth sections of $\Li_B(E)$. A connection is compatible with the inner product if the potentials are anti-self-adjoint. The bundle $\Li_B(E)$ inherits a connection from $E$ defined formally by $\nabla(a) := [\nabla,a]_-$.

An order $m$ differential operator on $E$ is a $B$-linear operator on smooth sections of $E$ which in any coordinate chart is given as a sum of smooth sections of $\Li_B(E)$ multiplied with partial derivative operators up to $m$'th order.

\begin{definition}
A \emph{$B$-$\frac12$-density bundle} over $\M$ is an equivalence class of pairs $(E,\mu)$ where $E$ is a Hilbert $B$-module bundle and $\mu$ is a smooth, nonsingular measure on $\M$. Such a pair has a $B$-valued inner product defined on sections of the bundle by integrating with $\mu$. Two such pairs are equivalent if this $B$-valued inner product is equal.
\end{definition}
In other words, this is a bundle of Hilbertian $B$-modules with a compatible inner product taking its values in $\Omega^n(\M,B)$. This is related to the concepts of ``half form'' and ``half density'' used in the geometric quantization literature \cite{woo}.

\begin{lem}\label{module.bundle2}
If $E$ is a $B$-$\frac12$-density bundle over $\M$, then the bundle of \cs-algebras $\Li_B(E)$ is well defined and
\[
\C_b(\M) \subset \Gb(\M,\Li_B[E]) \subset \Li_B[L^2(\M,E)]
\mbox.\]
\end{lem}
\begin{proof}
Since any two nonsingular measures differ by multiplying by a function, the Hilbert $B$-modules in equivalent pairs can only differ by multiplying the inner product by a function on $\M$. This does not change $\Li_B(E)$, thus $\Li_B(E)$ is well-defined for a $B$-$\frac12$-density bundle.

There always exists a strictly positive, smooth function on $\M$ with integral equal to $1$; by rescaling with such a function, we can always choose $\mu$ to be a probability measure.

The measure can also be thought of as a positive map $\mu:\C_0(\M)\to\co$. We have two inner products on $\Gz(\M,E)$. The $B$-valued inner product is equal to the $\C_0(\M)\otimes B$-valued inner product composed with $\mu\otimes\id$. So, the Hilbert module $L^2(\M,E)$ is just a KSGNS tensor product, $L^2(\M,E) = \Gz(\M,E)\otimes_{\mu\otimes\id}B$.

As with any KSGNS tensor product, there is a natural homomorphism 
\[
\Li_{\C_0(\M)\otimes B}(\Gz[\M,E]) \to \Li_B(L^2[\M,E])
\mbox.\]
 By Lem.~\ref{module.bundle}, we can apply this homomorphism to $\Gb(\M,\Li_B[E])$. Any nonzero section clearly gives a nonzero operator, so this map is injective.

Just as for a bundle of Hilbert $B$-modules, the fibers of $\Li_B(E)$ are unital algebras, so $\C_b(\M) \subset \Gb(\M,\Li_B[E])$.
\end{proof}

The concept of a $B$-$\frac12$-density bundle is useful when it is appropriate to regard the Hilbert $B$-module of sections as more fundamental than any specific Riemannian structure or connection. A connection on a $B$-$\frac12$-density bundle is given by a connection on an underlying Hilbert $B$-module bundle. The concept of compatibility between a connection and inner product is more awkward. 
If there is a canonical smooth measure, $\mu$, (such as a Riemannian volume form) then any $B$-$\frac12$-density bundle is given by a pair of the form $(E,\mu)$. The connection and inner product of this $B$-$\frac12$-density bundle are compatible, relative to $\mu$, if the connection is compatible with the inner product on $E$. For the rest of this section, $E$ will denote a $B$-$\frac12$-density bundle.

The rest of this section is mainly concerned with the analysis of Dirac operators. The main tool in this analysis will be Sobolev spaces.
Suppose that $E$ is a trivialized $B$-$\frac12$-density bundle over a flat $\R^n$ or $\T^n$ (torus); that is, $E$ is given by a pair such as $(\R^n\times\HH,dx^1\wedge\dots\wedge dx^n)$ and has the the trivial connection given by partial derivatives.
Then we can employ a Fourier transform or Fourier series. For $\psi \in \Gci(\R^n,E)$ (respectively $\Gi(\T^n,E)$) let $\tilde\psi(k)$ be the Fourier transform (or series coefficients). Normalize the inner product in the Fourier representation so that it coincides with the given inner product; that is, $\langle\tilde\psi\vert\tilde\varphi\rangle = \langle\psi\vert\varphi\rangle$. 
\begin{definition}
If $E$ is a trivialized $B$-$\frac12$-density bundle over $\R^n$ or $\T^n$ then the (order $m$) \emph{Sobolev inner product} is
\[
\langle\psi\vert\varphi\rangle_m := \braket{\tilde\psi}{1+\norm{k}^{2m}}{\tilde\varphi}
\]
for $\psi,\varphi\in\Gci(\R^n,E)$ or $\Gi(\T^n,E)$. The \emph{Sobolev norm} $\norm{\,\cdot\,}_m$ is the corresponding norm, and for any bounded open subset $U$, the  \emph{Sobolev space} $W^m(U,E)$ is the completion of $\Gzi(U,E)$ to a Hilbert $B$-module with this inner product.
\end{definition}
This is meaningful for any $m\in\R$, but we will only need it for $m\in\N$.

First, we need some properties of the Sobolev norms. The analysis in Lemmas \ref{Sobolev1} and \ref{Sobolev2} is based loosely on the proof of the basic elliptic estimate in \cite{fol1}. Use the symbol $\lesssim$ to denote that the ratio of the left hand side to the right hand side is bounded as a function of the obvious variable (usually $\psi$). This avoids a profusion of named constants. In this notation, two norms are equivalent if and only if both relations $\lesssim$ and $\gtrsim$ hold.
\begin{lem}\label{Sobolev1}
~
\begin{enumerate}
\item
The Sobolev norms for different choices of (constant) Riemannian metric are equivalent.
\item
If $U$ is a bounded open set and $\Xi$ is a differential operator of order $m$ defined over the closure $\bar U$ then $\norm{\Xi\psi} \lesssim \norm{\psi}_m$ for $\psi\in\Gzi(U,E)$.
\end{enumerate}
\end{lem}
\begin{proof}
The ratio of $1+\norm{k}^{2m}$ computed with two different metrics is bounded. Therefore Sobolev norms for different metrics are equivalent. This proves the first claim.

Let $\nabla$ be the trivial connection given directly by partial derivatives. The effect of partial derivatives on the Fourier transform is simply, $\widetilde{\nabla\psi}(k) = -i k\,\tilde\psi(k)$. Consider the second claim in a special case. Let $\Xi=\Xi(\nabla)$ be a degree $m$ polynomial in the partial derivative operators (with complex coefficients). This operator acts simply on the Fourier transform,
\[
\widetilde{\Xi\psi}(k) = \Xi(-ik)\tilde\psi(k)
\mbox.\]
So,
\[
\Norm{\Xi\psi} = \Norm{\langle \Xi\psi\vert \Xi\psi\rangle}^{1/2} = \Norm{\braket{\tilde\psi}{\Abs{\Xi(-ik)}^2}{\tilde\psi}}^{1/2} 
\mbox.\]
As $\Xi(-ik)$ is a polynomial of degree $m$, we have $\Abs{\Xi(-ik)}^2 \lesssim 1+\norm{k}^{2m}$. Therefore $\norm{\Xi\psi} \lesssim \norm{\psi}_m$.

An arbitrary $m$'th order differential operator over $\bar U$ can be written as a sum of sections of $\Li_B(E)$ multiplied with partial derivative operators up to $m$'th order. A section of $\Li_B(E)$ over $\bar U$ is bounded, thus the second claim follows.
\end{proof}

\begin{definition}
\label{Dirac.def}
Given a Riemannian metric on $\M$, a \emph{$B$-Dirac bundle} is a $B$-$\frac12$-density bundle $E$ with a connection, a $\Z_2$-grading, and a bundle homomorphism $c:T^*\M\to\Li_B(E)$ (the Clifford map) such that:
\begin{enumerate}
\item 
The connection is compatible with the inner product (relative to the Riemannian volume form) and the grading.
\item
For any $\xi\in\Omega^1(\M)$, 
\begin{enumerate}
\item $c(\xi)^2=\xi^2$, where $\xi^2$ is the Riemannian inner product of $\xi$ with itself,
\item if $\xi$ is real then $c(\xi)$ is self-adjoint, and
\item if $\M$ is even-dimensional then $c(\xi)$ is odd (with respect to the grading).
\end{enumerate}
\item
The map $c$ is parallel with respect to the Levi-Civita connection and the connection that $\Li_B(E)$ inherits from $E$.
\end{enumerate}
\end{definition}
This is equivalent to the definition given by Stolz \cite{sto1}, except that I do not assume the fibers of $E$ to be finitely generated. Several of his conditions are encompassed in my definitions of a bundle of Hilbert \cs-modules and a connection thereon. Because tangent vectors generate the Clifford algebra and $c(\xi)^2=\xi^2$ is the defining relation of the Clifford algebra, $c$ extends to a representation of the Clifford algebra.

If we view the connection $\nabla$ formally as a tangent vector, then we can construct $c(\nabla)$. For any section $\psi\in\Gi(\M,E)$, the connection gives $\nabla\psi \in \Omega^1(\M,E)$; we can think of $c$ as a tangent vector field with coefficients in $\Li_B(E)$; $c(\nabla)\psi$ is the contraction of the tangent vector $c$ with the $1$-form $\nabla\psi$.
\begin{definition}
If $E$ is a $B$-Dirac bundle then the \emph{Dirac operator} $D$ on $L^2(\M,E)$ is the closure of the operator defined on $\Gci(\M,E)$ by $ic(\nabla)$. 
\end{definition}
Note that a Dirac operator is a first order differential operator.

\begin{lem}\label{regularity}
If $\M$ is a complete Riemannian manifold then any Dirac operator is self-adjoint and regular.
\end{lem}
\begin{proof}
Roe gives a proof of this in \cite{roe2} using a generalization of a lemma of Chernoff \cite{che1} which says that  it is sufficient for the corresponding wave equation (on $\R\times\M$) to be solvable. The proof of the latter fact uses the following facts:

 Solvability is true in the classical case (a finitely generated $\co$-Dirac bundle). It is still true if we take the tensor product with a trivial bundle of Hilbert $B$-modules (with trivial connection). It is still true if we change the connection, because that only changes the Dirac operator by a self-adjoint section of $\Li_B(E)$. By local triviality and a partition of unity argument, this extends to the general case.
\end{proof}

The first objective here is to prove that the domain of $D^m$ is the same for any choice of $B$-Dirac bundle structure on a given $B$-$\frac12$-density bundle. I will do this by showing that the graph norm of $D^m$ is, at least locally, equivalent to an order $m$ Sobolev norm. The first step is the case of a trivial $B$-Dirac bundle. 
\begin{lem}\label{Dirac.Sobolev}
Let $E$ be a trivial $B$-Dirac bundle over $\R^n$ or $\T^n$ with trivial connection and a constant Riemannian metric. Let $D$ be the corresponding Dirac operator. The order $m$ Sobolev inner product (given by this flat metric) is the graph inner product for $D^m$. 
\end{lem}
\begin{proof}
Since in this case $D=ic(\nabla)$ with $\nabla$ the trivial connection, the action of $D$ in terms of the Fourier transform (or series) is simply $\widetilde{D\psi}(k) = c(k)\tilde\psi(k)$. The Clifford identity then shows that $\widetilde{D^2\psi}(k) = \norm{k}^2\tilde\psi(k)$. So, because $D^*=D$,
\[
\langle\psi\vert\varphi\rangle_{\dom D^m} 
= \braket\psi{1+D^{2m}}\varphi 
= \braket{\tilde\psi}{1+\norm{k}^{2m}}{\tilde\varphi} 
= \langle\psi\vert\varphi\rangle_m
\mbox.\]
\end{proof}
%: Stray QED

We need a couple of general inequalities for the ensuing computation.
\begin{lem}
For any unbounded operator $\Xi$ on a Hilbert $B$-module,
\beq
\norm\psi_{\dom \Xi} \leq \norm\psi + \norm{\Xi\psi} \leq 2 \norm\psi_{\dom \Xi}
\label{graph.norm.equivalent}
\mbox.\eeq
For any integer $m\geq 1$,
\beq
\norm\psi_{\dom D^{m-1}} \lesssim \norm\psi_{\dom D^m}
\label{Sobolev.embedding}
\mbox.
\eeq
\end{lem}
\begin{proof}
To prove \eqref{graph.norm.equivalent}, note that we can write
\[
\norm\psi_{\dom \Xi} = \Norm{a+b}^{1/2}
\]
and
\[
\norm\psi + \norm{\Xi\psi} = \norm{a}^{1/2} + \norm{b}^{1/2}
\]
where $a:= \langle\psi\vert\psi\rangle$ and $b:= \langle \Xi\psi\vert \Xi\psi\rangle$ are both positive elements of $B$.

Now, 
\[
\norm{a+b} \leq \norm{a} + \norm{b} \leq  \norm{a} + 2\norm{a}^{1/2}\norm{b}^{1/2} + \norm{b} = \left(\norm{a}^{1/2}+\norm{b}^{1/2}\right)^2
\]
gives the first inequality in \eqref{graph.norm.equivalent} (after taking a square root).

Because $a$ and $b$ are positive, $\norm{a},\norm{b}\leq \norm{a+b}$, so
\[
\norm{a}^{1/2} + \norm{b}^{1/2} \leq 2 \norm{a+b}^{1/2}
\]
which is the second inequality in \eqref{graph.norm.equivalent}.

To prove \eqref{Sobolev.embedding}, note that the polynomial $2z^{2m}-z^{2m-2}+1$ is positive for $z\in\R$. As $D$ is self-adjoint, this shows that $2D^{2m}-D^{2m-2}+1 \geq 0$. Equivalently,
\[
2(1+D^{2m}) \geq 1+D^{2m-2}
\mbox,\]
and so
\begin{multline*}
\norm\psi_{\dom D^{m-1}} 
= \Norm{\braket\psi{(1+D^{2m-2})}\psi}^{1/2} \\
\leq \sqrt2  \Norm{\braket\psi{(1+D^{2m})}\psi}^{1/2}
= \sqrt2\norm\psi_{\dom D^m}
\mbox.\end{multline*}
\end{proof}

\begin{lem}\label{Sobolev2}
Let $D$ be the Dirac operator for some $B$-Dirac bundle over a neighborhood of $0\in\R^n$. For any integer $m\geq 0$, there exists a neighborhood $U$ of $0$ such that the graph norm of $D^m$ is equivalent to the order $m$ Sobolev norms on $\Gzi(U,E)$.
\end{lem}
\begin{proof}
The trivial connection and the Clifford map at $0\in \R^n$ give another (trivial) $B$-Dirac structure to $E$. Let $D_0$ be the corresponding Dirac operator
and compute the Sobolev norms with this flat Riemannian metric.

The bound in one direction is easy. Because $D^m$ is of order $m$, Lem.~\ref{Sobolev1} and \eqref{graph.norm.equivalent} give
\[
\norm{\psi}_{\dom D^m} \lesssim \norm\psi + \norm{D^m\psi} \lesssim \norm\psi_m
\mbox.\]

I shall prove the other bound by induction on $m$. For $m=0$, $\norm{\,\cdot\,}_{\dom D^0} = \norm{\,\cdot\,} = \norm{\,\cdot\,}_0$, so the bound is trivially true.
Assume that it is true for $m-1$. That is, that $\norm\psi_{m-1}\lesssim \norm\psi_{\dom D^{m-1}}$ for $\psi\in\Gzi(U_{m-1},E)$ and $U_{m-1}\subset \R^n$ a neighborhood of $0$.

Because $D_0$ is the first order part of $D$ at the origin, $D^m$ is equal to $D_0^m$ plus lower order terms at the origin. In other words, we can decompose $D^m$ as $D^m = D_0^m + \Xi + \Upsilon$ where $\Xi$ is an order $m$ differential operator whose coefficients vanish at the origin and $\Upsilon$ is an order $m-1$ differential operator. 

By Lem.~\ref{Sobolev1}, $\norm{\Xi\psi} \lesssim \norm\psi_m$ for $\psi \in\Gzi(U_{m-1},E)$. In other words, $\Xi$ is a bounded operator from $\Gzi(U_{m-1},E)$ with the norm $\norm{\,\cdot\,}_m$ to $\Gzi(U_{m-1},E)$ with the norm $\norm{\,\cdot\,}$. The norm of $\Xi$ in this sense depends upon the norms of its coefficients. Because the coefficients vanish at the origin, we can arbitrarily decrease the norm of $\Xi$ by using a small enough neighborhood $U_m\subset U_{m-1}$. So, $U_m$ can be chosen small enough that the norm of $\Xi$ is strictly less than $1$. In other words,
\[
\norm\psi_m \lesssim \norm\psi_m - \norm{\Xi\psi}
\mbox,\]
for $\psi \in \Gzi(U_m,E)$.

Now, using Lem.~\ref{Dirac.Sobolev} and \eqref{graph.norm.equivalent},
\begin{align*}
\norm\psi_m 
= \norm\psi_{\dom D_0^m}
&\leq \norm\psi + \norm{D_0^m\psi} \\
&\leq
\norm\psi + \norm{D^m\psi} + \norm{\Xi\psi} + \norm{\Upsilon\psi}
\mbox.\end{align*}
Rearranging this and inserting it into the previous result gives,
\begin{align*}
\norm\psi_m &\lesssim \norm\psi_m - \norm{\Xi\psi} \\
&\quad\leq \norm\psi + \norm{D^m\psi} + \norm{\Upsilon\psi}
\mbox.\end{align*}
Since $\Upsilon$ is of order $m-1$, Lem.~\ref{Sobolev1} gives $\norm{\Upsilon\psi}\lesssim \norm\psi_{m-1}$, and with \eqref{Sobolev.embedding}
\begin{align*}
\norm\psi_m
&\lesssim  \norm\psi + \norm{D^m\psi} + \norm{\Upsilon\psi} \\
&\quad \lesssim \norm\psi_{\dom D^m} + \norm\psi_{m-1}
\mbox.\end{align*}
Finally, by the induction hypothesis and \eqref{Sobolev.embedding}, $\norm\psi_{m-1}\lesssim\norm\psi_{\dom D^{m-1}} \lesssim \norm\psi_{\dom D^m}$. Therefore,
\begin{align*}
\norm\psi_m 
&\lesssim \norm\psi_{\dom D^m} + \norm\psi_{\dom D^{m-1}} \\
&\quad \lesssim \norm\psi_{\dom D^m}
\mbox.\end{align*}
\end{proof}

In the remaining results, I assume that $\M$ is compact. 
\begin{thm}\label{domain.equivalence}
As a Hilbertian $B$-module, $\dom D^m$ is independent of the specific $B$-Dirac structure on $E$. That is, it does not depend on the Riemannian metric, connection, or Clifford map.
\end{thm}
\begin{proof}
Let $D_1$ and $D_2$ be the Dirac operators of two different $B$-Dirac bundle structures for $E$.
Around any point $x\in\M$ we can choose a neighborhood that can be identified with a neighborhood of $0\in\R^n$ ($n=\dim\M$). For sections supported in this neighborhood, we can compare the graph norms for $D_1^m$ and $D_2^m$ with the order $m$ Sobolev norms. By Lem.~\ref{Sobolev2}, there exists a neighborhood of $0$ for which these graph norms are both equivalent to the Sobolev norms, and thus equivalent to each another. Thus the graph norms are equivalent for sections supported in some neighborhood of $x$. The set of such neighborhoods gives an open cover of $\M$. By compactness we have a finite cover $\{U_\alpha\}$ of $\M$. such that for each $\Gzi(U_\alpha,E)$, the graph norms for $D_1^m$ and $D_2^m$ are equivalent.

The fact that the norms are equivalent means precisely that the identifying maps are bounded. Let $\id_\alpha$ be the identity map from $\Gzi(U_\alpha,E)$ with the $D_1^m$ graph norm to $\Gzi(U_\alpha,E)$ with the $D_2^m$ graph norm. Likewise, let $\id$ be the identity map for $\Gi(\M,E)$. Let $\{\chi_\alpha\}$ be a smooth partition of unity subordinate to the cover $\{U_\alpha\}$. So, the identity maps $\id_\alpha$ are all bounded, and we can write the identity map $\id$ as
\[
\id = \sum_\alpha \id\chi_\alpha = \sum_\alpha \id_\alpha \chi_\alpha
\mbox.\]
since this is a finite sum, this shows that $\id$ is bounded. In exactly the same way, the inverse identity map is also bounded, and so the $D_1^m$ and $D_2^m$ graph norms are equivalent on all of $\Gi(\M,E)$. The space $\Gi(\M,E)$ is a common core for $D_1^m$ and $D_2^m$, so by Lem.~\ref{equivalence.equivalence} this shows that $\dom D_1^m$ and $\dom D_2^m$ are equivalent as Hilbertian $B$-modules.
\end{proof}

This is a generalized Sobolev embedding theorem.
\begin{thm}\label{differentiable}
If $m>l+\frac12\dim\M$ then $\dom D^m \subset \Gamma^l(\M,E)$ ($l$ times continuously differentiable sections).
\end{thm}
\begin{proof}
Continuity and differentiability are local properties, so it is sufficient to prove this in any coordinate chart. Thus, let $U\subset\M$ be an open set which is identified with a (bounded) open set in $\R^n$ ($n=\dim\M$). By Thm.~\ref{domain.equivalence}, $\dom D^m$ is independent of the choice of Dirac operator. So, for convenience, we can choose a Riemannian metric on $\M$ in which $U$ is flat and a connection on $E$ which is flat over $U$. From this we see that the domain of $D^m$ over $U$ is equal to the Sobolev space $W^m(U,E)$.

Let $\HH$ be the fiber of $E$ over $U$.
We need to check that the order $m$ Sobolev space, $W^m(U,\HH)$, consists of continuously $l$-differentiable sections. The key step is the $l=0$ case.

Assume that $m>\frac12 n$. The Sobolev space $W^m(U,\HH)$ of $\HH$-valued functions is a Hilbert $B$-module and is simply the tensor product of $\HH$ with the scalar Sobolev space $W^m(U)$. In this way, there is an inner product $W^m(U)\times W^m(U,\HH) \to \HH$ which is continuous in both arguments. The evaluation of a section in $W^m(U,\HH)$ at a point can be rewritten as such an inner product as follows. In terms of the Fourier transform $\tilde\psi(k)$, the evaluation of $\psi\in W^m(U,\HH)$ at $x\in\R^n$ is
\[
\psi(x) = (2\pi)^{-n/2} \int e^{-ik\cdot x} \tilde\psi(k)\, d^nk 
= (2\pi)^{-n/2} \langle e^{ik\cdot x} \vert \tilde\psi\rangle
\mbox.\]
This is equal to a Sobolev inner product,
\[
\psi(x) = \langle\varphi_x\vert\psi\rangle_m
\mbox,\]
where $\tilde\varphi_x(k) = (2\pi)^{-n/2} (1+\norm{k}^{2m})^{-1}$. This means that $\psi : U\to\HH$ factors through $W^m(U)$. It will be continuous if the map $x\mapsto \varphi_x$, $U\to W^m(U)$ is continuous. It is sufficient to check that the pulled back metric function,
\[
\norm{\varphi_x-\varphi_y}_m = \left(\langle\varphi_x\vert\varphi_x\rangle_m + \langle\varphi_y\vert\varphi_y\rangle_m + \langle\varphi_x\vert\varphi_y\rangle_m + \langle\varphi_y\vert\varphi_x\rangle_m\right)^{1/2}
\]
is continuous on $U\times U$. In fact,
\[
\langle\varphi_x\vert\varphi_y\rangle_m = (2\pi)^{-n} \int e^{ik\cdot(x-y)} (1+k^{2m})^{-1}\, d^nk
\]
is uniformly convergent and continuous because the integral $\int_0^\infty r^{n-1}(1+r^{2m})^{-1}\,dr$ is convergent (since $m>\frac12n$).

So we see that $W^m(U,\HH) \subset \C(U,\HH)$ when $l=0$. Now let $l>0$. From the Fourier transform formula for the Sobolev norms, it is apparent that an $l$'th order derivative of an $m$'th order Sobolev section is an $m-l$'th order Sobolev section. This proves the more general result.
\end{proof}

This is a special case of a generalization of the Rellich lemma.
\begin{thm}\label{compact.resolvent}
If the fibers of $E$ are \emph{finitely generated} Hilbert $B$-modules then $D$ has compact resolvent, $(D-i)^{-1}\in\K_B(L^2[\M,E])$.
\end{thm}
\begin{proof}
In general, for $\psi,\varphi\in\dom D$,
\[
\langle\psi\vert\varphi\rangle_{\dom D} 
= \braket\psi{1+D^2}\varphi 
= \braket\psi{(D+i)(D-i)}\varphi
= \langle(D-i)\psi\vert(D-i)\varphi\rangle
\]
because $D$ is self-adjoint. So, 
\beq
\label{D-i}
(D-i) : \dom D \isom L^2(\M,E)
\eeq
is an isomorphism of Hilbert $B$-modules. (It is surjective because $D$ is regular and self-adjoint.)

First, consider the claim in the case that $\M$ is a flat torus, $\T^n$, and $E$ is trivial. This allows us to use Fourier series. Let $K=(1+D^2)^{-1}$.
In terms of the Fourier series, this is just $\widetilde{K\psi}(k) = (1+k^2)^{-1} \tilde\psi(k)$. 

Let $\HH$ be the fiber of $E$. The Hilbert $B$-module $L^2(\T^n,E)=L^2(\T^n)\otimes\HH$ is the direct sum of a copy of $\HH$ for every allowed value of $k$. We have assumed that $\HH$ is a finitely generated Hilbert $B$-module, so the identity operator on $\HH$ is compact. Because  $(1+k^2)^{-1}\in\C_0(\Z^n)$, the operator $K$ is compact, $K\in \K_B(L^2[\T^n,E])$. However,
\[
K = (D+i)^{-1}(D-i)^{-1} = \left[(D-i)^{-1}\right]^*(D-i)^{-1}
\mbox,\]
so this shows that $(D-i)^{-1}\in\K_B(L^2[\T^2,E])$ as well.

In the case of the flat torus, the inclusion $\iota : W^1(\T^n,E)=\dom D \into L^2(\T^n,E)$ is the composition of the isomorphism \eqref{D-i} with the compact operator $(D-i)^{-1}\in\K_B(L^2[\T^n,E])$. Thus this shows that $\iota\in\K_B(\dom D, L^2[\T^n,E])$.

The inclusion $W^1(U,E)\subset L^2(U,E)$ is also compact for sections supported over a subset $U\subset\T^n$, because that inclusion is just a restriction of $\iota$. Any bounded open subset of $\R^n$ can be identified with a subset of $\T^n$, so this property also holds for any such domain.

Now consider a general compact manifold $\M$. Let $\{U_\alpha\}$ be a cover of $\M$ by contractible, open sets. For each $U_\alpha$, we can compare the given Riemannian and $B$-Dirac structures with one in which $U_\alpha$ is flat and $E$ is trivial over $U_\alpha$. The domain of the latter Dirac operator over $U_\alpha$ is the Sobolev space $W^1(U_\alpha,E)$ (by Lem.~\ref{Dirac.Sobolev}) and Thm.~\ref{domain.equivalence} shows that this is Hilbertian equivalent to the domain of the given Dirac operator over $U_\alpha$. We have already seen that the inclusion $W^1(U_\alpha,E)\subset L^2(U_\alpha,E)$ is compact. The identification with the domain of $D$ is a bounded-adjoint\-able operator and the composition of a bounded-adjoint\-able with a compact operator is compact; therefore the inclusion of the domain of $D$ over $U_\alpha$ into $L^2(U_\alpha,E)$ is compact.

Let $\{\chi_\alpha\}$ be a partition of unity subordinate to $\{U_\alpha\}$. Let $\iota :\dom D \into L^2(\M,E)$ and $\iota_\alpha : W^1(U_\alpha,E) \into L^2(U_\alpha,E)$ be the inclusions.  We can write $\iota$ as a finite sum
\[
\iota = \sum_\alpha \iota_\alpha\chi_\alpha 
\mbox,\]
but we have already seen that the $\iota_\alpha$'s are compact and the $\chi_\alpha$'s are bounded-adjoint\-able (by Lem.~\ref{module.bundle2}) so $\iota$ is compact.

As an operator on $L^2(\M,E)$, $(D-i)^{-1}$ is the composition of the isomorphism $(D-i)^{-1}:L^2(\M,E)\isom\dom D$ with the inclusion $\iota : \dom D\into L^2(\M,E)$. This shows that $(D-i)^{-1}\in \K_B(L^2[\M,E])$.
\end{proof}

\subsection{Quantization Bundles}
In this section, I pick up where Section \ref{cocycle} left off. All but one of the results (Lem.~\ref{FieldBundle}) in this section will only require the following hypotheses: $\M$ is a connected, smooth manifold with universal covering, $\pi:\Mt\onto\M$ and fundamental group $\G:=\pi_1(\M)$; $\omega\in\Omega^2(\M)$ is a closed $2$-form; $L^s$ is a line bundle over $\Mt$ with connection $\nabla_s$ and curvature $s\pi^*\omega$ for some $s\in\R$; eq.~\eqref{adef} defines a projective right $\G$-action on $L^s$ with twist cocycle $\ss$; and $\Bs:=\csr(\G,\ss)$. These are not necessarily those constructed constructed in Section \ref{twisting}.

My main goal in this paper is to analyze the covering construction for quantizing $\M$ via $\Mt$. My main tool will be Dolbeault operators, which are particular Dirac operators.  However, Dirac operators on noncompact manifolds are not as well behaved as on compact manifolds. Dirac operators on noncompact manifolds do not generally have compact resolvent. It is preferable to work on $\M$ rather than $\Mt$. This can be accomplished by ``wrapping up'' the bundle $L^s$ into a bundle over $\M$. Since $L^s$ is a vector bundle over $\Mt$ which is itself a principal $\G$-bundle over $\M$. We can sum the fibers of $L^s$ over each point of $\M$ to get a bundle of Hilbert spaces over $\M$. In the terminology of the previous section, this is a bundle of Hilbert $\co$-modules. However, the fibers are infinite-dimensional, not finitely generated, so we still do not suit the hypotheses of Thm.~\ref{compact.resolvent} (compact resolvent). The solution is to form a bundle of Hilbert $\Bs$-modules rather than Hilbert spaces. The fibers of this bundle will actually be finitely generated.

The intention to employ Hilbert $\Bs$-modules is the reason for preferring a twist\-ed \emph{right} action of $\G$ on sections of $L^s$. A Hilbert $\Bs$-module has, in particular, a right projective $\G$-action.
\medskip

The bundle of Hilbert $\Bs$-modules I have in mind will be denoted as $\Ls$. The idea is that the more obvious bundle of Hilbert spaces is the completion, $\Ls\otimes_{\tau_s}\co$. The problem is how to reverse a mistake. Having gone too far and completed to a crude bundle of Hilbert spaces, how do we recover the more delicate underlying bundle of Hilbert $\Bs$-modules. 

If $\HH$ is a Hilbert $\Bs$ module, then using the tracial state $\tau_s$ and the KSGNS construction we can complete it to a Hilbert space, $\HH\otimes_{\tau_s}\co$. Suppose that $\HH_0\subset \HH\otimes_{\tau_s}\co$ is a dense, $\G$-invariant subspace, then it is sometimes possible to reconstruct $\HH$ from the structure of $\HH_0$ as a pre-Hilbert space and projective right $\G$-module. Specifically, this is possible if $\HH_0$ is a subspace of (the image of) $\HH$ in $\HH\otimes_{\tau_s}\co$ and is dense in the Hilbert $\Bs$-module topology.

Two facts are relevant. First, the $\co$-valued inner product (on $\HH\otimes_{\tau_s}\co$) is the trace, $\tau_s\left(\langle\,\cdot\,\vert\,\cdot\,\rangle\right)$, of the $\Bs$-valued inner product (on $\HH$). Second, the algebra is itself a pre-Hilbert space $\Bs\subset\Bs\otimes_{\tau_s}\co$, and $\G$ is an orthonormal basis. Thus, 
\[
b = \sum_{\g\in\G} [\g]\tau_s\left([\g]^*b\right) \in \Bs
\mbox.\] 
Using this identity, we can write the $\Bs$-valued inner product in terms of the $\co$-valued inner product as, 
\beq
\langle\psi\vert\varphi\rangle 
= \sum_{\g\in\G} [\g]\tau_s([\g]^*\langle\psi\vert\varphi\rangle) 
= \sum_{\g\in\G} [\g]\tau_s(\langle\psi\cdot\g\vert\varphi\rangle) 
\mbox.\label{reconstruct}
\eeq
where $\psi\cdot\g$ denotes the projective $\G$ action, as in Section \ref{cocycle}.
\medskip

We can construct a bundle of Hilbert spaces over $\M$ from $L^s\onto\Mt\onto\M$ by taking the Hilbert space direct sum of all fibers of $L^s$ over each point of $\M$.
Suppose that there exists a bundle of Hilbert $\Bs$-modules $\Ls\onto\M$ such that the completion $\Ls\otimes_{\tau_s}\co$ to a bundle of Hilbert spaces is that constructed from $L^s$. The space of sections $\Gz(\Mt,L^s)$ is a Hilbert $\C_0(\Mt)$-module, $\Gz(\M,\Ls)$ is a Hilbert $\C_0(\M)\otimes\Bs$-module, and $\Gz(\M,\Ls\otimes_{\tau_s}\co)$ is a Hilbert $\C_0(\M)$-module. These spaces should all be completions of $\Gc(\Mt,L^s)$ with the respective inner products. The first inner product is essentially pointwise multiplication, so denote it as $\bar\psi\varphi$ for $\psi,\varphi\in\Gc(\Mt,L^s)$. Denote the second inner product with the usual notation as $\langle\psi\vert\varphi\rangle\in\C_0(\M)\otimes\Bs$. I don't use any symbol for the last ($\C_0(\M)$-valued) inner product, but it can be written equivalently in terms of either of the others as,
\[
(\id\otimes\tau_s)\left(\langle\psi\vert\varphi\rangle\right) = \sum_{\g\in\G}\g^*(\bar\psi\varphi)
\mbox,\]
using the identification of $\C_0(\M)$ with the $\G$-invariant functions on $\Mt$ on the right hand side.
Using this and eq.~\eqref{reconstruct} (the additional factor of $\C_0(\M)$ makes no difference in the computation) we can determine what the $\C_0(\M)\otimes\Bs$-valued inner product must be. With the computation,
\begin{align*}
[\g_1\g_2^{-1}]\, \g_2^*\!\left(\overline{\psi\cdot\g_1\g_2^{-1}}\,\varphi\right)
&= [\g_1\g_2^{-1}] \,\overline{\psi\cdot\g_1\g_2^{-1}\cdot\g_2}\,\varphi\cdot\g_2 \\
&= \ss\left(\g_1\g_2^{-1},\g_2\right)^{-1} [\g_1\g_2^{-1}] \,\overline{\psi\cdot\g_1}\,\varphi\cdot\g_2 \\
&= \left([\g_1\g_2^{-1}][\g_2][\g_1]^{-1}\right)^{-1} [\g_1\g_2^{-1}] \,\overline{\psi\cdot\g_1}\,\varphi\cdot\g_2 \\
&= [\g_1][\g_2]^{-1} \,\overline{\psi\cdot\g_1}\,\varphi\cdot\g_2
\mbox,\end{align*}
we have
\beq
\langle\psi\vert\varphi\rangle = \sum_{\g_2,\g_1\in\G} [\g_1][\g_2]^{-1} \, \overline{\psi\cdot\g_1}\,\varphi\cdot\g_2
\mbox.\label{inner}
\eeq
For $\psi,\varphi\in\Gc(\Mt,L^s)$, $\langle\psi\vert\varphi\rangle \in \C_0(\M)\otimes\co[\G,\ss]$ (algebraic tensor product). Two simple consistency checks of this formula are that for $\Mt=\G$ it gives the correct formula for $a^*b\in\co[\G,\ss]$, and because $\tau_s\left([\g_1][\g_2]^{-1}\right) = \delta_{\g_1,\g_2}$,
\[
(\id\otimes\tau_s)\left(\langle\psi\vert\varphi\rangle\right) 
= \sum_{\g\in\G} \overline{\psi\cdot\g}\,\varphi\cdot\g
= \sum_{\g\in\G}\g^*(\bar\psi\varphi)
\mbox.\]

\begin{lem}\label{Lcompletion}
There exists a unique bundle of Hilbert $\Bs$-modules $\Ls$ over $\M$ such that:
\begin{enumerate}
\item
With the inner product \eqref{inner}, $\Gc(\Mt,L^s) \subset \Gz(\M,\Ls)$ and is dense.
\item
The fibers of $\Ls$ are isomorphic to $\Bs$. 
\item
There is a unique, connection on $\Ls$ with curvature $s\omega$ that coincides with the $L^s$ connection on $\Gci(\Mt,L^s)$.
\end{enumerate}
\end{lem}
\begin{proof}
I have already shown that $\Gc(\Mt,L^s)$ is a $\co[\G,\ss]$-module. It is of course a $\C_{\mathrm b}(\Mt)$-module and in particular a $\C_0(\M)$-module. Because $\C_0(\M) \subset \C_{\mathrm b}(\Mt)$ is $\G$-invariant, the $\C_0(\M)$ and $\co[\G,\ss]$ actions commute.

The sum in \eqref{inner} is well defined for $\psi,\varphi\in\Gc(\Mt,L^s)$ because it has finitely many nonzero terms at each point of $\M$. To check $\co[\G,\ss]$-linearity, a little algebra shows that
\[
\langle\psi\vert\varphi\cdot\g\rangle = \langle\psi\vert\varphi\rangle [\g]
\mbox.\]

The universal covering space, $\Mt$, is a principal $\G$ bundle over $\M$. It is thus locally trivial, so any point of $\M$ has an open neighborhood $U$ such that $\pi^{-1}U \cong U\times\G$. Given such a neighborhood, make such an identification.

Since $L^s$ is locally trivial, for a suitable neighborhood, $U$, we can identify $\C_0^\infty(U)$ with the smooth sections of $L^s$ supported over $U\times\{e\}\subset U\times\G = \pi^{-1}U$. With this identification, the $\co[\G,\ss]$-module structure defines a map, $\C_0^\infty(U)\otimes \co[\G,\ss] \to \Gzi(\pi^{-1}U,L^s)$. This map is an isomorphism; both sides have one copy of $\C_0^\infty(U)$ for every point of $\G$. 

This shows local triviality, thus $\Gc(\Mt,L^s) = \Gc(\M,\Ls_0)$ where $\Ls_0$ is a bundle of $\co[\G,\ss]$-modules over $\M$ whose fibers are isomorphic to $\co[\G,\ss]$.
The transition functions between different local trivializations preserve the $\co[\G,\ss]$-valued inner product. If we complete the fibers to Hilbert $\Bs$-modules (in fact, isomorphic to $\Bs$ itself) then these transition functions are unitary. Thus the completed bundle, $\Ls$, is a (smooth)  bundle of Hilbert $\Bs$-modules over $\M$.

The identification of $\Gzi(\pi^{-1}U,L^s)$ with $\C_0^\infty(U)\otimes\co[\G,\ss]\approx \C_0^\infty(U\times\G)$ is just a matter of multiplying by a phase-valued function on $U\times\G$. So, the connection coefficients are still just scalar functions. Compatibility with the projective action of $\G$ implies that the connection coefficients in this trivialization are $\G$-invariant. So, in this local trivialization, the connection coefficients are just scalar functions on $U$. These are certainly bounded-adjoint\-able, so this connection extends to a connection on the bundle of Hilbert $\Bs$-modules, $\Ls$.

Since the original connection on $L^s$ is essentially a restriction of this connection, the curvature must be the same, $s\omega$. 
\end{proof}

Except for the statements about smoothness and the connection, we do not need $\M$ to be a manifold for this result. It is true if $\Mt$ is a principal $\G$-space with compact quotient $\M=\Mt/\G$ and $L^s$ is a bundle over $\Mt$ with a projective $\G$-action on sections. 

When $s=0$, this construction reduces to the \Miscenko-Fomenko line bundle, $\Lb^0 = \Mt\times_\G\csr(\G)$. Indeed, $\Ls$ seems to be the natural twisted generalization of the \Miscenko-Fomenko bundle.
Because the fibers of $\Ls$ are isomorphic to $\Bs$ as Hilbert \cs-modules, $\Ls$ can be heuristically regarded as a line bundle in the same way. In this sense, the obstruction to the existence of a line bundle with curvature $s\omega$ has been circumvented by employing Hilbert \cs-modules. 

This is not the only bundle of this kind.
\begin{cor} 
Let $\Bsm :=\csm(\G,\ss)$.
There exists a unique bundle of Hilbert $\Bsm$-modules $\Lsm$ over $\M$ such that:
\begin{enumerate}
\item
With the inner product \eqref{inner}, $\Gc(\Mt,L^s) \subset \Gz(\M,\Lsm)$ and is dense.
\item
The fibers of $\Lsm$ are isomorphic to $\Bsm$. 
\item
There is a unique, holomorphic connection on $\Lsm$ with curvature $s\omega$ that coincides with the $L^s$ connection on $\Gci(\Mt,L^s)$.
\end{enumerate}
\end{cor}
\begin{proof}
The proof is identical to that of Lem.~\ref{Lcompletion}. We simply regard $\co[\G,\ss]\subset\Bsm$ and thus complete with a norm derived from the maximal (rather than reduced) norm on $\co[\G,\ss]$.
\end{proof}

Both of these bundles are generalizations of the quantization line bundle used in the standard geometric quantization construction. However, $\Ls$ is the more important of these two bundles. This has a lot to do with the existence of the canonical faithful state $\tau_s$ on $\Bs\equiv\csr(\G,\ss)$. It will be important in Thm.~\ref{standard} below which shows the relationship to geometric quantization of $\Mt$. It was also key in proving Thm.~\ref{Bfield} that the $\Bs$'s form a continuous field. There is no reason to expect the $\Bsm$'s to form a continuous field in general. 

The continuous field structure of $\B$ also extends to the $\Ls$'s. For this lemma, we do need the full assumptions of Section \ref{twisting}. In particular: the cocycle $\ss$ is given as an exponential \eqref{ss.def} of the cocycle $c$, and $L^s$ is globally trivial with connection $\nabla_s\equiv d-is A$.

\begin{lem}\label{FieldBundle}
There exists a bundle $\LL$ of Hilbert $\Gz(\R,\B)$-modules over $\M$ such that the push-forward to a bundle of Hilbert $\Bs$-modules is $\Ls$. For any compact interval, $J$, there is a connection on $\Lb^{J}$, the push-forward to a bundle of Hilbert $\Gamma(J,\B)$-modules. These are consistent in the sense that the connection for a subinterval is that induced from the connection for the larger interval. In particular, the connection induced for $s\in J$ is the connection of $\Ls$ 
\end{lem}
\begin{proof}
For clarity in this proof, I will use $\R$ to denote the set of real numbers as a group, and $\hat\R$ for the set of real numbers thought of as the Pontrjagin dual $\hat\R = \Spec[\cs(\R)]$.

The group $\tilde\G$ is again the extension of $\G$ by $\R$ defined in Section \ref{continuous.field}. Recall that the reduced \cs-algebra of $\tilde\G$ is the space of sections, $\csr(\tilde\G)=\Gz(\hat\R,\B)$, of the field of twisted group \cs-algebras. The space $\Mt\times\R$ is a principal $\tilde\G$-bundle over $\M$ with the action defined by, for $(\g,r)\in\tilde\G$ and $(x,r')\in\Mt\times\R$,
\[
(\g,r)(x,r') = \left(\g(x), r+r' + \phi_\g(x)\right)
\mbox.\]
From this, construct the bundle $\Lb:=(\Mt\times\R)\times_{\tilde\G}\csr(\tilde\G)$ over $\M$ of Hilbert $\csr(\tilde\G)$-modules. The space of continuous sections is a completion of $\C_{\mathrm c}^\infty(\Mt\times\R)$ to a Hilbert $\C_0(\M)\otimes\csr(\tilde\G)$-module. By Lem.~\ref{ModuleField}, this is a continuous field of Hilbert \cs-modules associated to $\B$. The identification of $\csr(\tilde\G)$ with $\Gz(\hat\R,\B)$ is a Fourier transform on $\R$. Applying this here shows that the push-forward to a Hilbert $\C_0(\M)\otimes\Bs$-module is $\Gz(\M,\LL)\otimes_{\C_0(\M)\otimes\csr(\tilde\G)}[\C_0(\M)\otimes\Bs] = \Gz(\M,\Ls)$.
 
Now regard the bundles $L^s$ as forming a bundle over $\Mt\times\hat\R$ such that the restriction to $\Mt\times\{s\}$ is $L^s$. The space of compactly supported, smooth sections of this is dense in $\Gz(\M,\LL)$. Define a connection on this by $\nabla = d - isA$, where $s$ is now the coordinate on $\hat\R$.  Clearly, this restricts to the correct connection on each $\Ls$. It does not extend to a connection on $\LL$ because the potential is unbounded. However, for any compact interval,  the restriction is bounded. 
Consider the quotient bundle $\LL\otimes_{\Gz(\hat\R,\B)} \Gz(J,\B)$ for a compact interval $J$. The connection coefficients are bounded and adjointable.
\end{proof}

There is a common structure to all these bundles. There are common properties as well. In order to discuss $\Ls$, $\Lsm$, $\Lb^{J}$, $\Ls\otimes_{\tau_s}\co$, and generalizations on an equal footing, define:
\begin{definition}
\label{Q.bundle.def}
A \emph{quantization bundle}, $Q$, with \emph{coefficient algebra}, $B$, is a bundle of Hilbert $B$-modules over a symplectic manifold, $\M$, with curvature $\zeta\omega$, where $B$ is a unital \cs-algebra, $\zeta\in B$ is self-adjoint and central.
\end{definition}
I shall use the symbols $Q$ and $B$ in this way throughout. The case that $\zeta=s\in\R$ is obviously of particular importance.

A standard quantization line bundle (and its tensor powers) fit this definition with coefficient algebra $B=\co$ and $\zeta$ an integer.
The quantization bundle $\Lsm$ is significant because it is universal in the following sense. 
\begin{thm}\label{universal1}
Let $s\in\R$. If  $Q$ is a quantization bundle with $\zeta=s$, then:
\begin{enumerate}
\item
There does exist a unique line bundle $L^s$ over $\Mt$ with curvature $s\pi^*\omega$ and a projective action of $\G$ on sections.
\item
There exists a Hilbert $\Bsm$-$B$-bimodule, $R$ such that $Q\cong \Lsm\otimes_\Bsm R$ as Hilbert $B$-module bundles with connections. 
\end{enumerate}
As a Hilbert $B$-module, $R$ is isomorphic to the fiber of $Q$.
\end{thm}
\begin{proof}
Because $Q$ has curvature $s\omega$, there is no local obstruction to the existence of a section $\psi$ of $Q$ such that $\nabla\psi = i\sigma\psi$ for \emph{some} $1$-form $\sigma\in\Omega^1(\M)$. Equivalently, the projectivised bundle $PQ$ is locally flat. If we lift this to the universal covering space $\Mt$, then there is no obstruction to the existence of globally constant sections. So, choose some globally constant section of the projectivised bundle over $\Mt$. The fibers of this (with $0$) form a line bundle of $\Mt$. This is a subbundle of $\pi^*Q$ which is preserved by the connection. Therefore it is a line bundle with curvature $s\pi^*\omega$.

This is not necessarily in the purview of Section \ref{cocycle} because the bundle $L^s$ does not necessarily exist for all possible values of $s$, merely for the given value. However, as I remarked in Section \ref{cocycle}, the construction of the cocycle $c$ there with real coefficients would work just as well with any coefficient group. If we apply that construction to the class of $s\omega$ in $H^2(\M,\T)$, then we will get the appropriate twist cocycle $\ss$. Similarly, the projective representation of $\G$ on sections of $L^s$ is obtained.

Consider the lift $\pi^*Q$ of $Q$ to $\Mt$. 
This has curvature $s\,\pi^*\omega$. We can cancel this curvature by taking the tensor product with the dual of $L^s$ to get a flat bundle,
\[
(L^s)^*\otimes \pi^*Q
\mbox.\]
Because $\Mt$ is simply connected, this bundle is globally trivialized by its connection. Thus $(L^s)^*\otimes \pi^*Q \cong \Mt\times R$, for some Hilbert $B$-module, $R$.
We can now write
\[
\pi^*Q = L^s\otimes (R\times\Mt)
\mbox.\]

Because $\pi^*Q$ is the lift of a bundle over $\M$, there is a true (not just projective) left representation of $\G$ on  $\Gc(\Mt,\pi^*Q)=\Gc(\Mt,L^s)\otimes R$; the action of $\g\in\G$ is the pull-back by $\g^{-1}$. Just as we summed up fibers of $L^s\onto\Mt\onto\M$ to get the bundle, $\Ls\otimes_{\tau_s}\co$, of Hilbert spaces, we can sum up fibers of $\pi^*Q$ to get a bundle of Hilbert $B$-modules over $\M$. The space of continuous sections of this is the Hilbert $\C_0(\M)\otimes B$-module $\Gz(\M,\Ls\otimes_{\tau_s}\co)\otimes R$. The representation of $\G$ on this is unitary (in the Hilbert \cs-module sense), therefore there must be a projective unitary representation of $\G$ on $R$; this must have the same twist cocycle $\ss$ to cancel the twist of $L^s$.  The maximal \cs-algebra $\Bsm\equiv\csm(\G,\ss)$ has the universal property that any projective unitary representation of $\G$ with twist $\ss$ must be a unitary representation of $\Bsm$. So, $R$ is a Hilbert $\Bsm$-$B$-bimodule.

There is an obvious surjection, $\Gc(\Mt,\pi^*Q)\onto\Gz(\M,Q)$, given by summing over $\G$-orbits. The kernel of this surjection consists of those sections which add up to $0$, so the image is the space of $\G$-coinvariants, which is the tensor product over $\co[\G,\ss]$. 
Now, because the fibers of $\Lsm$ are isomorphic to $\Bsm$, $\Gz(\M,\Ls_{\max})$ is generated as a $\Bsm$-module by the subspace $\Gc(\Mt,L^s)$, and so
\begin{align*}
\Gz(\M,Q) &= \Gc(\Mt,L^s)\otimes_{\co[\G,\ss]} R \\
&= \Gc(\Mt,L^s)\otimes_{\co[\G,\ss]} \Bsm \otimes_\Bsm R \\
&= \Gz(\M,\Ls_{\max}) \otimes_\Bsm R
\mbox.\end{align*}
Therefore $Q = \Lsm\otimes_\Bsm R$.

The connections must agree because $(L^s)^*\otimes\pi^*Q$ and $(L^s)^*\otimes\pi^*( \Lsm\otimes_\Bsm R) =\Mt\times R$ have the same trivial connection.

By definition, $R$ is the fiber of $(L^s)^*\otimes\pi^*Q$. Since $L^s$ is a line bundle, this isomorphic to any fiber of $\pi^*Q$ and thus to any fiber of $Q$.
\end{proof}

\section{The Toeplitz Construction}
\label{Toeplitz}
For this section, Let  $\M$ be a compact \Kahler\ manifold, $\omega$ its symplectic form, and $Q$ some quantization bundle. Let $2n=\dim\M$. As in Section \ref{classical}, let $L_0$ be a holomorphic line bundle over $\M$ with a $2n$-form valued inner product; that is, a $\co$-$\frac12$-density bundle. Aside from the results in Section \ref{injectivity}, it is not actually necessary to assume that $L_0$ is a line bundle, only that it is of finite rank. 

The combination $Q\otimes L_0$ will now play the preeminent role, so I denote it as $\Qb=Q\otimes L_0$. 
\begin{definition}
\[
\Hi_Q := \Gh(\M,\Qb)
\]
is the space of holomorphic sections of $\Qb$.
\end{definition}
The idea is to mimic the standard Toeplitz quantization construction (Section \ref{classical}) with $Q$ replacing the quantization line bundle $L^{\otimes N}$. The Toeplitz operator of a function $f\in\C(\M)$ will be defined by composing with the orthogonal projection from the Hilbert $B$-module $L^2(\M,\Qb)$ down to the Hilbert $B$-submodule, $\Hi_Q$. Before we can do this, we must first establish that $\Hi_Q$ really is a Hilbert $B$-module and that the orthogonal projection exists.

It will be necessary to lower bound the curvature $\zeta\omega$ of $Q$ as $\zeta\geq s$ for some $s\in\R$. The fundamental results about $\Hi_Q$ are true for $s$ sufficiently large, and the asymptotic results concern the limit of $s\to\infty$.

\subsection{Dolbeault Operators}
There are several ways of viewing this space of holomorphic sections. For one thing, it is the degree $0$ cohomology of the sheaf, $\Or_\M(\Qb)$, of local holomorphic sections. However, for sufficiently large $s$, the higher degree cohomology groups are trivial, so $\Hi_Q$ is the \emph{total} cohomology of $\Or_\M(\Qb)$ (see Lem.~\ref{holomorphic} below). The cohomology of a sheaf of local holomorphic sections can also be computed as the cohomology of a Dolbeault complex. 

In the simplest Dolbeault complex, which computes the cohomology of $\Or_\M$, the space of $p$-chains is $\Omega^{0,p}(\M)$. The coboundary operator is $\bar\partial$, the antiholomorphic component of the exterior derivative operator $d=\partial+\bar\partial$. 

In the Dolbeault complex that computes $H^*[\Or_\M(\Qb)]$ the space of $p$-chains is $\Omega^{0,p}(\M,\Qb)$. I will still denote the coboundary as $\bar\partial$, but it is of course constructed using the connection on $\Qb$.

This Dolbeault complex is an elliptic complex. As such, the cohomology is isomorphic to the kernel of a certain elliptic operator, namely the sum of the coboundary and its adjoint, $\bar\partial + {\bar\partial}^*$. 

\begin{definition}
The Dolbeault operator $D$ is the (unbounded) closed operator on the Hilbert $B$-module $\h_Q := L^2(\M,\Qb\otimes\wedge^{0,*}\M)$  which restricts to $\bar\partial+{\bar\partial}^*$ on the smooth forms $\Omega^{0,*}(\M,\Qb)\subset\h_Q$. 
\end{definition}

So, for $s$ sufficiently large, the space of holomorphic sections $\Hi_Q$ should be realized as the kernel of this Dolbeault operator. I shall prove this below in Lem.~\ref{holomorphic}.
 
\begin{lem}\label{is.Dirac}
$\sqrt2\,D$ is a Dirac operator. 
\end{lem}
\begin{proof}
To deduce the structure of $\Qb\otimes\wedge^{0,*}\M$ as a $B$-Dirac bundle, we need to rewrite $\bar\partial+{\bar\partial}^*$ in the form $\frac{i}{\sqrt2}c(\nabla)$. Any $1$-form on $\M$ is the sum of some $\upsilon\in\Omega^{0,1}(\M)$ and $\xi\in\Omega^{1,0}(\M)$. For $\psi\in\Omega^{0,*}(\M,\Qb)$, we need $\frac{i}{\sqrt2}c(\upsilon)\psi = \upsilon\wedge\psi$. So, $c(\upsilon)\psi := -i\sqrt2\,\upsilon\wedge\psi$. The other part of the Clifford action can be deduced as $c(\xi)= c(\bar\xi)^*$. To be precise, use the metric to identify $\xi$ with a section of $\TM$. The Clifford action is then given by $c(\xi)\psi := i\sqrt2\,\xi\cdot\psi$, where the dot denotes the contraction of a tangent vector and form. Direct computation then shows that $c(\upsilon)^2=0$, $c(\xi)^2=0$ and $[c(\upsilon),c(\xi)]_+=2\xi\cdot\upsilon$, as they should be.

Tensor products of bundles are defined in local trivializations in the usual way. By definition, $Q$ is a Bundle of Hilbert $B$-modules and $L_0$ is a $\co$-$\frac12$-density bundle; the fibers of the tensor product $\Qb\equiv Q\otimes L_0$ inherit $B\otimes\wedge^{2n}\M$-valued inner products. These integrate to give a $B$-valued inner product, thus $\Qb$ is a $B$-$\frac12$-density bundle. The \Kahler\ structure gives a $\co$-valued inner product on each fiber of $\wedge^{0,*}\M$, so the tensor product $\Qb\otimes\wedge^{0,*}\M$ is also a $B$-$\frac12$-density bundle.

The connections on all the factors are compatible with the inner products, so the same is true of $\Qb\otimes\wedge^{0,*}\M$. The $\Z_2$ grading is given be the grading of $\wedge^{0,*}\M$ into even and odd degrees; this is compatible with the connection because $\wedge^{0,*}\M$ is a direct sum of bundle with connections. The action of $c(\upsilon)$ increases degree by $1$ and $c(\xi)$ decreases degree by $1$, so the action of $c$ is odd.

By this construction, $c$ really only deals with $\wedge^{0,*}\M$, so 
\[
c :T^*\M \to \End(\wedge^{0,*}\M) \subset \Li_B(\Qb\otimes\wedge^{0,*}\M)
\mbox.\]
This map is parallel because the wedge product, inner product, and complex structure are parallel with respect to the Riemannian connection.
\end{proof}
Of course, the factor of $\sqrt2$ doesn't effect any of the important properties of Dirac operators discussed in Section \ref{B.bundles}. $D$ is regular, self-adjoint, and has compact resolvent if the fibers of $Q$ are finitely generated.

The main technical lemma we need is a Weitzenbock identity and the resulting lower bound on $D^2$. Let $\db : \Gamma^\infty(\M,\Qb\otimes\wedge^{0,*}\M) \to \Omega^{0,1}(\M,\Qb\otimes\wedge^{0,*}\M)$ be the antiholomorphic gradient operator; this is identical to $\bar\partial$ only in degree $0$. 

\begin{lem}
\label{D2lemma}
Let $K$ be the curvature of $L_0\otimes \wedge^{0,n}\M$,
\beq
(\bar\partial+\bar\partial^*)^2 =  \db^*\db + \zeta\delta + \hat{K}
\mbox,\label{Dsquared}\eeq
where $\hat K$ is a section of $\Li(L_0\otimes \wedge^{0,*}\M)$ constructed $\C(\M)$-linearly from $K$, and $\delta$ is the degree operator for $\wedge^{0,*}\M$. Consequently,
\beq
D^2 \geq s\delta - \norm{\hat K}
\mbox.\label{Dgreater}
\eeq
\end{lem}
Note that $\Khat$ \emph{does not} depend upon $Q$.
\begin{proof}
This is closely related to the computation needed in the proof of the Kodaira-Nakano vanishing theorem.

% I denote the bundle of $(q,p)$-forms as $\wedge^{q,p}\M$.
The canonical volume form of the \Kahler\ structure identifies $\wedge^{n,n}\M = \wedge^{n,0}\M\otimes \wedge^{0,n}\M$ with the trivial bundle. Using this, we can identify $\wedge^{0,*}\M = \wedge^{0,n}\M\otimes \wedge^{n,*}\M$. In this way, regard our Dolbeault complex as $\Omega^{0,*}(\M,\Qb) = \Omega^{n,*}(\M,\Qb\otimes \wedge^{0,n}\M)$.

 Let $\Lambda$ be the adjoint of the operator of exterior multiplication by $\omega$. Using the commutator/anticommutator notation, $[A,B]_\pm = AB\pm BA$,  a standard computation\footnote{See \cite{g-h1}, but ignore the erroneous $\frac12$ there.}  
gives the identity,
\[
[\Lambda,\bar\partial]_- = -i\partial^*
\mbox,\]
and likewise,
\[
[\Lambda,\partial]_- = i{\bar\partial}^*
\mbox.\]
However, $\partial$ of an $(n,*)$-form would be an $(n+1,*)$-form; hence $\partial=0$ on $(n,*)$-forms and the second identity simplifies to,
\[
{\bar\partial}^* = i \partial \Lambda
\mbox.\]
The definition of curvature gives,
\[
(\partial+\bar\partial)^2 = \nabla^2 = -i\zeta\omega-iK
\mbox.\]
Because $L_0\otimes \wedge^{0,n}\M$ is a holomorphic bundle, $K$ is a $(1,1)$-form. Thus ${\bar\partial}^2=({\bar\partial}^*)^2=0$.

We wish to compute 
\begin{align*}
(\bar\partial+{\bar\partial}^*)^2
&=[\bar\partial,{\bar\partial}^*]_+ \\
&= i[\bar\partial,\partial\Lambda]_+
\\&= i \partial[\Lambda,\bar\partial]_- +i [\bar\partial,\partial]_+\Lambda
\\&= \partial \partial^* + (\zeta\omega+K)\Lambda
\mbox.\end{align*}
On $(n,p)$-forms, $\delta:=\omega\Lambda = p$ and $\partial \partial^*= \db^*\db$. This gives eq.~\eqref{Dsquared} with $\Khat=K\Lambda$.

Inequality \eqref{Dgreater} is immediate from \eqref{Dsquared} because smooth forms are dense in the domain $\dom D^2$, $\db^*\db$ is explicitly positive, and $\zeta\geq s$.
\end{proof}
Note that more generally, on $(q,p)$-forms, $[\bar\partial,{\bar\partial}^*]_+ = [\partial, \partial^*]_++ [\zeta\omega+K,\Lambda]_+$ and $[\omega,\Lambda]_- = n-p-q$.

\begin{lem}\label{holomorphic}
For $s>\norm{\Khat}$, the kernel of $D$ is the space of holomorphic sections,
\[
\Hi_Q  = \ker D \subset \Omega^{0,0}(\M,\Qb)
\mbox.\]
\end{lem}
\begin{proof}
Since $D$ is self-adjoint, $\ker D = \ker D^2$. The operator $D^2$ preserves the $\Z$-grading of $\h_Q$. By \eqref{Dgreater}, the restriction of $D^2$ to positive degree is lower bounded by $s-\norm{\Khat}$. So, if $s>\norm{\Khat}$ then $\ker D$ can only be of degree $0$. 

By self-adjointness, $\ker D^m = \ker D$ for any $m>0$. Obviously, the kernel is contained in the domain, so $\ker D \subset \dom D^m$. However, Thm.~\ref{differentiable} shows that $\dom D^m$ consists of $(m-n-1)$-times differentiable sections. Thus $\ker D$ must consist of smooth sections.
The restriction of $D$ to degree $0$ smooth sections is $\bar\partial$. Thus,
\[
\ker D = \ker \bar\partial \equiv \Gh(\M,\Qb) \equiv \Hi_Q
\mbox.\]
\end{proof}

\begin{lem}\label{parametrix}
For $s>\norm{\Khat}$:
\begin{enumerate}
%\item
%If $0\in\Spec D$, it is an isolated point.
\item
There exist a self-adjoint kernel projection and a parametrix, $\Pi$ and $P\in\Li_B(\h_Q)$, such that $\Im\Pi = \Hi_Q$ and $PD=DP=1-\Pi$.
\item
$\norm P \leq (s-\norm\Khat)^{-1/2}$
\item
If the fibers of $Q$ are finitely generated, then $\Pi\in\K_B(\h_Q)$.
\end{enumerate}
\end{lem}
\begin{proof}
Rather than using the full $\Z$-grading of $\h_Q$, it is convenient to just use the $\Z_2$-grading by even and odd degree. In this way, write $D$ as $D=\left(\begin{smallmatrix}0&D_+\\ D_-&0\end{smallmatrix}\right)$, where $D_+$ maps even degree to odd and $D_-$ maps odd to even. The self-adjointness of $D$ means that $D_-=D_+^*$.

The restriction of $D^2$ to odd degree is $D_+D_-$. Since $1$ is the smallest possible odd degree, $\delta\geq1$ when restricted to the odd subspace and \eqref{Dgreater} gives
\beq
D_+D_-\geq s - \norm{\hat K}
\mbox.\label{DBound}
\eeq
So, if $s>\norm{\hat K}$ then there exists an inverse $S = (D_+D_-)^{-1}$ and $\norm{S}\leq (s-\norm{\hat K})^{-1}$. Define,
\[
P := \begin{pmatrix}0&SD_+\\ D_-S&0\end{pmatrix}
\mbox.\]
Direct computation shows that $PD=DP$, $PDP=P$, and $DPD=D$. From this we deduce that $PD$ is an idempotent with the same kernel and image as $D$. 
$P$ is bounded as 
\[
\norm P = \norm{D_-S} = \norm{SD_+}
= \norm{SD_+D_-S}^{1/2} = \norm S^{1/2} \leq (s-\norm{\hat K})^{-1/2}
\mbox.\] 

For $\lambda\in\co$ with $0<\abs\lambda<\norm{P}^{-1}$, the series,
\[
\lambda^{-1}(PD-1) + P + \lambda P^2 + \lambda^2P^3+\dots
\mbox,\]
is norm-convergent and gives the inverse of $D-\lambda$. This shows that $0$ is isolated in the spectrum of $D$; it is at least a distance $\norm P^{-1}\geq (s-\norm{\Khat})^{1/2}$ from any other point. This means that $P$ and $\Pi$ can be constructed by bounded functional calculus from the regular operator $D$. Therefore they are bounded-adjoint\-able, $\Pi,P\in \Li_B(\h_Q)$.

Because $0\in\Spec D$ is isolated, we can construct the kernel projection by a contour integral,
\[
\Pi = \frac1{2\pi i}\oint_\C \frac{d\lambda}{\lambda - D}
\mbox,\]
where the contour encloses $0$ within the open disc of radius $\norm P^{-1}$ around $0$. 

If the fibers of $Q$ are finitely generated, then by Thm.~\ref{compact.resolvent}, $D$ has compact resolvent. This means that $(D-\lambda)^{-1}$ is compact (when bounded) and so $\Pi\in\K_B(\h_Q)$.
\end{proof}

\begin{lem}\label{compact}
$\Hi_Q\subset\h_Q$ is an orthogonally complemented Hilbert $B$-submodule (a direct summand).
If the fibers of $Q$ are finitely generated (as $B$-modules) then for $s>\norm{\Khat}$, $\Hi_Q$ is finitely generated and $\Li_B(\Hi_Q)=\K_B(\Hi_Q)$.
\end{lem}
\begin{proof}
Since $\Hi_Q = \Im \Pi$, it must be closed. So $\Hi_Q$ is a Hilbert $B$-submodule. The image of $1-\Pi$ is its orthogonal complement.

If $s>\norm\Khat$ then by Lem.~\ref{parametrix}, $\Pi\in\K_B(\h_Q)$. Since $\Hi_Q$ is a complemented Hilbert submodule of $\h_Q$, there is a natural inclusion of \cs-algebras $\Li_B(\Hi_Q)\into\Li_B(\h_Q)$ and $\Pi$ is the image of $1$ in this inclusion. This shows that $1\in\K_B(\Hi_Q)$ and so $\Li_B(\Hi_Q)\subset \K_B(\h_Q)$, which proves that $\Li_B(\Hi_Q)=\K_B(\Hi_Q)$. This can only be true if $\Hi_Q$ is finitely generated. 
% All Hilbert \cs-modules are projective.
\end{proof}

In cases when it is given that the fibers of $Q$ are finitely generated, we can speak of $\Li_B(\Hi_Q)$ and $\K_B(\Hi_Q)$ interchangeably. However, I prefer to regard this as the algebra of compact operators. The reason for this is the generalization to noncompact $\M$. The standard Toeplitz quantization of $\C_0(\M)$ gives compact operators on a Hilbert space.

\subsection{Toeplitz Maps}
Because the image of $\Pi$ is $\Hi_Q$, we can regard $\Pi$ as a bounded-adjoint\-able map from $\h_Q$ to $\Hi_Q$; that is, $\Pi\in\Li_B(\h_Q,\Hi_Q)$. This is the main way in which I will use it.

Identify a function $f\in\C(\M)$ with the operator of multiplication by $f$ on $\h_Q$.
\begin{definition}
The Toeplitz operator of a function $f\in\C(\M)$ is,
\[
T_Q(f) := \Pi f : \Hi_Q\to\Hi_Q
\mbox.\]
\end{definition}
\begin{lem}
This defines a unital, completely positive linear map,
\[
T_Q:\C(\M)\to\Li_B(\Hi_Q)
\mbox.\]
\end{lem}
\begin{proof}
By Lem.~\ref{module.bundle2}, (multiplication by) $f$ is a bounded-adjoint\-able operator on $\h_Q$; by Lem.~\ref{parametrix}, $\Hi_Q\subset\h_Q$ is direct summand, so $f$ restricts to a bounded-adjoint\-able map from $\Hi_Q$ to $\h_Q$. By Lem.~\ref{parametrix}, $\Pi$ is bounded-adjoint\-able. A composition of bounded-adjoint\-able maps is bounded-adjoint\-able. A map of this form is automatically completely positive and unital.
\end{proof}
\begin{definition}
\emph{The Toeplitz algebra} is the \cs-subalgebra, $\A_Q\subseteq\Li_B(\Hi_Q)$, generated by the image of $T_Q$.
\end{definition}

\begin{definition}
The special cases of these constructions for $Q=\Ls$ and $\Lsm$ will be denoted as $\hs$, $\Hs$, $T_s$, $\As$, $\Hsm$, $T_s^{\max}$, and $\Asm$ respectively.
\end{definition}

\subsection{Asymptotic Properties}
\label{asymptotic}
The most important property of the Toeplitz maps is asymptotic multiplicativity; the product $T_Q(f)T_Q(g)$ is approximately $T_Q(fg)$ for large $s$. We can actually bound the difference by a function of $s$, independent of a specific choice of $B$ and $Q$. In order to prove quantization, we need to compute the first order correction to $T_Q(f)T_Q(g)$. However, first consider the much simpler proof of asymptotic multiplicativity.

The symbol $\Or^{-p}(s)$ denotes some expression of order $-p$. That is, something bounded by $Cs^{-p}$ for large $s$, with $C$ a constant independent of $Q$.

\begin{lem}\label{0th.order}
    For $f\in\C^1(\M)$ and $g\in\C(\M)$,
    \[
    \Norm{T_Q(f)T_Q(g)-T_Q(fg)} = \Or^{-1/2}(s)
    \mbox.\]
\end{lem}
\begin{proof}
We can assume that $s>\norm\Khat$.
For any $\psi,\varphi\in\Hi_Q$, the identity $DP=1-\Pi$ gives
\begin{align*}
\langle\psi\rvert T_Q(f)T_Q(g)-T_Q(fg) \lvert\varphi\rangle 
&= \langle\psi\rvert (f\Pi g - fg) \lvert\varphi\rangle\\
&= -\langle\psi\rvert fDPg \lvert\varphi\rangle \\
&= \langle\psi\rvert [D,f]_-Pg \lvert\varphi\rangle
\mbox.\end{align*}
Since $\sqrt2\,D$ is  a Dirac operator, $\norm{[D,f]_-}=\frac1{\sqrt2}\norm{\nabla f}$. Using this, we can bound the norm of this ``multiplicative error'' as
\begin{align*}
\Norm{T_Q(f)T_Q(g)-T_Q(fg)} 
&\leq 
\Norm{[D,f]_-}\norm{P}\,\norm{g} \\
&\leq \tfrac1{\sqrt2}\norm{\nabla f} \left(s-\norm{\Khat}\right)^{-1/2}\norm{g}
\mbox.\end{align*}
This bound is of order $\Or^{-1/2}(s)$.
\end{proof}

Now restrict to the case that $Q$ has curvature exactly $s\omega$.
In order to prove that Toeplitz maps generate a quantization, we need to compute the first order correction to the product $T_Q(f)T_Q(g)$ as $s\to\infty$,  in order to check that this corresponds to the Poisson bracket constructed from the symplectic structure.
\begin{lem}\label{1st.order}
For $f,g\in\C^\infty(\M)$,
\beq
\Norm{T_Q(f)T_Q(g)-T_Q\left(fg+s^{-1}\dbb f\cdot \db g\right)} = \Or^{-2}(s)
\mbox,\label{1st.order.eq}
\eeq
where $\dbb f\cdot \db g$ is the metric inner product of the holomorphic and antiholomorphic gradients of $f$ and $g$.
\end{lem}
\begin{proof}
We can again assume $s>\norm{\Khat}$, so that by Lem.~\ref{holomorphic}, 
\[
\Hi_Q=\ker D \subset\Omega^{0,0}(\M,\Qb)
\mbox.\]
Recall that $\h_Q$ inherits the $\Z$-grading from $\wedge^{0,*}\M$. Again, let $P\in\Li_B(\h_Q)$ be the parametrix from Lem.~\ref{parametrix}. Let $\db : \Gamma^\infty(\M,E) \to \Omega^{0,1}(\M,E)$ denote the antiholomorphic gradient (as opposed to exterior derivative) operator for any vector bundle and $\db^*$ its adjoint. 

For any $\psi,\varphi\in\Hi_Q$,
\begin{align*}
\langle\psi\rvert T_Q(f)T_Q(g)-T_Q(fg) \lvert\varphi\rangle
&= \langle\psi\rvert f[\Pi-1]g \lvert\varphi\rangle
= -\langle\psi\rvert fDP^2Dg \lvert\varphi\rangle \\
&= -\langle\psi\rvert f \bar\partial^* P^2 \bar\partial g \lvert\varphi\rangle
\mbox,
\end{align*}
using in the last step that $\psi$ and $\varphi$ are of degree $0$.
The $P^2$ here acts only on degree $1$. Lemma \ref{parametrix} shows that the restriction of $D^2$ to degree $1$ is invertible (since $1$ is odd). The inverse is the restriction of $P^2$. The formula \eqref{Dsquared} gives the restriction of  $D^2$ to smooth forms. This shows that the restriction of $P^2$ to $\Omega^{0,1}(\M,\Qb)$ is
\beq
P^2 = \frac1{D^2} = \frac1{\db^*\db+s+\Khat} 
= \frac1s - \frac1s \frac{\db^*\db+\Khat}{\db^*\db+s+\Khat} 
\mbox.\label{Psquared}
\eeq
Because $\psi,\varphi\in\Hi_Q$ are holomorphic, $\bar\partial\psi=\bar\partial\varphi=0$ and so
\beq
\braket{\psi}{f \bar\partial^*\bar\partial g}{\varphi} = -\braket{\psi}{(\dbb f)\cdot (\db g)}{\varphi}
\mbox.\label{derivative}\eeq
So, the $\frac1s$ in eq.~\eqref{Psquared} provides the desired first order term in \eqref{1st.order.eq}. With this, we can express an arbitrary ``matrix element'' of the remainder term as,
\begin{multline}
\braket{\psi}{T_Q(f)T_Q(g)
-T_Q\left(fg+s^{-1}\dbb f\cdot \db g\right)}{\varphi}
\\
= \braket{\psi}{f(1-\Pi + s^{-1}\db^*\db)g}{\varphi}
= \frac1s \braket{\psi}{f \bar\partial^*
\frac{\db^*\db+\Khat}{\db^*\db+s+\Khat} \bar\partial g}{\varphi}
\mbox.\label{remainder}
\end{multline}
The numerator and denominator of the fraction in the right hand side commute, thus there is no ambiguity in writing it as a fraction. This fraction is a bounded operator on $\Omega^{0,1}(\M,\Qb)$. Indeed, because $\bar\partial g \lvert\varphi\rangle = (\db g) \lvert\varphi\rangle$ and $\db g$ is bounded (likewise with $f$), this shows that \eqref{remainder} is bounded by a term of order $s^{-1}$\/.

However, this is not good enough. In order to get a better bound we need to rewrite the fraction as an expression in which $\db^*$ is to the right and $\db$ to the left of bounded terms. If we write the fraction as 
\[
(\db^*\db+\Khat)(\db^*\db+s+\Khat)^{-1}
\mbox,\]
then we need to move $\db$ past (the reciprocal of) the denominator. Note that in so doing, the denominator goes from acting on sections of $\Qb\otimes\wedge^{0,1}\M$ to acting on sections of $\Qb\otimes(\wedge^{0,1}\M)^{\otimes 2}$.
This requires computing the commutator of $\db$ with the denominator, $\db^*\db+s+\Khat$. 

First consider the commutator of $\db$ with the Laplace-like term $\db^*\db$. It is necessary to resort to index notation in which $\db$ is $\nabla_\ib$ and $\db^*$ is contraction with $-\nabla^\ib$. Because the bundles concerned are holomorphic, $\nabla_\ib\nabla_\jb=\nabla_\jb\nabla_\ib$ and
\[
[\nabla_\ib,\nabla^\jb\nabla_\jb]_- 
= [\nabla_\ib,\nabla^\jb]_-\nabla_\jb
\mbox.\]
So, by a computation similar to Lem.~\ref{D2lemma} we can write $[\db,\db^*\db]_- = (s+\Khat'+\hat R)\db$, where $\Khat'$ (like $\Khat$) is constructed from $K$, and $\hat R$ is a Riemann curvature term.

The commutator of $\db$ with $\Khat$ is simpler. It is just the antiholomorphic derivative of $\Khat$ and is a $0$'th order operator. Write this as $\bar\partial\Khat := [\db,\Khat]_-$.

Putting this together and writing $K_2:=\Khat+\Khat'+\hat R$, we have
\[
\db(\db^*\db+s+\Khat)=(\db^*\db+2s+K_2)\db + \bar\partial\Khat
\mbox.\]
Dividing by $\db^*\db+2s+K_2$ on the left and $\db^*\db+s+\Khat$ on the right gives,
\begin{multline*}
\db (\db^*\db+s+\Khat)^{-1} \\
= (\db^*\db+2s+K_2)^{-1}\db 
- (\db^*\db+2s+K_2)^{-1} \bar\partial\Khat (\db^*\db+s+\Khat)^{-1}
\end{multline*}
So,
\begin{multline*}
\frac{\db^*\db+\Khat}{\db^*\db+s+\Khat} 
= \db^*(\db^*\db+2s+K_2)^{-1}\db  + \Khat (\db^*\db+s+\Khat)^{-1}\\ 
 - \db^*(\db^*\db+2s+K_2)^{-1} \bar\partial\Khat (\db^*\db+s+\Khat)^{-1} 
\mbox.\end{multline*}

By the same principle as in eq.~\eqref{derivative}, $\db\bar\partial\,f\lvert\varphi\rangle = (\db^2 f) \lvert\varphi\rangle$.
Inserting this last expression for the fraction into eq.~\eqref{remainder} then gives a long expression for a matrix element of the remainder. This shows that,
\begin{multline*}
\Norm{T_Q(f)T_Q(g)-T_Q\left(fg+s^{-1}\dbb f\cdot \db g\right)} \\
\leq \frac{\norm{\dbb^2 f} \norm{\db^2 g}} {s 
\left(2s-\norm{K_2}\right)}
+ \frac{\norm{\dbb f}\norm{\db g}\norm{\Khat} }{s
\left(s-\norm{\Khat}\right)} 
+  \frac{\norm{\dbb^2 f}\norm{\db g} 
\norm{\bar\partial\Khat}}{s \left(s-\norm{\Khat}\right) \left(2s-\norm{K_2}\right)} 
\end{multline*}
The first two terms are of order $\mathcal O^{-2}(s)$, the last term 
is of order $\mathcal O^{-3}(s)$.
\end{proof}

\begin{cor}\label{T.commutator}
For $f,g\in\C^2(\M)$,
\beq
\Norm{[T_Q(f),T_Q(g)]_- - is^{-1}T_Q\left(\{f,g\}\right)} = \Or^{-2}(s)
\mbox.\label{commutator.eq}
\eeq
\end{cor}
\begin{proof}
This is a matter of computing $\dbb f\cdot \db g - \db f\cdot \dbb g$. Using the identities relating the complex structure, metric, and symplectic form, along with the definition \eqref{Poisson} of the Poisson bracket, we can compute this to be $-i\{f,g\}$.
\end{proof}

\section{Topological Properties}
\label{topological}
In this section, I restrict to the case that the fibers of $Q$ are finitely generated, $\M$ is a compact \Kahler\ manifold, the coefficient algebra $B$ is separable, and $s>\norm\Khat$. With these assumptions, Lemmas \ref{holomorphic} and \ref{compact} show that $\Hi_Q=\ker D$ and is a finitely generated Hilbert $B$-module. Any Hilbert \cs-module is projective, so $\Hi_Q$ is a finitely generated, projective $B$-module. We can thus consider its $K$-theory class, $[\Hi_Q]\in K_0(B)$. 

Because $\ker D$ is finitely generated, the Dolbeault operator, $D$, is an (un\-bound\-ed) Fredholm operator. It defines a class in bivariant $K$-theory, $[D] \in   KK_0(\C(\M),B)$. There is a natural index map,
\[
\ind : KK_0(\C(\M),B) \to K_0(B)
\mbox.\]
This is given by the formal difference of the even and odd graded parts of the kernel. However, by Lem.~\ref{holomorphic}, $\Hi_Q=\ker D$ is entirely of degree $0$, thus 
\[
[\Hi_Q] = \ind [D] 
\mbox.\]
This has a couple of significant implications. 

The $K$-theory class specifies $\Hi_Q$ modulo stable equivalence. Since it is given by an index, it is essentially topological. If we deform the geometry of $\M$, then as long as $s>\norm\Khat$, $[\Hi_Q]$ will be unchanged.

Because, $[\Hi_Q]$ is given by an index, there is some possibility of computing it through topological methods. With Thm.~\ref{Morita2} below, this means that the Toeplitz algebra $\A_Q$ is essentially determined by topology.

Because $\Hi_Q$ is a finitely generated Hilbert $\A_Q$-$B$-bimodule, it gives a  cycle in $KK_0(\A_Q,B)$ with trivial Fredholm operator.

\subsection{Twisted Baum-Connes}
If $\G$ is a countable, discrete group, $A$ is a  $\G$-\cs-algebra, and $\sigma:\G\times\G\to\T$ is a $2$-cocycle, then Packer and Raeburn \cite{p-r} define a twisted crossed product \cs-algebra, $A\rtimes_\sigma\G$ analogous to the (maximal) crossed product \cs-algebra $A\rtimes\G$. This is constructed with the universal property that any representation of $A$ which is equivariant with a $\sigma$-twisted projective-unitary representation of $\G$, is a representation of $A\rtimes_\sigma\G$. The twisted group \cs-algebra is the special case $\cs(\G,\sigma) := \co\rtimes_\sigma\G$. There is also a reduced version, $A\rtimes_{\sigma,\mathrm r}\G$, which generalizes the reduced crossed product.

For two  $\G$-\cs-algebras $A$ and $B$, the twisted descent (or induction) map, 
\[
j_{\G,\sigma} : KK^*_\G(A,B) \to KK^*(A \rtimes_\sigma \G,B \rtimes_\sigma\G)
\mbox,\]
was constructed by Baum and Connes \cite{b-c}. There is also a reduced version, defined identically, but with reduced crossed products.

If $X$ is a topological space with a proper $\G$-action and compact quotient, $X/\G$\/, then we can set $A=\C_0(X)$ and $B=\co$. The equivariant $K$-homology of $X$ is defined as $K_*^\G(X) := KK^*_\G(\C_0(X),\co)$, so the twisted descent map is
\[
j_{\G,\sigma} : K_*^\G(X) \to KK^*(\C_0(X)\rtimes_\sigma\G,\cs(\G,\sigma))
\mbox.\]
The (unreduced) classical Baum-Connes assembly map is a composition,
\[
\mu : K_*^\G(X) \xrightarrow{j_\G} KK^*(\C_0(X)\rtimes\G,\cs(\G)) \xrightarrow{[e]\otimes_{\C_0(X)\rtimes\G}}
 K_*[\cs(\G)]
\mbox,\]
where the latter map is the Kasparov product with $[e]\in K_0[\C_0(X)\rtimes\G]$, a canonical class in the $K$-theory. The $K$-theory of $\C_0(X)\rtimes\G$ is the Grothendieck group of $\G$-equivariant bundles over $X$; $[e]$ is the class of the trivial line bundle with the trivial $\G$-action.

To twist this, we need a generalization of $[e]$, but there is no canonical choice in general. The $K$-theory, $K_0[\C_0(X)\rtimes_\sigma\G]$ is the Grothendieck group of vector bundles over $X$ with $\sigma$-twisted projective right $\G$-actions (see \cite{b-c}). There is no canonical choice of such a line bundle. However, any two choices will be related by taking the tensor product with a line bundle with trivial $\G$-action. This ambiguity does not make any difference to the properties of the resulting assembly map.

However, if we specialize to the case at hand, in which the cocycle is of the exponential form, $\ss := e^{isc}$, then we can make a natural choice for any $s\in\R$ by continuously deforming from $[e]$. The result is a topologically trivial line bundle over $X$ with a projective $\G$-action just like that constructed for $L^s$ in Section \ref{cocycle}. With this in mind, call this class $[L^s]\in K_0[\C_0(X)\rtimes_\sigma\G]$ and define the twisted Baum-Connes assembly map as
\[
\mu_s : K_*^\G(X) \xrightarrow{j_{\G,\ss}} KK^*(\C_0(X)\rtimes_\ss\G,\Bsm) \xrightarrow{[L^s]\otimes_{\C_0(X)\rtimes_\ss\G}}
 K_*(\Bsm)
\mbox,\]
where once again, $\Bsm\equiv\cs(\G,\ss)$.
The reduced version is defined identically, with reduced crossed products and $\Bs$.

Twisted Baum-Connes maps can be constructed even when the cocycle is not an exponential; however this requires the $K$-homology on the left hand side to be twisted using the Dixmier-Douady invariant of the cocycle.  The twisted Baum Connes conjecture \cite{b-c} asserts that the reduced assembly map should be an isomorphism for $X=\Ebar\G$, the universal space for proper $\G$-actions.
\medskip

Now consider the case $X=\Mt$ and $\G\equiv \pi_1(\M)$. Any element of $K_0(\M)$ is given by a differential operator and can be lifted to a $\G$-invariant differential operator on $\Mt$. This gives a canonical isomorphism $K_0(\M)\cong K_0^\G(\Mt)$. Let $\varepsilon \in K_0(\M)$ be the \Kahler\ $K$-orientation; that is, the class of the Dolbeault operator for the complex $\Omega^{0,*}(\M)$. Let $\varepsilon\otimes L_0 \in K_0(\M)$ be the class of the ($L_0$-twisted) Dolbeault operator for the complex $\Omega^{0,*}(\M,L_0)$. We can apply the twisted assembly map to $\varepsilon\otimes L_0$, $\mu_s(\varepsilon\otimes L_0)$ is a Kasparov product of $[L^s]$ with $j_{\G,\ss}(\varepsilon\otimes L_0)$. This is also the index of $j_{\G,\ss}(\varepsilon\otimes L_0)$ twisted by $L^s$. So, we lift the $L_0$-twisted Dolbeault operator to $\Mt$, twist by $L^s$ and descend to a class in $KK^0(\C(\M),\Bsm)$. This is just $[D_{\Lsm}]$ and the index is $[\Hsm]$. So, $\mu_s(\varepsilon\otimes L_0) = [\Hsm]$.

Because the universal covering space $\Mt$ is a proper $\G$-space, there is an equivariant  classifying map $\kbar : \Mt \to \Ebar\G$. The assembly map for $\M$ factors as,
\[
\mu_{s} : K_0(\M) \xrightarrow{\kbar_*} K^\G_0(\Ebar\G) \xrightarrow{\mu_{s}} K_0(\Bs)
\mbox.\]
So, $[\Hsm] = \mu_s(\kbar_*(\varepsilon\otimes L_0))$.

By Thm.~\ref{universal1}, if $Q$ is a quantization bundle with curvature $s\omega$ and coefficient algebra $B$, then there exists a Hilbert $\Bsm$-$B$-bimodule, $R$, such that $Q=\Lsm\otimes_\Bsm R$. Provided that the fibers of $Q$ are finitely generated, $R$ will be finitely generated and defines a Kasparov bimodule with trivial Fredholm operator and class $[R]\in KK^0(\Bsm,B)$. The tensor product with $R$ corresponds to the Kasparov product with $[R]$. So, we have a general formula for $[\Hi_Q]$,
\begin{align}
[\Hi_Q] &= [\Hsm\otimes_\Bsm R] = [\Hsm] \otimes_\Bsm [R] \nonumber\\
&= \mu_s(\kbar_*(\varepsilon\otimes L_0)) \otimes_\Bsm [R]
\label{BC}
\mbox.\end{align}

In the special case that $Q=\Ls$, this is given by the reduced twisted assembly map,
\beq
\label{BCr}
[\Hs] = \mu_{s,\mathrm{r}}(\kbar_*(\varepsilon\otimes L_0))
\mbox.\eeq
If the twisted Baum-Connes conjecture is true, then this is equivalent to $\kbar_*(\varepsilon\otimes L_0) \in K^\G_0(\Ebar\G)$.

If $\G$ is torsion free, then $\Ebar\G= E\G$ and this simplifies to $k_*(\varepsilon\otimes L_0) \in K_0(B\G)$. If $\Mt$ is contractible then $\G$ must be torsion free and $B\G \cong \M$; the conjecture reduces to $\mu_{s,\mathrm r} : K_0(\M) \isom K_0(\Bs)$. 

We should expect that as $s$ increases, the space of holomorphic sections $\Hs$ will grow larger. It thus seems paradoxical that the classes $[\Hs]$ are all given by a fixed element of $K^\G_0(\Ebar\G)$. This happens because the twisted assembly maps are not quite canonical. For example, if $[\omega]$ is integral, then $\Bs$ depends periodically on $s$; for $s\in\Z$, $\Bs=\csr(\G)$. There is a sort of holonomy to the twisted assembly map; as we follow $\mu_{s,\mathrm r}(\varepsilon\otimes L_0)$ around continuously, each time we return to $K_0(\csr(\G))$ we get a different element.
\medskip 

I discuss the possible relevance of my construction to the untwisted Baum-Connes conjecture in the conclusions (Section \ref{conclusions}).

\subsection{Trace}
\label{trace.index}
Let $D$ be the Dolbeault operator constructed from $L_0$ and $\Ls$, whose kernel is $\Hs$. The trace $\tau_s$ on the coefficient algebra $\Bs\equiv \csr(\G,\ss)$ extends canonically to matrices over $\Bs$. Since any class in $K_0(\Bs)$ can be written as a formal difference of self-adjoint idempotent matrices, this gives a map,
\[
[\tau_s] : K_0(\Bs) \to \R
\mbox.\]
The twisted $L^2$-index theorem (see \cite{mat1,gro1}) provides a way to compute $[\tau_s]\circ\ind$ topologically, thus giving the invariant $[\tau_s][\Hs]$ of $[\Hs]\in K_0(\Bs)$.

This twisted $L^2$-index is a generalization of Atiyah's $L^2$-index. The theorem applies to a $\G$-invariant Dirac-type operator when the $\G$ action is projective and the twist is the exponential of a real cocycle $sc$. The theorem gives the index as an integral (over the quotient space) of a product of characteristic classes identical to that in the Atiyah-Singer index theorem, but for an additional factor,\footnote{My factors of $2\pi$ are placed differently than in \cite{mat1}} $e^{s\omega/2\pi}$. This $\omega$ is a closed $2$-form whose cohomology class is the pullback of $[c]\in H^2(\G;\R)$ by the classifying map $k:\M\to B\G$. In the case at hand, $c$ is the cocycle constructed in Section \ref{cocycle}, and $\omega$ can be taken to be the symplectic form. The twisted $L^2$-index theorem gives,
\beq\label{L2.index}
[\tau_s][\Hs] = [\tau_s](\ind D) = \intM \td(\TM)\wedge e^{\frac{s\omega}{2\pi} - \frac12 c_1(L_0)}
\mbox.\eeq

This can be proven by adapting the heat kernel proof of the Atiyah-Singer index theorem, almost \emph{mutatis mutandem}. The local trace is replaced with one constructed using $\tau_s$. Equivalently, one can use the heat kernel for the $L^s\otimes\pi^*L_0$-twisted Dolbeault operator on $\Mt$. This is a $\G$-invariant integral kernel on $\Mt$. Its trace is the integral of the trace of the diagonal over a fundamental domain of the $\G$ action. The expression for the index is thus formally identical to that for the classical index of a Dolbeault operator over $\M$ twisted by $L_0$ and a line bundle with curvature $s\omega$. 

\section{Morita Equivalence}
\label{Morita}
There are several equivalent definitions of a Morita equivalence of \cs-algebras. This is the definition given in \cite{con1}.
\begin{definition}
Two \cs-algebras $A$ and $B$ are Morita equivalent if there exist a Hilbert $A$-$B$-bimodule $R_1$ and a Hilbert $B$-$A$-bimodule $R_2$ such that,
\beq
\label{left}
R_1\otimes_B R_2 \cong A
\eeq
and
\beq
\label{right}
R_2 \otimes_A R_1 \cong B
\eeq
as Hilbert \cs-bimodules.
\end{definition}

This looks very much like the conditions for isomorphism in some category with bimodules as morphisms.

The bivariant $K$-theory functor $KK$ can be viewed in terms of an abelian category ``$KK$''\/. The objects in $KK$ are the separable \cs-algebras. However, the morphisms are not the $*$-homomorphisms of the usual category of \cs-algebras. The group $KK(A,B)$ is the set of $KK$-morphisms from $A$ to $B$. A $KK$-morphism from $A$ to $B$ is an equivalence class of Kasparov $A$-$B$-bimodules. A Kasparov $A$-$B$-bimodule is a Hilbert $A$-$B$-bimodule equipped with a grading and a bounded ``Fredholm'' operator. I take this as a guideline in defining a category for studying Morita equivalence. 
Instead of using Kasparov bimodules, simply use Hilbert bimodules. This will not give an abelian category, so it does not define a bivariant homological theory.

Define a putative category, \bimod\ as follows.
\begin{definition}
The objects of \bimod\ are \cs-algebras. The \bimod-morphisms $A\to B$ are the Hilbert $A$-$B$-bimodules, modulo the relation $\approx$ generated by isomorphism and $R\approx A\otimes_A R$. The composition, $R_2\circ R_1:A\to C$, of morphisms $R_1 : A\to B$ and $R_2 :B \to C$ is (the class of) the tensor product $R_1\otimes_B R_2$.
\end{definition}
Note that the composition $R_2\circ R_1$ and the tensor product $R_1\otimes_B R_2$ are written in the opposite order, just as in $KK$.

\begin{thm}
\bimod\ is a category. The identity morphism at $A$ is given by $A$ as a Hilbert $A$-bimodule.
\end{thm}
\begin{proof}
We need to check that the composition of morphisms is well defined and associative, and respects the identity morphisms.
Note that any Hilbert $A$-$B$-bi\-mod\-ule can be reduced to one such that $A\otimes R \cong R$, because $A\otimes R$ has this property.

To check that the composition $R_2\circ R_1$ is well defined, replace $R_1$ and $R_2$ with their reduced equivalents,
\[
(A\otimes_AR_1)\otimes_B(B\otimes_BR_2) 
= A\otimes_AR_1\otimes_BB\otimes_BR_2
\cong A \otimes_A (R_1\otimes_B R_2)
\approx R_1\otimes_B R_2
\mbox.\]

That composition respects identity morphisms is almost dictated by the definition. In fact, for $R:A\to B$.
\[
A\otimes_A R \approx R
\]
and
\[
R\otimes_BB \cong R
\mbox.\]

The associativity of composition is immediate from the associativity of the tensor product.
\end{proof}

\begin{thm}
Morita equivalence is the relation of isomorphism in \bimod.
\end{thm}
\begin{proof}
First, suppose that there exist bimodules $R_1$ and $R_2$ satisfying the conditions of a Morita equivalence. Then this gives directly that $R_2\circ R_1 \approx A$ and $R_1\circ R_2\approx B$, so the morphisms given by $R_1$ and $R_2$ are inverse to one another and are isomorphisms.

Conversely, suppose that $R_1$ and $R_2$ give a pair of inverse isomorphisms. 
We have $A\approx R_1\otimes_AR_2$, so
\[
A \cong A\otimes_A(R_1\otimes_AR_2) \cong 
(A\otimes_AR_1)\otimes_B(B\otimes_BR_2) 
\]
and likewise in the other order. Thus the bimodules $A\otimes_AR_1$ and $B\otimes_BR_2$ give a Morita equivalence.
\end{proof}

Because \bimod\ is a category, there is the usual property of uniqueness of inverses. We can thus think of the conditions for Morita equivalence as two separate conditions on $R_1$; \eqref{left} is left invertibility and \eqref{right} is right invertibility. Left invertibility is analogous to injectivity and right invertibility is analogous to surjectivity.
% By analogy with the category of sets, the first condition can be thought of as injectivity and the second as surjectivity.

Morita equivalence can be defined by a different pair of conditions. 
\begin{definition}
A Hilbert $B$-module, $\HH$ is \emph{full} if the \cs-ideal generated by inner products from $\HH$ is $B$ itself.
\end{definition}
A Hilbert $A$-$B$-bimodule $R$ gives a Morita equivalence (i.~e., is a \bimod\ isomorphism) if it is a full Hilbert $B$-module and the representation of $A$ on $R$ gives an isomorphism $A\cong\K_B(R)$. The latter two conditions are \emph{not} separately equivalent to the two conditions of left and right invertibility on \bimod. 
These are thus two different ways of dividing Morita equivalence into two separate conditions. There is the following relationship.
\begin{prop}
\label{correspondence}
Let $R$ be a Hilbert $A$-$B$-bimodule. If $A=\K_B(R)$ then $R$ determines a left invertible \bimod\ morphism. If $R$ determines a right invertible \bimod\ morphism, then $R$ is a full Hilbert $B$-module.
\end{prop}
\begin{proof}
In general, $\K_B(R,B)$ is a Hilbert $B$-$\K_B(R)$-bimodule and $R \otimes_B \K_B(R,B) \cong \K_B(R)$. If $A=\K_B(R)$ this means that $\K_B(R,B)$ is a Hilbert $B$-$A$-bimodule and is left inverse to $R$.

Suppose that $R$ is not full. Then the inner product on $R$ takes its values in a proper ideal $J\subset B$. So, $R$ is really a Hilbert $J$-module and $RJ=R$. Suppose that $R$ is also right invertible. Then there exists 
a Hilbert $B$-$A$-bimodule $R_2$ such that $R_2\otimes_A R \cong B$, but
\[
B \cong R_2\otimes_A R = R_2\otimes_A RJ \cong BJ = J
\mbox,\]
which is a contradiction. Therefore right invertibility implies fullness.
\end{proof}

The conditions of left and right invertibility are obviously very symmetric, so each is half of Morita equivalence. This proposition shows that fullness is a little bit less than  half and $A=\K_B(R)$ is a little more than half. We can cut Morita equivalence into two conditions either directly down the middle or a little to the right.

There is a functor from the category of \cs-algebras and $*$-homomorphisms to \bimod. In fact, it is the restriction of a functor from the larger category of  \cs-algebras and strict completely positive maps. If $\rho:A\to B$ is a strict completely positive map, then $A\otimes_\rho B$ is the corresponding Hilbert $A$-$B$-bimodule.
\medskip

Since, by construction, the Toeplitz algebra is a subalgebra $\A_Q\subseteq\Li_B(\Hi_Q)$, $\Hi_Q$ is a Hilbert $\A_Q$-$B$-bimodule and defines a \bimod-morphism, $\Hi_Q:\A_Q\to B$.
The obvious question is whether this is a Morita equivalence, or if not, which of the above properties is true.

As it is, all that can be said \emph{a priori} is that if the fibers of $Q$ are finitely generated, then $\A_Q$ is a unital subalgebra of an algebra ($\K_B(\Hi_Q)$) which is Morita equivalent to an ideal in $B$ (the \cs-ideal generated by inner products in $\Hi_Q$). Another trivial, but not altogether useless, fact is that if $B$ is simple (has no proper \cs-ideals) then $\Hi_Q$ is full unless it is $0$.

None of the conditions for $\Hi_Q$ to be a Morita equivalence are always satisfied. Suppose that for some quantization bundle $\Hi_Q:\A_Q\to B$ is a Morita equivalence and consider the quantization bundle $Q\oplus Q$. The space of holomorphic sections is $\Hi_{Q\oplus Q} = \Hi_Q\oplus\Hi_Q$. The Toeplitz algebra is,
\[
\A_{Q\oplus Q} \cong \A_Q \subsetneq \K_B(\Hi_{Q\oplus Q}) \cong M_2(\A_Q)
\]
($2\times2$-matrices).
The Hilbert $B$-module $\Hi_{Q\oplus Q}$ \emph{is} full. Suppose that $R$ is a Hilbert $B$-$\A_Q$-bimodule, then
\[
R\otimes_{\A_Q}\Hi_{Q\oplus Q} \cong (R\otimes_{\A_Q}\Hi_Q)^{\oplus 2}
\]
which cannot be isomorphic to $B$; thus $\Hi_{Q\oplus Q}$ is not right invertible. Similarly,
\[
\Hi_{Q\oplus Q} \otimes_B R \cong (\Hi_Q\otimes_B R)^{\oplus 2}
\]
cannot be isomorphic to $\A_Q$; thus $\Hi_{Q\oplus Q}$ is not left invertible.

We can also consider $Q$ as a quantization bundle with coefficient algebra $B\oplus \co$. Since $\Hi_Q$ is a Hilbert $B$-module and $B\subset B\oplus \co$ is a proper ideal, $\Hi_Q$ is not a full Hilbert $B\oplus\co$-module. For the same reason, $\Hi_Q : \A_Q\to B\oplus\co$ cannot be right invertible. It is left invertible with the same inverse as $\Hi_Q:\A_Q\to B$, and 
\[
\A_Q = \K_B(\Hi_Q) = \K_{B\oplus\co}(\Hi_Q)
\mbox.\]

\subsection{Change of Coefficient Algebras}
There is an action of the category \bimod\ on the collection of possible quantization bundles over $\M$.
\begin{lem}\label{coefficients}
Let $Q$ be a quantization bundle with coefficient algebra $B$, $B'$ another unital \cs-algebra, and $R$ a Hilbert $B$-$B'$-bimodule. Then $Q':=Q\otimes_B R$ is a quantization bundle with coefficient algebra $B'$, and there is a surjective homomorphism $p :\A_Q\onto \A_{Q'}$ such that $T_{Q'}=p\circ T_Q$. This process is functorial in the sense that if $B''$ is another \cs-algebra and $R'$ is a Hilbert $B'$-$B''$-bimodule, then $p_{R\otimes_{B'}R'} = p_{R'}\circ p_R$.
\end{lem}
\begin{proof}
The action of $\otimes_B R$ is a functor from the category of Hilbertian $B$-modules (see p.~\pageref{Hilbertian.def}) to that of Hilbertian $B'$-modules. We can apply this in local trivializations of $Q$ to get a natural connection for $Q'$, so $Q'$ has the correct curvature and is a quantization bundle. In particular $\h_{Q'}=\h_Q\otimes_BR$. Locally, $\otimes_BR$, takes the $Q$ Sobolev spaces to the $Q'$ Sobolev spaces, and so (globally) it takes the domain of the $Q$ Dolbeault operator to the $Q'$ Dolbeault operator. In this sense, $D_{Q'} = D_Q\otimes_BR$. So, $\Hi_{Q'}=\Hi_Q\otimes_B R$ and $\Pi_{Q'} = \Pi_Q\otimes_BR$.

Let $p$ be the restriction to $\A_Q$ of the functorial map $\otimes_BR: \Li_B(\h_Q) \to \Li_{B'}(\h_{Q'})$. For $f\in\C(\M)$,
\[
p[T_Q(f)] = [T_Q(f)]\otimes_BR = (\Pi_Qf)\otimes_BR = (\Pi_Q\otimes_BR)f = \Pi_{Q'}f = T_{Q'}(f)
\]
so $p\circ T_Q = T_{Q'}$. 
Because the Toeplitz algebras are (by definition) generated by the images of the Toeplitz maps, $\A_{Q'}$ must be generated by $p(\A_Q)$. Therefore $\A_{Q'}=p(\A_Q)$, i.~e., $p:\A_Q\onto\A_{Q'}$ is surjective.

The last claim (functoriality) is true simply because if we compose the functors $\otimes_BR$ and $\otimes_{B'}R'$ then
\[
(\otimes_{B'}R')\circ(\otimes_BR) = \otimes_BR \otimes_{B'}R' = \otimes_B (R \otimes_{B'}R')
\mbox.\]
\end{proof}

Let $\As^{\max}\equiv \A_{\Lsm}$ be the Toeplitz algebra constructed with the maximal quantization bundle.
\begin{cor}
\label{universal2}
Under the general assumptions of Section \ref{Q.bundle.def} (i.~e., $\exists L^s$) for any quantization bundle $Q$ with curvature $s\omega$, the Toeplitz algebra is a quotient, $p:\Asm\onto\A_Q$, of the maximal Toeplitz algebra.
\end{cor}
\begin{proof}
By Thm.~\ref{universal1}, there exists a Hilbert $\Bsm$-$B$-bimodule $R$ such that 
\[
Q\cong \Lsm\otimes_\Bsm R
\mbox.\]
 The result is then immediate from Lem.~\ref{coefficients}.
\end{proof}

\begin{prop}
If $R$ is a Hilbert $B$-$B'$-bimodule which is left invertible in the above sense, then $R$ determines an isomorphism, $\A_Q\cong \A_{Q\otimes_B R}$.
\end{prop}
\begin{proof}
Left invertibility means that there exists a Hilbert $B'$-$B$-bimodule, $R'$ such that $R\otimes R' \cong B$. So 
\[
(Q\otimes_B R)\otimes_{B'} R' = Q
\mbox.\]
Lemma \ref{coefficients} shows that there exist natural surjective homomorphisms $p:\A_Q\onto\A_{Q\otimes_B R}$ and $p':\A_{Q\otimes_B R}\onto\A_Q$ such that the composition $p'\circ p$ is the identity. So, $p$ is injective and therefore an isomorphism.
\end{proof}

\subsection{Fullness}
Not much can be said generally about whether $\Hi_Q$ is a full Hilbert $B$-module.
\begin{lem}\label{Morita1}
$\Hi_Q$ is a full Hilbert $B$-module if and only if, for any nonzero homomorphism $\rho:B\to B'$, to another  \cs-algebra $B'\neq0$, the push-forward of $\Hi_Q$ is nonzero:
\[
\Hi_Q\otimes_\rho B' \neq 0
\mbox.\]
It is sufficient to consider just surjective homomorphisms. Equivalently, $\Hi_Q$ is full if and only if, for any nonzero representation, $\rho:B\to\Li(H)$,
\[
\Hi_Q\otimes_\rho H \neq 0
\mbox.\]
\end{lem}
\begin{proof}
Let $J\subseteq B$ be the \cs-ideal generated by inner products in $\Hi_Q$. $\Hi_Q$ is full if and only if $J=B$. If $\Hi_Q$ is not full then $B/J$ is a nonzero \cs-algebra, but the push forward of $\Hi_Q$ to a Hilbert $B/J$-module is trivial, $\Hi_Q\otimes_B(B/J)=0$. 

Conversely, suppose that $\Hi_Q$ is full. It is sufficient to check surjective homomorphisms, since any homomorphism is surjective onto its image. For any $\rho:B\onto B'$, the push-forward, $\Hi_Q\otimes_\rho B'$, is a full Hilbert $B'$-module and is therefore nonzero.

It is also sufficient to consider representations, since any \cs-algebra is isomorphic to a subalgebra of bounded operators on a Hilbert space. For any representation $\rho:B\to\Li(H)$,
\[
\Hi_Q\otimes_\rho\Li(H) = \Li(H,\Hi_Q\otimes_\rho H)
\]
so this is nonzero if and only if $\Hi_Q\otimes_\rho H$ is.
\end{proof}

The quantization bundles $\Ls$ and $\Lsm$ are the most important. Any representation of the coefficient algebras $\Bs$ or $\Bsm$ is a projective-unitary representation of the fundamental group $\G$. Given a projective-unitary representation $\rho:\G\to\Li(H)$ with twist $\ss$, we can construct a quantization bundle of Hilbert spaces, $\Lsm\otimes_\rho H$.
\begin{cor}
$\Hs$ is a full Hilbert $\Bs$-module if and only if, for any subrepresentation $\rho:\G\to\Li(H)$ of the $\ss$-twisted regular representation (see p.~\pageref{reg.rep}), the bundle $(\Ls\otimes_\rho H)\otimes L_0$ has a nonzero holomorphic section. Likewise, $\Hsm$ is a full Hilbert $\Bsm$-module if and only if, for any projective unitary representation $\rho:\G\to\Li(H)$ with twist $\ss$, the bundle $(\Lsm\otimes_\rho H)\otimes L_0$ has a nonzero holomorphic section.
\end{cor}
\begin{proof}
Any representation of $\Bsm$ is given by a projective-unitary representation of $\G$ with twist $\ss$. So, $\Hsm$ is full if and only if, for any such representation, $\Hsm\otimes_\rho H \neq 0$, but
\[
\Hsm\otimes_\rho H = \Hi_{\Lsm\otimes_\rho H} = \Gh[\M,(\Lsm\otimes_\rho H)\otimes L_0]
\mbox.\]

Any representation of $\Bs$ is given by a projective-unitary representation of $\G$ which factors through the $\ss$-twisted regular representation. It is thus sufficient to consider subrepresentations of this regular representation. In that case, $\Lsm\otimes_\rho H = \Ls\otimes_\rho H$.
\end{proof}

\subsection{Completeness}
\label{injectivity}
In this section I will prove that $\As=\K_\Bs(\Hs)$ and provide a way of computing any Toeplitz algebra $\A_Q$.

The idea of this proof is as follows. If $\A_Q=\K_B(\Hi_Q)$, then by Prop.~\ref{correspondence}, $\Hi_Q : \A_Q\to B$ is a left invertible \bimod\ morphism. This means that there exists a Hilbert $B$-$\A_Q$-bimodule, $R$, such that $\A_Q\cong \Hi_Q\otimes_B R$. Where can we get such a bimodule? If $Q'$ is a quantization bundle with coefficient algebra $B'$ and curvature $s\omega$, then Thm.~\ref{universal1} gives us a Hilbert $\Bsm$-$B'$-bimodule, $R$, such that
\[
\Hsm \otimes_\Bsm R = \Hi_{\Lsm\otimes_\Bsm R} = \Hi_{Q'}
\mbox.\]
If we can construct a quantization bundle, $Q'$, with coefficient algebra $\Asm$ and $\Hi_{Q'}=\Asm$ then we will at least prove that $\Hsm:\Asm\to\Bsm$ is left invertible.

If $\Asm=\K_\Bsm(\Hsm)$ then $\K_\Bsm(\Hsm,\Lsm)$ will be such a bundle. \emph{A priori}, $\K_\Bsm(\Hsm,\Lsm)$ is a quantization bundle with coefficient algebra $\K_\Bsm(\Hsm)$, and $\Hi_{\K_\Bsm(\Hsm,\Lsm)}=\K_\Bsm(\Hsm)$. We can try to show that there at least exists a subbundle $Q' \subseteq \K_\Bsm(\Hsm,\Lsm)$ of Hilbert $\Asm$-modules such that
\[
\Asm = \Hi_Q \subseteq \Hi_{\K_\Bsm(\Hsm,\Lsm)}=\K_\Bsm(\Hsm)
\mbox.\]

The key to this is that, for $s$ sufficiently large, we can reconstruct a quantization bundle $Q$ along with its connection entirely from the algebraic structures $\Hi_Q$ and $T_Q:\C(\M)\to \K_B(\Hi_Q)$. In this way, we can recognize when a subspace of $\Hi_Q$ is in fact the space of holomorphic sections of a subbundle.

The first fact I will need is that there exist enough holomorphic sections of $\Qb$ that $\Hi_Q$ actually generates the space of continuous sections, $\Gamma(\M,\Qb)$ as a $\C(\M)$-module. It will be sufficient to show that for every point, $x\in\M$, every vector in the fiber $\Qb_x$ is the value of some holomorphic section at $x$. This is a matter of ``ampleness''\/. If $\Qb$ were simply a line bundle (in the classical sense) then we would only need to show that there exists no point of $\M$ where all holomorphic sections of $\Qb$ vanish. This is a simpler version of the key lemma for the Kodaira embedding theorem, which says that for a line bundle of sufficiently positive curvature, there do not exists points $x\neq y\in\M$ such that any holomorphic section vanishing at $x$ also vanishes at $y$. I closely follow the proof of that lemma as given in \cite{g-h1}.

Again denote $\Qb:= Q\otimes L_0$ and $2n=\dim\M$.

\begin{lem}\label{ample}
Let $Q$ be a quantization bundle with finitely generated fibers and curvature $s\omega$. For $s$ sufficiently large, $\Gamma(\M,\Qb) = \C(\M)\Hi_Q \subset\h_Q$. 
\end{lem}
\begin{proof}
Assume that at least $s>\norm\Khat$, so that by Lem.~\ref{holomorphic},
$\Hi_Q\equiv\Gh(\M,\Qb)$ is the kernel of $D$. Note that holomorphic sections are continuous, so $\C(\M)\Hi_Q\subseteq\Gamma(\M,\Qb)$. 
I will first prove that for any $x\in\M$, the evaluation map $\Po_x:\Hi_Q\to \Qb_x$ is surjective. 

For any point $x\in\M$, let $\beta_x:\M_x\onto\M$ be the ``blow-up'' of $\M$  at $x$. The space $\M_x$ is a complex manifold which is a (not disjoint) union $\M_x=(\M\smallsetminus \{x\})\cup E_x$, $E_x\cong \co P^{n-1}$ and $\beta_x$ is the quotient map that identifies $E_x$ (the ``exceptional divisor'') to the point $x$. Associated to $E_x$ there is a holomorphic line bundle $[E_x]$ over $\M_x$ which has a canonical holomorphic section $\sigma_x\in\Gh(\M_x,[E_x])$ that vanishes precisely on $E_x$.

The restriction to $E_x$ of the pull-back $\beta_x^*\Qb$ is trivial, since it is the pull back of the restriction of $\Qb$ to the point $x$.  So, holomorphic sections of $\beta_x^*\Qb$ over $E_x$ must be constant.  Therefore the space of holomorphic sections of $\beta_x^*\Qb$ over $E_x$ is (naturally) isomorphic to the fiber of $\Qb$ over $x$. Stated in terms of sheaf cohomology, $H^0[\Or_{E_x}(\Qb] = \Qb_x$.

Because holomorphic sections are continuous, any holomorphic section of $\Qb$ is uniquely specified by its restriction to $\M\smallsetminus\{x\}$, but this extends uniquely to a holomorphic section of $\beta_x^*\Qb$ over $\M_x$. This gives an isomorphism of the spaces of holomorphic sections which is consistent with the identification of sections over $E_x$ with the fiber over $x$. In terms of sheaf cohomology, 
\[
H^0[\Or_{\M_x}(\beta_x^*\Qb)] = \Gh(\M_x,\beta_x^*\Qb) = \Gh(\M,\Qb) \equiv \Hi_Q
\mbox.\]

Let $[-E_x]=[E_x]^*$ be the dual bundle. Any local holomorphic section of $\beta_x^*\Qb\otimes [-E_x]$ can be multiplied with $\sigma_x\in\Gh(\M_x,[E_x])$ to give a local holomorphic section of $\beta_x^*\Qb$ which vanishes at $E_x$. This gives an isomorphism of the sheaf $\Or_{\M_x}(\beta_x^*\Qb\otimes [-E_x])$ of (local) holomorphic sections with the sheaf of holomorphic sections of $\beta_x^*\Qb$ which vanish at $E_x$. For this reason, there is a short exact sequence of sheaves,
\[
0 \to \Or_{\M_x}(\beta_x^*\Qb\otimes [-E_x]) \stackrel{\cdot\sigma_x}\longrightarrow \Or_{\M_x}(\beta_x^*\Qb) \longrightarrow \Or_{E_x}(\Qb) \to 0
\mbox.\]
This leads to a long exact sequence of sheaf cohomology groups. 
The degree $0$ cohomology groups are just the respective spaces of (global) holomorphic sections. Two of these $H^0$'s can be identified as described in the preceding paragraphs.  The beginning of the long exact sequence is thus,
\[
0 \to \Gh(\M_x,\beta_x^*\Qb\otimes [-E_x]) \longrightarrow
\Hi_Q \stackrel{\Po_x}\longrightarrow \Qb_x \longrightarrow H^1[\Or_{\M_x}(\beta_x^*\Qb\otimes [-E_x])] \to \dots
\mbox.\]
The map $\Po_x$ is the evaluation of holomorphic sections at $x$. To show that $\Po_x$ is surjective, it is sufficient to show that this $H^1$ vanishes.

The connection of $[E_x]$ can be constructed such that along $E_x$, the curvature of $[-E_x]$ is positive in directions tangent to $E_x$.  
The symplectic form $\omega$ is positive, so its pullback $\beta_x^*\omega$ is positive everywhere except $E_x$, where it is positive in all transversal directions. So, for some $m>0$, the sum of $m\beta_x^*\omega$ and the curvature of $[-E_x]$  will be positive.

The blow-up $\M_x$ embeds in the blow-up of $\M\times\M$ along the diagonally embedded $\M$. The bundle $[E_x]$ is just the restriction of a bundle on this larger blow-up. Since this larger blow-up is a compact manifold, we can actually choose $m$ sufficiently large for any $x\in\M$.

A \Kahler\ structure on $\M_x$ can be chosen such that 
\[
\wedge^{0,n}\M_x = \beta_x^*(\wedge^{0,n}\M)\otimes [-(n-1)E_x]
\mbox,\]
 where $[-kE_x]:=[-E_x]^{\otimes k}$. So,
\beq
\beta_x^*\Qb\otimes [-E_x] \otimes \wedge^{0,n}\M_x =
\beta_x^*(\Qb\otimes \wedge^{0,n}\M)\otimes [-nE_x]
\mbox.\label{positive.bundle}
\eeq
For $s$ large enough, $\Qb\otimes \wedge^{0,n}\M$ will have curvature more positive than $mn\omega$, therefore \eqref{positive.bundle} will have positive curvature. 

When the curvature of \eqref{positive.bundle} is positive, we can use it (the curvature) as a (different) \Kahler\ symplectic form for $\M_x$ and compute the cohomology of $\Or_{\M_x}(\beta_x^*\Qb\otimes [-E_x])$ as the kernel of a Dolbeault operator constructed with this \Kahler\ structure. 
Apply Lem.~\ref{holomorphic} with the bundle \eqref{positive.bundle}\footnote{Although $\Qb$ is a $B$-$\frac12$-density bundle, we don't need that structure here. This is the only point where some inner product is needed. We should ``divide'' $\Qb$ by the \Kahler\ volume form on $\M$, so as to give \eqref{positive.bundle} an inner product compatible with its connection.}  for $Q$, $\wedge^{n,0}\M_x$ for $L_0$, and $s=1$. By definition, $\hat K$ is constructed linearly from the curvature of $\wedge^{n,0}\M_x \otimes \wedge^{0,n}\M_x$ which is flat (and trivial), so $s=1>\norm{\hat K} = 0$ is sufficient to guarantee that the kernel of the Dolbeault operator is entirely of degree $0$. In other words, the higher degree cohomology groups of $\Or_{\M_x}(\beta_x^*\Qb\otimes [-E_x])$ vanish. In particular, $H^1$ vanishes and thus the evaluation of holomorphic sections of $\Qb$ at $x$ is surjective.

Given any point $x\in\M$, choose some (finite) set of generators $\{\psi_i\}\subset\Qb_x$ of the fiber (as a $B$-module). For each of these, choose a holomorphic section $\varphi_i$ such that $\varphi_i(x)=\psi_i$. There exists an open neighborhood of $x$ over which $\Qb$ is trivial. Within this there exists a neighborhood $U\ni x$ for which $\{\varphi\}$ generates the space of continuous sections $\Gz(U,\Qb)$ as a $\C_0(U)\otimes B$-module; therefore $\Hi_Q$ generates $\Gz(U,\Qb)$ as a $\C_0(U)$-module. 

Choosing such an open neighborhood for every $x\in\M$ we get an open cover --- and by compactness a finite cover --- of $\M$ by such neighborhoods. A partition of unity argument then shows that $\C(\M)\Hi_Q=\Gamma(\M,\Qb)$.
\end{proof}

\begin{cor}
$\C(\M)\otimes_{T_Q}\Hi_Q = L^2(\M,\Qb)$.
\end{cor}
\begin{proof}
The KSGNS tensor product $\C(\M)\otimes_{T_Q}\Hi_Q$ is the completion of the algebraic tensor product $\C(\M)\algtensor \Hi_Q$ with the inner product,
\[
\langle f\otimes \psi \vert g\otimes \varphi\rangle := \braket{\psi}{T_Q(\bar f g)}{\varphi} = \braket\psi{\bar f g}\varphi = \langle f\psi\vert g\varphi\rangle
\mbox.\]
This is thus the completion of the continuous sections of $\Qb$ and so it is $L^2(\M,\Qb)$.
\end{proof}

Classically, if we are given a holomorphic vector bundle with an inner product, then the connection can be reconstructed from the sheaf of local holomorphic sections. If the bundle has ``enough'' global holomorphic sections, then the connection can be reconstructed from the space of global holomorphic sections. We can do the same thing with $\Qb$. This makes it possible to recognize the space of holomorphic sections of a subbundle.

\begin{lem}\label{subbundle1}
In addition to the hypotheses of Lem.~\ref{ample}, suppose that $B'\subset B$ is a unital \cs-subalgebra. Let $\Hi'\subset \Hi_Q$ be a subspace which is a Hilbert $B'$-module (with the same inner product) and such that $\Hi'B=\Hi_Q$ and $\Hi_Q \cap (\C(\M)\Hi')= \Hi'$\/. Then there exists a subbundle $Q'\subset Q$ of Hilbert $B'$-modules preserved by the connection and such that $\Hi_{Q'}=\Hi'$.
\end{lem}
\begin{proof}
First note that the subspace $\C^\infty(\M) \Hi_Q\subset L^2(\M,\Qb)$, spanned by products of smooth functions and holomorphic sections, is dense and consists of smooth sections of $\Qb$. Within this subspace, we can construct inner products with derivatives. 

Let $f,g\in\C^\infty(\M)$, $\psi,\varphi\in\Hi_Q$, and $\xi\in\Gamma^\infty(\M,\TM)$ a smooth vector field. Decompose $\xi$ as a sum, $\xi=\xi_{\partial}+\xi_{\bar\partial}$, of sections of the holomorphic and antiholomorphic tangent bundles, respectively. Because $\Hi_Q$ consists of holomorphic sections,
\[
\nabla_{\xi_{\bar\partial}}\lvert\varphi\rangle =0
\]
so,
\[
\nabla_{\xi_{\bar\partial}} \left(g \lvert\varphi\rangle\right)
= \xi_{\bar\partial}(g)\lvert\varphi\rangle
\mbox.\]
The connection on $\Qb$ is compatible with the inner product relative to the \Kahler\ volume form. The adjoint of $\nabla_{\xi_\partial}$ involves the divergence $\nabla\cdot\xi_\partial$. Since $\psi$ is holomorphic, this shows that,
\[
\langle\psi\rvert \bar f \nabla_{\xi_\partial} = \langle\psi\rvert \left[\bar f (\nabla\cdot\xi_\partial) - \xi_\partial(\bar f)\right]
\mbox.\]
Using these identities, we have,
\begin{align}
\braket{f \psi}{\nabla_\xi}{g\varphi} 
&= \braket{\psi}{\bar f \nabla_\xi g}{\varphi} 
\nonumber\\
&= \braket{\psi}{T_Q\left[\bar f \,\xi_{\bar\partial}(g) + \bar f(\nabla\cdot\xi_\partial)g - \xi_{\partial}(\bar f) g \right]}{\varphi}
\mbox.\label{recover.connection}
\end{align}

The Hilbert $B'$-module $\C(\M)\otimes_{T_Q}\Hi'$ is the completion of $\C(\M)\Hi'$. By the Kasparov stabilization theorem, this embeds as a direct summand of the canonical Hilbert $B'$-module, ${B'}^{\oplus\infty}$\/. Let $\{e_i\}$ be the canonical basis of ${B'}^{\oplus\infty}$\/. By assumption, $\Hi_Q = \Hi' B = \Hi'\otimes_{B'}B$, so 
\[
(\C(\M)\otimes_{T_Q}\Hi')\otimes_{B'}B = \C(\M)\otimes_{T_Q}(\Hi'\otimes_{B'}B) = L^2(\M,\Qb)
\mbox,\]
 and we can use the same canonical basis for $B^{\oplus\infty} \supset L^2(\M,\Qb)$.

Suppose that $\varphi\in\Hi'$. Let $\varphi' := \nabla_\xi f\varphi$. Since $f\varphi$ is a smooth section of $\Qb$, $\varphi'$ must be a smooth section of $\Qb$. Because $L^2(\M,\Qb)\subset B^{\oplus\infty}$ is a direct summand, we can project any basis vector $e_i$ to the subspace $\C(\M)\otimes_{T_Q}\Hi'$ and then express this as a limit of vectors in $\C^\infty(\M)\Hi'$. In this way, eq.~\eqref{recover.connection} shows that $\langle e_i\vert \varphi'\rangle \in B'$.
The inner product $\langle e_i\vert\varphi'\rangle$ is the $i$'th component of $\varphi'$\/, this shows that $\varphi'\in {B'}^{\oplus\infty}$. We thus have,
\[
\varphi' \equiv \nabla_\xi f\varphi\in \Gamma^\infty(\M,\Qb) \cap (\C(\M)\otimes_{T_Q}\Hi')
\mbox.\]

This shows that any derivative $\nabla_\xi$ preserves the space $\Gamma^\infty(\M,\Qb) \cap (\C(\M)\otimes_{T_Q}\Hi')$. We can use the connection to construct a local trivialization. This shows that this is the space of smooth sections of a subbundle $\Qb'\subset\Qb$ preserved by the connection. Tensoring with the dual of $L_0$ gives the desired subbundle $Q' := (L_0)^*\otimes \Qb' \subset Q$.

By definition, $\Hi_{Q'}$ is the space of holomorphic sections,
\[
\Hi_{Q'} = \Gh(\M,\Qb') 
= \Hi_Q \cap (\C(\M)\Hi')
\mbox.\]
By assumption, this is $\Hi'$.
\end{proof}

Now consider the Hilbert $\C(\M)$-$\A_Q$-bimodule, $\C(\M)\otimes_{T_Q}\A_Q$. This is a subspace of the $\C(\M)$-$\K_B(\Hi_Q)$-bimodule,
\begin{align*}
\C(\M)\otimes_{T_Q} \K_B(\Hi_Q) 
&= \K_B(\Hi_Q,\C(\M)\otimes_{T_Q}\Hi_Q) \\
&= \K_B(\Hi_Q,L^2[\M,\Qb]) \\
&= L^2[\M,\K_B(\Hi_Q,Q)\otimes L_0]
\mbox,\end{align*}
if $s$ is large enough for Lem.~\ref{subbundle1}.

\begin{lem}\label{subbundle}
Under the hypotheses of Lem.~\ref{ample}, there exists another quantization bundle, $Q'$, with coefficient algebra $\A_Q$ and curvature $s\omega$, such that $\Hi_{Q'} \cong \A_Q$ (as a Hilbert $\A_Q$-bimodule).
\end{lem}
\begin{proof}
The bundle $\K_B(\Hi_Q,Q) = Q \otimes_B \K_B(\Hi_Q,B)$ is a quantization bundle of Hil\-bert $\K_B(\Hi_Q)$-modules with
$\Hi_{\K_B(\Hi_Q,Q)} = \Hi_Q \otimes_B \K_B(\Hi_Q,B) = \K_B(\Hi_Q)$.

We can view $\A_Q$ as a Hilbert $\A_Q$-module and it is a subspace of the Hilbert $\K_B(\Hi_Q)$-module 
\[
\Hi_{\K_B(\Hi_Q,Q)} = \K_B(\Hi_Q)
\mbox.\] 
Because $1\in \A_Q$, 
\[
\A_Q\, \K_B(\Hi_Q) = \K_B(\Hi_Q)
\mbox.\] 
The inner product on the subspace $\C(\M)\A_Q \subset L^2[\M,\K_B(\Hi_Q,Q)\otimes L_0]$ takes its values in $\A_Q$. If $a \in \Hi_{\K_B(\Hi_Q,Q)}$, then $\langle a \vert a \rangle = a^*a \in \A_Q$ if and only if $a\in \A_Q$. Therefore 
\[
\Hi_{\K_B(\Hi_Q,Q)} \cap (\C(\M)\A_Q) = \A_Q
\mbox.\]
So, $\A_Q \subseteq \Hi_{\K_B(\Hi_Q,Q)}$ satisfies the hypothesis of Lem.~\ref{subbundle1}, and thus there exists a subbundle $Q' \subseteq \K_B(\Hi_Q,Q)$ of Hilbert $\A_Q$-modules with $\Hi_{Q'} = \A_Q$.
\end{proof}

\begin{thm}\label{Morita2}
For $s$ sufficiently large: 
\begin{enumerate}
\item
$\As = \K_{\Bs}(\Hs)$.
\item
$\Asm = \K_{\Bsm}(\Hsm)$.
\item
For any quantization bundle $Q$ with curvature $s\omega$, $\A_Q$ is the image of the functorial map,
\[
\otimes_\Bsm R : \K_\Bsm(\Hsm) \to \Li_B(\Hi_Q)
\mbox,\]
where $R$ is the Hilbert $\Bsm$-$B$-bimodule in Thm.~\ref{universal1}.
\end{enumerate}
\end{thm}
\begin{proof}
The second claim is the key.

Apply Lem.~\ref{subbundle} with $Q=\Lsm$. Since $Q'$ satisfies the hypothesis of Thm.~\ref{universal1}, there exists a Hilbert $\Bsm$-$\Asm$-bimodule $R$ such that $Q'=\Lsm\otimes_{\Bsm}R$, and thus, $\Asm \cong \Hi_{Q'} =\Hsm\otimes_{\Bsm}R$.

By construction, $Q' \subseteq \Lsm \otimes_\Bsm \K_\Bsm(\Hsm,\Bsm)$, so $R \subseteq \K_\Bsm(\Hsm,\Bsm)$ is a $\Bsm$-submodule. This latter space can be identified with the complex conjugate of $\Hsm$. It is a Hilbert $\Bsm$-$\K_\Bsm(\Hsm)$-bimodule and it is convenient to represent the vectors as ``bra''s. Since $\Hsm\otimes_\Bsm\overline\Hsm \cong \Hsm\overline\Hsm = \K_\Bsm(\Hsm)$, we can consider the tensor product $\Hsm\otimes_\Bsm R$ as a subspace of $\K_\Bsm(\Hsm)$. We know that,
\[
1\in \Asm \cong \Hsm \otimes_\Bsm R \cong \Hsm R \subseteq \K_\Bsm(\Hsm)
\mbox.\]
So, we can write the identity operator on $\Hsm$ as a finite sum,
\[
1 = \sum_i \lvert\psi_i\rangle\langle\varphi_i\rvert
\mbox,\]
with each $\lvert\psi_i\rangle \in \Hsm$ and $\langle\varphi_i\rvert \in R$. However, any $\langle\chi\rvert \in \overline\Hsm$ can be rewritten as,
\[
\langle\chi\rvert = \langle\chi\rvert 1 = \sum_i \langle\chi\vert\psi_i\rangle\langle\varphi_i\rvert 
\mbox.\]
Since each $\langle\chi\vert\psi_i\rangle \in \Bsm$, this shows that $\langle\chi\rvert \in R$, and therefore $\overline\Hsm = R$. Thus $\Asm=\K_\Bsm(\Hsm)$, which was the second claim.

With that done, now let $R$ be the $\Bsm$-$B$-bimodule constructed in Thm.~\ref{universal1} for some arbitrary quantization bundle, $Q$, with curvature $s\omega$. By Lem.~\ref{coefficients}, there is a surjective homomorphism,
\[
p : \Asm = \K_\Bsm(\Hsm) \onto \A_Q
\mbox,\]
which is a restriction of the functorial map $\otimes_\Bsm R$. The third claim follows.

The quantization bundle $\Ls$ is a quotient of $\Lsm$ and can be constructed as $\Ls = \Lsm\otimes_\Bsm \Bs$. Thus $\Hs = \Hsm\otimes_\Bsm\Bs$ is a quotient of $\Hsm$. Any compact operator on $\Hs$ is in the image of the functorial map,
\[
\otimes_\Bsm\Bs : \K_\Bsm(\Hsm) \onto \K_\Bs(\Hs)
\mbox.\]
So, $\As = \K_\Bs(\Hs)$, which was the first claim.
\end{proof}

This shows in particular that if $\Hs$ is full, then it is a Morita equivalence between $\As$ and $\Bs$. 

In the standard Toeplitz quantization of a noncompact manifold, the Toeplitz algebra is simply the algebra of compact operators on the Hilbert space. If we generalize the present construction to noncompact manifolds, then the Hilbert $\Bs$-module, $\Hs$, will not be finitely generated. I expect that Thm.~\ref{Morita2} will continue to hold as stated. The Toeplitz algebra should be exactly $\K_\Bs(\Hs)$.

\section{Quantization}
\label{quantization}
Let us now turn to properties which are specific to the most important choice of quantization bundle; namely $Q=\Ls$. We now use the full set of assumptions about $\M$. It is compact and \Kahler\ and the lift of the symplectic form to $\Mt$ is exact.

Recall again that the coefficient algebra in this case is the reduced, twisted, group \cs-algebra, $\Bs=\csr(\G,\ss)$. First, I shall demonstrate the significance of $\Ls$.

\subsection{The Covering Construction}
We can define the standard Toeplitz construction for the universal covering space as follows. 
\begin{definition}
Let $L_s:=L^s\otimes\pi^*L_0$. Let $\Hs^\Mt\subset L^2(\Mt,L_s)$ be the Hilbert space of square integrable holomorphic sections of $L_s$ over $\Mt$. Let $\Pi: L^2(\Mt,L_s) \onto \Hs^\Mt$ be the orthogonal projection. The Toeplitz operator of a bounded, continuous function $f\in\Cb(\Mt)$ is,
\[
T_s^\Mt(f) := \Pi f : \Hs^\Mt\to \Hs^\Mt
\mbox.\]
\end{definition}
\begin{lem}
$T_s^\Mt:\Cb(\Mt)\to \Li(\Hs^\Mt)$.
\end{lem}
\begin{proof}
Being an orthogonal projection, $\Pi$ is bounded. As a multiplication operator, $f$ is bounded.
\end{proof}

Regard $\C(\M)\subset\Cb(\Mt)$ as the subalgebra of $\G$-invariant functions. 
The \cs-algebra generated by the image $T_s^\Mt[\C(\M)]$  is the natural generalization of Klimek and Lesniewski's quantum Riemann surfaces \cite{k-l2,k-l4}.
Recall that $T_s$ and $\As$ are the Toeplitz map and algebra constructed using the quantization bundle $\Ls$.
\begin{thm}\label{standard}
If $s>\norm\Khat$ then the \cs-algebra generated by the image $T_s^\Mt[\C(\M)]$ is isomorphic to $\As$ and can be identified with it in such a way that $T_s=T_s^\Mt|_{\C(\M)}$.
\end{thm}
\begin{proof}
We can complete the bundle $\Ls$ of Hilbert $\Bs$-modules to a bundle of Hilbert spaces, $Q:=\Ls\otimes_{\tau_s}\co$, and use this as a quantization bundle with coefficient algebra $\co$. By Lem.~\ref{Lcompletion}, the Hilbert space, $\h_Q=\hs\otimes_{\tau_s}\co$ is a completion of $\Gc(\Mt,L_s\otimes\wedge^{0,*}\Mt)$; the inner product, $\tau_s(\langle\,\cdot\,|\,\cdot\,\rangle)$, is the canonical one, so $\h_Q = L^2(\Mt,L_s\otimes\wedge^{0,*}\Mt)$.  $\Hi_Q=\Hs\otimes_{\tau_s}\co$ is the Hilbert space of holomorphic sections of $Q\otimes L_0$. However, since the connection on $Q$ coincides with that on $L^s$, this is simply the Hilbert space of holomorphic sections of $L_s$. So, $T_Q = T_s^\Mt|_{\C(\M)}$.

The KSGNS tensor product can be factorized as, $\Hs\otimes_{\tau_s}\co = \Hs\otimes_{\Bs}\Bs\otimes_{\tau_s}\co$.  The Hilbert space $\Bs\otimes_{\tau_s}\co\cong l^2(\G)$ is the GNS Hilbert space of $\Bs$ constructed with the state $\tau_s$. It can be thought of as a Hilbert $\Bs$-$\co$-bimodule, so Lem.~\ref{coefficients} applies and we have a natural surjective homomorphism, $p:\As\onto\A_Q$. 

Recall that $p$ is the restriction of the functorial map, 
\[
\otimes_\rho\co : \Li_{\Bs}(\Hs)\to\Li(\Hs\otimes_{\tau_s}\co)
\mbox.\]
 The Hilbert space $\Hs\otimes_{\tau_s}\co$ is simply a closure, $\Hs\otimes_{\tau_s}\co\supset\Hs$, so a nonzero bounded-adjoint\-able operator on $\Hs$ must give a nonzero operator on $\Hs\otimes_{\tau_s}\co$, thus $p$ must be injective. So, $p$ is an isomorphism. 
\end{proof}

\subsection{Continuous Field}
I have already shown in Section \ref{asymptotic} that the asymptotic multiplicativity properties required of a strict quantization are satisfied by the Toeplitz maps $T_Q$ for any quantization bundles as $s\to\infty$. What remains is to show that in the special case of the quantization bundles $\Ls$, the Toeplitz construction actually gives a continuous field of \cs-algebras. The key properties that make this possible are that the coefficient algebras $\Bs$ form a continuous field, and that the quantization bundles, $\Ls$, have a compatible structure. Even without the Toeplitz maps, there is a continuous field structure to the Hilbert \cs-modules, $\Hs$, and the algebras of compact operators on these.

By Lemmas \ref{Lcompletion} and \ref{compact}, the fibers of $\Ls$ are finitely generated $\Bs$-modules and  $\Li_\Bs(\Hs)=\K_\Bs(\Hs)$. 
Choose some $s_0>\norm\Khat$ and denote the semi-infinite interval $I:=[s_0,\infty)$. This is simply because it is easier to work with the closed interval $I$, than the open interval $(\norm\Khat,\infty)$. Since we need to deal with all values of $s$, I will now denote the projection onto $\Hs$ as $\Pi_s$.

Although the bundle $\Lb$ in Lem.~\ref{FieldBundle} is not a quantization bundle (it has no connection) it is nevertheless a bundle of Hilbert $\Gz(\R,\B)$-modules and  we can still construct the Hilbert $\Gz(\R,\B)$-module, $\h_\Lb$. By Lem.~\ref{ModuleField}, this is the space of $\C_0$ sections of a continuous field of Hilbert \cs-modules. By Lem.~\ref{FieldBundle}, the fiber at $s\in\R$ is just $\hs$. Denote this field as $\h$. Note that the space of sections over some $J\subset\R$ is, $\Gz(J,\h) = \h_\Lb\otimes_{\Gz(I,\B)} \Gz(J,\B)$.

\begin{lem}\label{continuous.projection}
The projections, $\Pi_s$, give a continuous section 
\[
\Pi \in \Gb[I,\K_\B(\h)] \subset \Li_{\Gz(I,\B)}(\Gz[I,\h])
\mbox,\]
 which is a projection.
\end{lem}
\begin{proof}
We can define $\Pi$ as the section given by the kernel projections, $\Pi_s$. We need to prove that $\Pi$ is a \emph{continuous} section.

Continuity is a local issue. It is sufficient to prove continuity for every compact interval $J\subset I$. For such an interval, the restriction of $\Pi$ is the kernel projection of the Dolbeault operator associated with $\Lb^{J}$. By Lem.~\ref{parametrix}, this is compact. So, by Lem.~\ref{ModuleField}, it is a continuous section of $\K_\B(\h)$ over $J$. Since these projections all have norm $1$, $\Pi$ is a bounded section of $\K_\B(\h)$ over $I$.

 Lemma \ref{ModuleField} shows that $\Gb[I,\K_\B(\h)] \subset \Li_{\Gz(I,\B)}(\Gz[I,\h])$. The section, $\Pi$, is thus bounded-adjoint\-able; it is a projection because it is a section of projections.
\end{proof}

\begin{cor}
The Hilbert \cs-modules $\Hs$ form a continuous field of Banach spaces, $\Hi$, over $I$.
\end{cor}
\begin{proof}
Because $\Pi$ is a projection, its image is a Hilbert $\Gz(I,\B)$-submodule of $\Gz(I,\h)$. Lem.~\ref{ModuleField} shows that this must be the space of $\C_0$-sections of a continuous field of Banach spaces. The fiber over $s\in I$ is the image of $\Pi_s$, which is just $\Hs$.
\end{proof}

\begin{lem}
\label{adjointable.continuous}
The spaces $\K_\Bs(\Hs)$ form a continuous field of \cs-algebras $\K_{\B}(\Hi)$ over $I$.  Over any closed interval $J \subseteq I$, the $\C_0$ sections are
\beq
\label{C0.sections}
\Gz[J,\K_{\B}(\Hi)] = \K_{\Gz(J,\B)}(\Gz[J,\Hi])
\mbox,\eeq
and the continuous, bounded sections are
\beq
\label{Cb.sections}
\Gb[J,\K_{\B}(\Hi)] = \Li_{\Gz(J,\B)}(\Gz[J,\Hi])
\mbox.\eeq
\end{lem}
\begin{proof}
The existence of $\K_\B(\Hi)$ comes directly from Lem.~\ref{ModuleField}, as does eq.~\eqref{C0.sections} for $J=I$. For $J\subset I$, eq.~\eqref{C0.sections} comes by applying Lem.~\ref{ModuleField} to the Hilbert $\Gamma(J,\B)$-module, $\Gamma(J,\Hi)$.

Consider some $a\in \Li_{\Gz(I,\B)}(\Gz[I,\Hi])$. For any compact interval $J\subset I$, the space of sections of $\Hi$ over $J$ is the push forward, $\Gamma(J,\Hi) = \Gz(I,\Hi) \otimes_{\Gz(I,\B)}\Gamma(J,\B)$. So, there is a functorial map, $\Li_{\Gz(I,\B)}(\Gz[I,\Hi]) \onto \Li_{\Gamma(J,\B)}(\Gamma[J,\Hi])$. However, because $J$ is compact, $\Gamma(J,\Hi)$ is finitely generated and,
\[
\Li_{\Gamma(J,\B)}(\Gamma[J,\Hi]) 
= \K_{\Gamma(J,\B)}(\Gamma[J,\Hi])
= \Gamma[J,\K_\B(\Hi)]
\mbox.\]
Thus $a$ gives a continuous section of $\K_\B(\Hi)$ over any compact interval. These fit together, because the push-forwards are functorial. So, $a$ defines a continuous (and certainly bounded) section of $\K_\B(\Hi)$ over $I$. We already have $\Gb[J,\K_{\B}(\Hi)] \subseteq \Li_{\Gz(J,\B)}(\Gz[J,\Hi])$ from Lem.~\ref{ModuleField}, so this gives eq.~\eqref{Cb.sections}.
\end{proof}

The natural question now is whether the Toeplitz maps $T_s$ are consistent with this continuous field structure. They are:
\begin{lem}\label{big.Toeplitz}
For any function $f\in\C(\M)$,
\[
s\mapsto T_s(f)
\]
defines a continuous, bounded section of $\K_\B(\Hi)$.
\end{lem}
\begin{proof}
Since $\Im \Pi = \Gz(I,\Hi)$, we can regard $\Pi$ as a bounded-adjoint\-able map from $\Gz(I,\h)$ to $\Gz(I,\Hi)$.
The ``big'' Toeplitz operator of $f\in \C(\M)$ is,
\[
T(f) := \Pi f : \Gz(I,\Hi) \to \Gz(I,\Hi)
\mbox.\]

The Hilbert $\Gz(I,\B)$-module $\Gz(I,\h)$ is the space of $L^2$ sections of a $\Gz(I,\B)$-$\frac12$-density bundle. So, by Lem.~\ref{module.bundle2}, $f$ acts as a bounded-adjoint\-able operator. Since $\Pi$ is also bounded-adjoint\-able, $T(f)$ is bounded-adjoint\-able. Using Lem.~\ref{adjointable.continuous},
\[
T(f) \in \Li_{\Gz(I,\B)}(\Gz[I,\Hi]) = \Gb[I,\K_{\B}(\Hi)]
\mbox.\]

Now consider the diagram,
\[
\begin{CD}
\Gz(I,\Hi) @>{f}>> \Gz(I,\h) @>{\Pi}>> \Gz(I,\Hi) \\
@VVV @VVV @VVV \\
\Hs @>{f}>> \hs @>{\Pi_s}>> \Hs
\end{CD}
\]
The vertical maps are the evaluations at $s\in I$ for the continuous fields $\Hi$ and $\h$. This diagram is commutative because (by Lem.~\ref{continuous.projection}) $\Pi_s$ is the evaluation of $\Pi$ at $s$. The composition of the top row is $T(f)$. The composition of the bottom row is $T_s(f)$. This shows that the evaluation of $T(f)$ at $s$ is $T_s(f)$. 
\end{proof}

In order to prove that the Toeplitz algebras form a continuous field over $[\norm\Khat,\infty]$ using Lem.~\ref{field}, we will need a family of faithful states on the algebras. In the case of $\B$, we used the canonical tracial states $\tau_s$ for this purpose. The algebras $\K_\Bs(\Hs)$ also have canonical traces.
Define $\tr_s$ on ``rank one'' operators in $\K_\Bs(\Hs)$ by,
\beq\label{trace2}
\tr_s\left(\lvert\psi\rangle\langle\varphi\rvert\right) 
:= \tau_s\left(\langle\varphi\vert\psi\rangle\right)
\mbox.\eeq
\begin{lem}
This defines a positive, faithful trace, $\tr_s : \K_\Bs(\Hs) \to\co$.
\end{lem}
\begin{proof}
By Lem.~\ref{compact}, $\Hs$ is finitely generated. So,  $\Hs$ is a complemented submodule of $\B_s^m$ for some finite $m$. The algebra $\K_\Bs(\Bs^m)=M_m(\Bs)$ is the set of $m\times m$ matrices over $\Bs$ and so consists entirely of finite-rank operators. Therefore $\K_\Bs(\Hs)$ consists entirely of finite-rank operators, and $\tr_s$ is defined (and finite) on all of $\K_\Bs(\Hs)$.

If we regard this as a subalgebra, $\K_\Bs(\Hs)\subset M_m(\Bs)$, then $\tr_s$ is the restriction of the canonical trace constructed from $\tau_s$ on the latter algebra. That canonical trace is positive and faithful, therefore $\tr_s$ is.
\end{proof}

Let $\Po_s:\Gb[I,\K_\B(\Hi)]\onto \K_\Bs(\Hs)$ be the evaluation homomorphism at $s\in I$.

\begin{lem}\label{big.trace}
For any $a\in \Gb[I,\K_\B(\Hi)]$, $\tr_s[\Po_s(a)]$ is a continuous function of $s$.
\end{lem}
\begin{proof}
It is sufficient to check continuity on any compact interval $J\subset I$.

Let $\tr$ be the $\C(J)$-linear trace defined by,
\[
\tr(a) : s \mapsto \tr_s[\Po_s(a)]
\mbox.\]
Since $\tr_s$ is uniquely defined by eq.~\eqref{trace2}, $\tr$ is uniquely defined as the $\C(J)$-linear map satisfying the direct analogue and consequence of eq.~\eqref{trace2}. That is, for any $\psi,\varphi\in\Gz(J,\Hi)$,
\beq
\tr\left(\lvert\psi\rangle\langle\varphi\rvert\right) 
= \tau\left(\langle\varphi\vert\psi\rangle\right)
\mbox.\label{trace3}
\eeq

$\Gamma(J,\Hi)$ is a finitely generated $\Gamma(J,\B)$-module. Therefore,  $\Gamma[J,\K_\B(\Hi)]$ consists entirely of finite-rank operators. So, $\tr$ is defined on $\Gamma[J,\K_\B(\Hi)]$ and is $\C(J)$-valued.

Thus $\tr_s[\Po_s(a)]$ is a continuous (but not usually bounded) function of $s\in I$.
\end{proof}

Now consider the normalized traces defined by $\ntr_s(a) := \tr_s(a)/\tr_s(1)$ for any $a\in \K_{\Bs}(\Hs)$. Because $\tr_s$ is positive and faithful, $\ntr_s$ is a faithful state. 
\begin{lem}
\label{trace.limit}
For any $f\in\C(\M)$,
\[
\lim_{s\to\infty} \ntr_s\left[T_s(f)\right] = \frac1{\vol\M}\intM f \frac{\omega^n}{n!}
\mbox,\]
where $\frac{\omega^n}{n!}$ is the \Kahler\ volume form.
\end{lem}
\begin{proof}
Because $\ntr_s$ is a trace, for any $f,g\in\C^\infty(\M)$,
\[
\ntr_s [T_s(f),T_s(g)]_- = 0
\mbox.\]
Because $\ntr_s$ is a state, it is norm-contracting and so Lem.~\ref{1st.order} shows that,
\[
\ntr_s\left[T_s (\{f,g\})\right] = \Or^{-1}(s)
\mbox.\]

The Poisson bracket $\{f,g\}$ is the divergence of $g\,\xi_f$, where $\xi_f$ is the Hamiltonian vector field for $f$ (such that $\xi_f\inner\omega = df$).
Any vector field can be written as a finite sum of vector fields of this form. Simply take a finite cover of $\M$ by Darboux coordinate patches. The $f$'s will be coordinate functions and the $g$'s will be components. This shows that for any smooth vector field $\xi\in\Gi(\M,\TM)$ the divergence $\nabla\cdot\xi$ is a finite sum of Poisson brackets and so
\[
\ntr_s\left[T_s(\nabla\cdot\xi)\right] = \Or^{-1}(s)
\mbox.\]

Because $\M$ is compact and connected, any smooth function $f\in\C^\infty(\M)$ can be written as the sum of its mean value and the divergence of some vector field,
\[
f = \nabla\cdot\xi + \frac1{\vol\M}\intM f\frac{\omega^n}{n!}
\mbox.\]
Applying $\ntr_s$ to this equation establishes the lemma for smooth functions, and because $\C^\infty(\M)\subset\C(\M)$ is dense and $\ntr_s$ is norm-contracting, this proves the lemma.
\end{proof}

\begin{thm}
\label{A.field}
There is a unique continuous field of \cs-algebras $\A$ over $(\norm\Khat,\infty]$ such that:
\begin{enumerate}
\item 
The fiber of $\A$ over $s$ is $\As$.
\item 
The fiber over $\infty$ is $\C(\M)$.
\item
For any $f\in\C(\M)$,
\[
s \mapsto 
\begin{cases}
T_s(f) & : \norm\Khat<s<\infty\\
f & : s=\infty
\end{cases}
\]
defines a continuous section of $\A$ over $(\norm\Khat,\infty]$.
\end{enumerate}
\end{thm}
\begin{proof}
Again, $s_0$ is any number $s_0>\norm\Khat$ and $I$ is the interval $I=[s_0,\infty)$.

For any $f\in\C(\M)$, let $T(f)$ be the section of $\K_\B(\Hi)$ defined in Lem.~\ref{big.Toeplitz} by $s\mapsto T_s(f)$. This defines a map $T:\C(\M)\to\Gb[I,\K_\B(\Hi)]$. Let 
\[
\AA_0 := \{a\in \Gz(I,\K_\B[\Hi]) \mid \forall s\in I,\, a(s) \in \As\}
\mbox.\] 
This is not empty; at least $0\neq \C_0(I)\Im T \subseteq \AA_0$.

Any Cauchy sequence, $\{a_i\}_{i=1}^\infty \subset \AA_0 \subseteq \Gz(I,\K_\B[\Hi])$ must have a limit $a\in \Gz(I,\K_\B[\Hi])$ because $\Gz(I,\K_\B[\Hi])$ is closed. For any $s\in I$, the sequence of evaluations 
\[
\{a_i(s)\}_{i=1}^\infty \subset \As
\]
 must converge to $a(s)$, but because $\As$ is closed, $a(s)\in \As$. Thus $a\in\AA_0$. Therefore $\AA_0$ is closed.

Let $\AA := \Im T + \AA_0 \subset \Gb(I,\K_\B[\Hi])$. Let $\hat I := [s_0,\infty]$ be the one-point compactification of $I$. $\AA$ will be the algebra of continuous sections of $\A$ over $\hat I$; I first verify that it is a \cs-algebra.

For $f\in\C^1(\M)$ and $g\in\C(\M)$, Lem.~\ref{0th.order} shows that
\[
\lim_{s\to\infty}\Norm{T_s(f)T_s(g)-T_s(fg)} = 0
\mbox.\]
So $T(f)T(g)-T(fg)$ is a continuous section of $\K_\B(\Hi)$ which goes to $0$ at $\infty\in\hat I$ and by construction takes values in the algebras $\As$. That is,
\beq
T(f)T(g)-T(fg)\in\AA_0
\label{A0}\mbox.\eeq
For fixed $g$, the map $f \mapsto T(f)T(g)-T(fg)$ does not increase norm by more than a factor of $2\norm g$. Therefore it is continuous and \eqref{A0} holds for all $f,g\in\C(\M)$. 

This shows that $\AA$ is algebraically closed. The map $T$ is contractive, so $\AA$ is topologically closed. $T$ is also ${}^*$-invariant, so $\AA$ is a \cs-algebra. The algebra $\C(\hat I)$ is contained in $\AA$ by construction; it is central in $\AA$ since it is central in $\Gb[I,\K_\B(\Hi)]$. So, by Lem.~\ref{field}, $\AA$ must be the algebra of continuous sections of some upper semicontinuous field of \cs-algebras over $\hat I$. Denote this as $\A$. 

Since $\AA \subset \Gb(I,\K_\B[\Hi])$, the fiber of $\A$ over $s\in I$ is a \cs-subalgebra of $\K_\Bs(\Hs)$. It must contain the image of $T_s$ and it must be contained in $\As$. Therefore the fiber of $\A$ over $s$ is $\As$.

Let $\Po_\infty:\AA\onto\A_\infty$ be the evaluation homomorphism at $\infty\in\hat I$. By definition, $\ker \Po_\infty = \C_0(I)\AA$.
However,
\[
\AA_0 = \C_0(I) \AA_0 \subseteq \C_0(I) \AA  
\mbox.\]
and conversely, any section in $\C_0(I) \AA$ vanishes at $\infty$ and takes its values in $\As$. Thus $\ker \Po_\infty = \AA_0$, and
\eqref{A0} shows that $\Po_\infty\circ T:\C(\M)\to\A_\infty$ is a homomorphism. The definition of $\AA$ shows that $\Po_\infty\circ T$ is surjective, so $\A_\infty$ is a quotient of $\C(\M)$.
 
Lemma \ref{big.trace} says that, for $a\in\Gb(I,\K_\B[\Hi])$, $\tr_s[\Po_s(a)]$ is a continuous function of $s$. In particular $\tr_s 1$ is continuous and by fidelity, nonvanishing. Therefore the normalized trace $\ntr_s[\Po_s(a)]$ is also a continuous function of finite $s$. We can thus put the normalized traces together into a $\C(\hat I)$-linear map, $\ntr:\AA\to\C_{\mathrm b}(I)$.  
Lemma \ref{trace.limit} implies that $\ntr\circ T : \C(\M) \to \C(\hat I)$, and by continuity, $\ntr : \AA_0 \to \C_0(I)$. Therefore $\ntr : \AA \to \C(\hat I)$. 

By Lem.~\ref{field}, $\ntr$ determines a map $\ntr_\infty : \A_\infty \to \co$ such that, by Lem.~\ref{trace.limit}, the composition $\ntr_\infty\circ\Po_\infty\circ T : \C(\M)\to\co$ is
\[
\ntr_\infty\circ\Po_\infty\circ T(f) = \frac1{\vol\M}\intM f \frac{\omega^n}{n!}
\mbox.\]
Because the \Kahler\ volume form is strictly positive, this is a faithful state on $\C(\M)$. Since this state factors through $\A_\infty$,
\[
\Po_\infty\circ T : \C(\M) \isom \A_\infty
\]
must be an isomorphism.

We have now verified the full hypotheses of Lem.~\ref{field} with the algebra $\AA$ and the linear map $\ntr : \AA \to \C(\hat I)$, therefore $\A$ is a continuous field of \cs-algebras with $\Gamma(\hat I,\A) = \AA$.
Since $\A$ is a continuous field over $\hat I = [s_0,\infty]$ for \emph{any} $s_0>\norm\Khat$, it is a continuous field over $(\norm\Khat,\infty]$.
\end{proof}

Let $\I := [0,\norm\Khat^{-1})$. In order to connect with the notation in the definition of strict quantization in Section \ref{quantization.definition}, we need to identify the parameter $s$ with $\hbar^{-1}$. Obviously, the interval $(\norm\Khat,\infty]$ is homeomorphic to $\I$ by this identification. 

\begin{thm}
\label{quantization.thm}
    $T:\C(\M)\to\Gamma(\hat I,\A)$ satisfies the definition of a strict quantization (p.~\pageref{quantization.def}) of $\M$ if we identify $\hbar=s^{-1}$.
\end{thm}
\begin{proof}
I have just proven that $\A$ is a continuous field. The first condition, that $\Im T_s$ generates $\As$, is true by construction. 
The third condition is given by Cor.~\ref{T.commutator}, eq.~\eqref{commutator.eq}. In fact, for $f,g\in\C^\infty(\M)$,
\[
\Norm{[T_s(f),T_s(g)]_- -is^{-1}T_s(\{f,g\})} = \Or^{-2}(s)
\mbox.\]
So, with $\hbar=s^{-1}\in\C_0\left(I\right)$,
\[
[T(f),T(g)]_- -i\hbar T(\{f,g\})
\]
is a section of $\A$ that is goes to $0$ as $\hbar\to0$ at least as fast as $\hbar^2$.  To prove quantization we merely needed this to go to $0$ faster than $\hbar$. So, this is more than enough.
\end{proof}

We know from Thm.~\ref{Morita2}, that the Toeplitz algebras are $\As=\K_\Bs(\Hs)$ for $s$ sufficiently large. Lemma \ref{big.Toeplitz} implies that the continuous field structures coincide. That is, 
\[
\Gz[I,\A] = \Gz[I,\K_\B(\Hi)]
\mbox.\]
In fact, using eq.~\eqref{C0.sections}, we have an extension,
\[
0 \to \K_{\Gz(I,\B)}(\Gz[I,\Hi]) \longrightarrow \Gamma(I,\A) \longrightarrow \C(\M) \to 0
\mbox.\]

\section{Surfaces}
\label{examples}
\subsection{The Torus}
I began by motivating this construction with the idea of generalizing the noncommutative torus. Let's return there and consider what happens if this construction is applied to a torus.

Let $\M$ be the torus obtained as the quotient of the complex plane, $\co$, by the lattice $\Z+i\Z$. The fundamental group is $\Z^2$ and the universal covering space is $\Mt=\co$. If we write the complex coordinate on $\co$ as $z=x+iy$, then the symplectic form is,
\[
\omega = 2\pi\,dx\wedge dy = i\pi\, dz\wedge d\bar z
\mbox.\]
This is normalized so that $\intM\frac{\omega}{2\pi} = 1$.

We can now follow through the computations in Section \ref{cocycle} explicitly.
First, we need a $1$-form $A$ such that $dA=\omega$. Let's use,
\[
A = \pi(x\,dy-y\,dx)
\mbox.\]
With this, we can define $L^s$ as the trivial line bundle over $\co$ with connection $\nabla_s = d-isA$. Let $L_0$ be the trivial flat line bundle over the torus with the trivial inner product multiplied by the canonical volume form, $dx\wedge dy$. So, $L_s:=  L^s\otimes \pi^*L_0$ is the $\co$-$\frac12$-density bundle given by the pair $(L^s,dx\wedge dy)$.  
For $(n,m)\in\Z^2$, we need $\phi_{(n,m)}$ to satisfy,
\[
d\phi_{(n,m)} = A - (n,m)^*A = -\pi(n\,dy-m\,dx) = \pi\, d(mx-ny)
\mbox.\]
So, we can take,
\[
\phi_{(n,m)} = \pi(mx-ny)
\mbox.\]
This defines the projective action of $\Z^2$ on sections of $L^s$ by eq.~\eqref{adef}.
From $\phi$, compute the cocycle $c$ on $(n,m),(n',m')\in\Z^2$ as,
\begin{align}
c[(n,m),(n',m')] &= \phi_{(n',m')} + (n',m')^*\phi_{(n,m)} - \phi_{(n+n',m+m')}
\nonumber \\
&= \pi(mn'-nm') 
\label{torus.cocycle}
\mbox.\end{align}
As expected, $x$ and $y$ cancel from this expression. 

\begin{lem}
In this case, the coefficient algebra is $\Bs=\T_s$, a noncommutative torus algebra.
\end{lem} 
\begin{proof}
Because $\Z^2$ is amenable, the reduced and maximal constructions coincide. The coefficient algebra is the twisted group \cs-algebra $\Bs = \cs(\Z^2,\ss)$ where $\ss=e^{isc}$. Let $u=(1,0)$ and $v=(0,1)$ be the generators of $\Z^2$. The algebra $\Bs$ is generated by the unitaries $[u]$ and $[v]$. Their products are,
\[
[u][v] = e^{isc(u,v)}[(1,1)] = e^{-i\pi s}[(1,1)]
\]
and
\[
[v][u] = e^{isc(v,u)}[(1,1)] = e^{i\pi s}[(1,1)]
\mbox.\]
Thus $[v][u]=e^{2\pi i s}[u][v]$. This relation defines the noncommutative torus algebra, $\T_s$.
\end{proof}

Note that because $\M$ and $L_0$ are flat, the curvature term in eq.~\eqref{Dsquared} is $\Khat=0$. So, $s>\norm\Khat=0$ is large enough for most of the results in earlier sections.

\begin{thm}
For $s>0$, $\Hs$ is a Morita equivalence of $\As$ and $\Bs$.
\end{thm}
\begin{proof}
We first need to check that $\Hs\neq 0$. The function, $\psi = e^{-\frac{\pi s}2(x^2+y^2)}$ is a holomorphic section of $L_s$ with the connection $\nabla_s$. Using eq.~\eqref{inner}, the inner product of this with itself is,
\[
\langle\psi\vert\psi\rangle = \frac1s \sum_{n,m=-\infty}^\infty e^{-\frac{\pi s}2 (n^2+m^2)} [(n,m)]
\mbox.\]
Since the coefficients fall off faster than exponentially, this sum is convergent in $\T_s$. Therefore $\psi$ is square integrable in the Hilbert $\T_s$-module sense. So, $\psi\in\Hs\neq 0$.

The algebras $\Bs=\T_s$ are very well studied. For $s$ irrational, $\Bs$ is a simple \cs-algebra and $\Hs$ is automatically full, since it is not $0$. If $s=n/m$ is a fraction in simplest form, then $\Bs$ is isomorphic to the algebra of $m\times m$-matrix-valued functions on a $2$-torus.  By Lem.~\ref{compact}, $\Hs$ is finitely generated and as a Hilbert $\Bs$-module it is projective. As such, $\Hs$ must be the space of continuous sections of a bundle over the $2$-torus. Referring to Lem.~\ref{Morita1}, any surjection $p:\Bs\onto B'$, must be the restriction of $m\times m$-matrix-valued functions to some closed subset of the torus. The push-forward of any bundle will still be a nonzero bundle, therefore $p_*\Hs$ is a nonzero $B'$-module. Therefore $\Hs$ is a full Hilbert $\Bs$-module.

With Thm.~\ref{Morita2}, fullness implies Morita equivalence for $s$ sufficiently large. However, because the torus and $L_0$ are flat, $\Khat=0$, and the blow-up construction used in Lem.~\ref{ample} is unnecessary for a surface, so any $s>0$ is sufficiently large.
\end{proof}

\begin{thm}
\label{torus}
For $s>0$, the Toeplitz algebra is a noncommutative torus algebra: $\As=\T_{1/s}$.
\end{thm}
\begin{proof}
It will be convenient in this proof to take advantage of Thm.~\ref{standard} 
and identify $\As$ with its isomorphic image in the representation on $\Hs^\Mt = \Hs\otimes_{\tau_s}\co = L^2_{\mathrm{hol}}(\Mt,L_s)$.

The cocycle $c$ extends trivially (by the same formula \eqref{torus.cocycle}) to the group $\R^2\supset\Z^2$. There is a projective action of $\R^2$ on sections of $L_s$ which is compatible with the connection for the same reasons that the $\Z^2$ action is. This can also be viewed as a representation of the Heisenberg group, $G$, the central extension 
\[
0\to \R \to G \to \R^2\to 0
\]
determined by $c$. This is a connected Lie group and we can work with its Lie algebra generators. The Heisenberg Lie algebra is generated by $X$, $Y$, and $\zeta$ with the relations,
\[
[X,Y] = -2\pi i\zeta,\quad [\zeta,X]=[\zeta,Y]=0
\mbox.\]
The unitary representation (generators are self-adjoint) on $\Hs^\Mt$ is given by,
\[
X = i \tfrac{\partial}{\partial x} +\pi s y,
\quad
Y = i \tfrac{\partial}{\partial y} - \pi s x,
\quad
\zeta = s
\mbox.\]
To minimize notation, I am using the same symbols for the generators and their representatives on $\Hs^\Mt$.
% The generators of the \emph{right} representation of $\Bs$ can be constructed from these;  $[u]$ is represented by $e^{iX}$ and $[v]$ by $e^{iY}$. The relation between $[u]$ and $[v]$ follows from this, but with the order reversed because of the right representation.

Because the generator $\zeta$ is central, it is represented by a constant in any irreducible representation. For any real $s\neq 0$, there is a unique irreducible unitary representation of $G$ with $\zeta$ represented by $s$. So, the representation on $\Hs^\Mt$ is either irreducible, or the direct sum of copies of a single irreducible representation.

Suppose that $\psi \in \Hs^\Mt$ satisfies $0=(X-iY)\psi$. This is the differential equation,
\[
0 = \left(i\tfrac\partial{\partial x} - \tfrac\partial{\partial y} + \pi s [ y+ix]\right) \psi
= 2i\left(\tfrac\partial{\partial z} + \tfrac\pi2 s \bar z\right) \psi
\mbox.\]
The condition that $\psi$ is holomorphic (as a section of $L_s$) is the differential equation,
\[
0 = \left(\tfrac\partial{\partial \bar z} + \tfrac\pi2 s z\right) \psi
\mbox.\]
The unique simultaneous solution of these two equations is,
\[
\psi = e^{-\frac\pi2 s z\bar z} = e^{-\frac\pi2 s (x^2+y^2)}
\mbox.\]
This is indeed square integrable for $s>0$. The existence of a unique solution implies that the representation of $G$ on $\Hs^\Mt$ is irreducible.

Because of this irreducibility, the image of $\csm(G)$ is $\K(\Hs^\Mt)$. So, any element of $\Li(\Hs^\Mt)$ is the weak limit of a sequence of linear combinations of (representatives of) elements of $G$. The set  
\[
\left\{e^{i(\alpha X+\beta Y)} \bigm| \alpha,\beta \in \R\right\}
\]
is thus a weak basis of $\Li(\Hs^\Mt)$. 

There is a representation of $G$ on $\C(\M)$, although it factors through the group $\T^2$. The representative of $X$ is $i\frac\partial{\partial x}$, the representative of $Y$ is $i\frac\partial{\partial y}$, and $\zeta$ is represented by $0$.
In this way, the Toeplitz map, $T_s^\Mt$, is equivariant. 

The commutative \cs-algebra, $\C(\M)$, of continuous functions on the torus is densely spanned by the functions $e^{-2\pi i (jx+ky)}$ for $j,k\in\Z$. Let $\psi = e^{-\frac{\pi s}2(x^2+y^2)} \in \Hs^\Mt$ again. We compute
\begin{align*}
\braket\psi{T_s^\Mt(e^{-2\pi i (jx+ky)})}\psi
&= \braket\psi{e^{-2\pi i (jx+ky)}}\psi \\
&= \int e^{-\pi s(x^2+y^2) - 2\pi i (jx+ky)} dx\, dy \\
&= s^{-1}e^{-\pi\, s^{-1}(j^2+k^2)}
\neq 0
\mbox,\end{align*}
and see that $T_s^\Mt(e^{-2\pi i (jx+ky)}) \neq 0$.

The Toeplitz operator $T_s^\Mt(e^{-2\pi i (jx+ky)})$ must be a multiple of the unique element of $\Li(\Hs^\Mt)$ which transforms in the same way as $e^{-2\pi i (jx+ky)}$ under $G$. We see that this is $e^{-i(kX-jY)/s}$ because,
\[
i\tfrac\partial{\partial x}e^{-2\pi i (jx+ky)} = 2\pi\, j e^{-2\pi i (jx+ky)}
\]
compares to
\[
[X,e^{-i(kX-jY)/s}]_- = i s^{-1}j [X,Y] e^{-i(kX-jY)/s}
= 2\pi\,  j e^{-i(kX-jY)/s}
\]
and
\[
i\tfrac{\partial}{\partial y}e^{-2\pi i (jx+ky)} = 2\pi k e^{-2\pi i (jx+ky)}
\]
compares to 
\[
[Y,e^{-i(kX-jY)/s}]_- = -i s^{-1}k [Y,X] e^{-i(kX-jY)/s}
= 2\pi k e^{-i(kX-jY)/s}
\mbox.\]

Define the elements $U := e^{iY/s}$ and $V := e^{-iX/s}$. These satisfy the noncommutative torus relation \eqref{nct} for $\theta=-1/s$ because 
\[
VUV^{-1}U^{-1}
= e^{[-iX/s,iY/s]_-} 
= e^{-s^{-2}[X,Y]} 
= e^{2\pi i/s}
\mbox.\]
The elements $U^jV^k$ and $e^{-i(kX-jY)/s}$ only differ by a phase. Therefore $T_s^\Mt(e^{-2\pi i (jx+ky)})$ is proportional to $U^jV^k$.

The subspace (algebraically) spanned by the functions $e^{-2\pi i (jx+ky)}$ for $j,k\in\Z$ is dense in $\C(\M)$, and the map $T_s^\Mt$ is norm-contracting. Therefore the image of this subspace is dense in $\Im T_s^\Mt$ and generates $\As$ (as a \cs-algebra). Therefore the set $\{U^jV^k \mid j,k\in\Z\}$ generates $\As$. Therefore $U$ and $V$ generate $\As$. Therefore $\As = \T_{1/s}$.
\end{proof}

If we make the identification $\hbar=1/s$, then $\As=\T_\hbar$ and $\Bs=\T_{1/\hbar}$.

Because $\Hs$ is essentially topological, and $\As=\K_\Bs(\Hs)$, the algebras $\As$ should be unchanged if we deform the geometry of the torus. So, for any $2$-torus and $s$ sufficiently large, the Toeplitz algebra $\As$ should be $\T_{1/s}$ as long as the symplectic form is normalized so that $\intM\omega = 2\pi$.
\medskip

Since this quantization construction gives such a different result from the standard Toeplitz construction, it is reasonable to ask how the two are related. Let $s=N\in\N$. It is not difficult to determine from the generators and relations that $\A_N = \T_{1/N}$ is isomorphic to the algebra of continuous $N\times N$ matrix-valued functions on a $2$-torus. However, it is also possible to ``explain'' this fact in terms of my construction. 

Because the central extension of $\Z^2$ is amenable, the maximal and reduced constructions coincide in this case. So, $\Lsm=\Ls$. Because $s=N$ is integral, $\ss = e^{iNc} = 1$ is the trivial cocycle. So, 
\[
\Bs = \Bsm \equiv \cs(\Z^2,\ss) = \cs(\Z^2) \cong \C(\T^2)
\mbox.\]
Any nonzero Hilbert $\C(\T^2)$-module is full, thus $\Hs$ is full. Because $\Bs$ is commutative, this implies that $\C(\T^2) \subset \K_{\Bs}(\Hs)$, and is central.

Let $L^N$ be a quantization line bundle over $\M$. By Thm.~\ref{universal1}, there exists a Hilbert $\Bs$-$\co$-bimodule, $R$, such that $L^N = \Ls\otimes_{\Bs}R$. However, $R$ is just a Hilbert space with a representation of $\Bs=\C(\T^2)$, and because $L^N$ is a line bundle, $R$ must be $1$-dimensional. This representation is just a character of $\Z^2$, given by a point of $\T^2$. This corresponds to the fact that quantization line bundles over the torus are not unique. A given one can be multiplied by a flat line bundle and these are classified by characters of the fundamental group $\Z^2$.

Lemma \ref{coefficients} shows that the Toeplitz algebra constructed with $L^N$ is a quotient $p : \As \onto \A_{L^N}$. The algebra $\A_{L^N}$ is just the full matrix algebra over the finite-dimensional Hilbert space $\Hi_{L^N}$ (see \cite{haw}). The dimension of $\Hi_{L^N} = \Gh(\M,L^N)$ is easily computed by the Riemann-Roch theorem to be $N$. So, $\A_{L^N}$ is the algebra of $N\times N$ matrices.

Since $\Hs$ is a finitely generated Hilbert $\C(\T^2)$-module, it is projective and is the space of continuous sections of a vector bundle over $\T^2$. Since $\Hs\otimes_{\Bs}R = \Hi_{L^N} \cong \co^N$, this is a rank $N$ vector bundle. This means that $\As=\K_\Bs(\Hs)$ is the algebra of continuous matrix-valued functions on this bundle, hence $\As=M_N[\C(\M)]$.

\subsection{Higher Genus Surfaces}
In \cite{k-l2,k-l4}, Klimek and Lesniewski considered a Riemann surface, $\M$, with constant negative curvature. They applied the Toeplitz quantization map for $\Mt$ to continuous functions on $\M$. By Thm.~\ref{standard}, this is equivalent to the construction I have presented here (for compact $\M$). 

Let $\M$ be a compact Riemann surface of genus $g\geq 2$ (I don't assume constant curvature). Because the universal covering space, $\Mt$, is contractible, $\M$ is a classifying space for its fundamental group $\G\equiv \pi_1(\M)$. The group cohomology can thus be computed as,
\[
H^*(\G) = H^*(\M)
\mbox.\]
The cohomology of interest here is in degree $2$, $H^2(\G;\Z) = H^2(\M,\Z) \cong \Z$. The group $\G$ has a universal central extension by $\Z$. This is the fundamental group of $\M$ with one point removed. This central extension is classified by the generator of $H^2(\G;\Z)$. 

I constructed the continuous field $\B$ of reduced twisted group \cs-algebras in Section \ref{continuous.field} by constructing an extension, $0\to\R\to\tilde\G\to\G\to0$, using the cocycle $c$. In this case $c$ must be proportional to the generator of $H^2(\G;\Z)$, so we can also use the universal central extension by $\Z$. The reduced \cs-algebra of this group is the algebra of continuous sections of a continuous field of \cs-algebras over a circle. These are the same twisted group \cs-algebras as before. This construction merely reflects the fact that $\Bs$ only depends upon $s$ modulo $1$. If this continuous field of \cs-algebras over $S^1$ is ``unrolled'', it becomes the same field $\B$ over $\R$ constructed before. Periodicity is the consequence of integrality.
\medskip

Natsume and Nest \cite{n-n} constructed a continuous field of ``noncommutative Riemann surfaces'' over $(2,\infty]$. Their starting point was the assumption that their algebra, $R_s$, should be Morita equivalent to the reduced twisted group \cs-algebra $\Bs$. The obvious question is of course whether the algebra $R_s$ is isomorphic to the Toeplitz algebra $\As$ given by my construction for a Riemann surface.  They construct $R_s$ as the ``full corner'' subalgebra of $\K(\Hs\otimes_{\tau_s}\co)\otimes\Bs$ given by a projection $e_s$; however, because of Morita equivalence, this can also be expressed as $R_s = \K_{\Bs}(H_{e_s})$, where $H_{e_s}$ is a finitely generated Hilbert $\Bs$-module. 

Recall from Section \ref{trace.index} that $[\tau_s]:K_0(\Bs)\to\R$ is the trace applied to $K$-theory. Some insight is gained here by applying the twisted $L^2$-index theorem.

\begin{prop}
\label{higher}
$[\tau_s][H_{e_s}] = [\tau_s][\Hs]$ and for $s$ irrational $[H_{e_s}]=[\Hs]$.
\end{prop}
\begin{proof}
It is first necessary to synchronize parameters. According to Prop.~3.3 of \cite{n-n}, Natsume and Nest's cocycle ``represents $e^{2\pi i s(1-g)}$ via the canonical
isomorphism of $H^2(\G;\T)$ with $\T$.'' By my geometrical construction of the twist cocycle from $\omega$, this means that 
\[
s(g-1) = \intM s\omega
\]
(and that $\M$ has mean curvature $-2$).

In the translation between different characterizations of Morita equivalence, $H_{e_s}$ is the image 
\[
e_s\cdot\left[(\Hs\otimes_{\tau_s}\co)\otimes\Bs\right]
\mbox.\]
 So, $[H_{e_s}]=[e_s]\in K_0(\Bs)$ and $[\tau_s][H_{e_s}] = \tau_s(e_s)$. This was computed in \cite{n-n} to be,
\[
\tau_s(e_s) = (s-1)(g-1)
\mbox.\]

The other side can be computed using the twisted $L^2$-index theorem discussed in Section \ref{trace}. Equation \eqref{L2.index} gives,
\begin{align*}
[\tau_s][\Hs] &= \intM \td(\TM)\wedge e^{s\omega/2\pi} \\
&= \intM \left(\frac{s\omega}{2\pi} - c_1\TM\right) \\
&= (g-1)s - (g-1) \\
&= (s-1)(g-1)
\end{align*}

The relevant $K$-theory group has been computed as $K_0(\Bs) \cong \Z^2$. The range of the trace on $K$-theory is $\Im[\tau_s] = \Z + s(g-1)\Z$. So, for $s$ irrational, $[\tau_s] : K_0(\Bs) \into \R$ is injective and the second claim follows from the first.
\end{proof}

It seems quite likely that this is not just a stable equivalence --- that indeed $\Hs \cong H_{e_s}$ for all $s$. If this is so, then by Thm.~\ref{Morita2}, $R_s\cong\As$. Proving this requires a more meticulous deconstruction of $e_s$. 

\section{Conclusions}
\label{conclusions}
My starting point here was the idea that a compact \Kahler\ manifold might be quantized for continuous values of the parameter $s=\hbar^{-1}$ if we use Toeplitz maps constructed for the universal covering space. What have we found?

This covering construction does give such a quantization, provided that the lift of the symplectic form to the universal covering space is exact. The Toeplitz algebra given by this construction at $s$ is $\As = \K_\Bs(\Hs)$ (Thm.~\ref{Morita2}), the algebra of compact operators on a Hilbert $\Bs$-module, where $\Bs=\csr(\G,\ss)$ is a reduced, twisted \cs-algebra of the fundamental group. The twist, $\ss$, is determined by the cohomology class of the symplectic form (Thm.~\ref{cohomology.isomorphism}). The Hilbert $\Bs$-module, $\Hs$, is the space of holomorphic sections (Lem.~\ref{holomorphic}) of a bundle whose fibers are Hilbert $\Bs$-modules. Modulo stable equivalence, $\Hs$ is determined by topology, not the specific geometry of the manifold, and its $K$-theory class can be constructed (eq.~\eqref{BC}) using a twisted version of the Baum-Connes assembly map. Thus, modulo the issue of stable equivalence, the algebra $\As$ can in principle be determined topologically. The continuous field structure of the collection of Toeplitz algebras is also mostly determined (Lem.~\ref{big.Toeplitz} and Thm.~\ref{A.field}) by a continuous field structure for the Hilbert \cs-modules $\Hs$, which is in turn related to an easily constructed continuous field structure of the reduced, twisted group \cs-algebras (Thm.~\ref{Bfield}).

When this construction is applied to a $2$-dimensional torus, it yields the ``noncommutative torus'' algebras (Thm.~\ref{torus}). This covering construction coincides with that of Klimek and Lesniewski in the case of Riemann surfaces. It also appears to reproduce the quantum surfaces of Natsume and Nest (Prop.~\ref{higher}), and ``explains'' their assumption that the algebras should be Morita equivalent to reduced, twisted group \cs-algebras of the fundamental group.

This covering construction motivates a generalized Toeplitz construction, based on a bundle, $Q$, of Hilbert $B$-modules with curvature $s\omega$. Any such ``quantization bundle'' can be constructed (Thm.~\ref{universal1}) from a unique maximal quantization bundle, $\Lsm$, of Hilbert $\cs(\G,\ss)$-modules. If we choose some point $x$ of the manifold, then there is a projective representation of the fundamental group, $\G$, on the fiber $R=Q_x$. Using this structure, the Toeplitz algebra for $Q$ can be constructed (Thm.~\ref{Morita2}) from the Hilbert $\cs(\G,\ss)$-module, $\Hsm$, constructed with $\Lsm$. Moreover, $\Hsm$ is itself determined modulo stable equivalence by topology and can potentially be computed topologically.
\medskip

This construction also seems enticingly pertinent to the (untwisted) Baum-Con\-nes conjecture. The reduced Baum-Connes map for the group $\G$ and the space $\Mt$ is a homomorphism,
\[
\mu_{\mathrm r} : K_*(\M) \to K_*[\csr(\G)]
\mbox.\]
On the left, the $K$-theory of $\M$ is of course the $K$-theory of the \cs-algebra $\C(\M) \cong\A_\infty$. On the right, the reduced group \cs-algebra is a special case  $\csr(\G)=\B_0$ of the  twisted group \cs-algebra. This Baum-Connes map can thus be thought of as
\[
\mu_{\mathrm r} : K^*(\A_\infty) \to K_*(\B_0)
\mbox.\]

Both of these \cs-algebras are fibers of the two continuous fields of \cs-algebras that I have constructed. Although $\A_0$ and $\B_\infty$ are not defined,  for large finite $s$ both fields are defined and the algebras are related by $\As=\K_\Bs(\Hs)$; if the Hilbert $\Bs$-module $\Hs$ is full, then these algebras are Morita equivalent.

This quantization construction thus provides an indirect connection between the left and right hand sides of the Baum-Connes map. In the case that the Hilbert \cs-modules are full, these \cs-algebras are (up to Morita equivalence) interpolated by a continuous field. This construction may therefore be useful for understanding the Baum-Connes conjecture for the fundamental groups of some compact \Kahler\ manifolds.

\section*{Acknowledgements}
This paper is essentially my thesis in mathematics at Pennsylvania State University; I would
thus like to thank my advisor Nigel Higson, and the other members of my committee: Paul Baum, Victor Nistor, and John Collins. I also wish to thank Ugo Bruzzo, Cesare Reina, and Marius Wodzicki for advice with parts of this work. This research was carried out at The Pennsylvania State University, la Scuola Internazionale Superiore di Studi Avanzati (Trieste), and l'Institut des Hautes \'{E}tudes Scientifiques (Paris).

%: Bibliography

\end{document}